\documentclass[12pt,a4paper]{article}
\usepackage{amsfonts}

\usepackage{graphicx,multicol}
\usepackage{placeins}
\usepackage[a4paper,margin=1in]{geometry}
\usepackage{amsmath, amsthm, amssymb}
\usepackage{setspace}
\usepackage{subcaption}
\usepackage[round,authoryear]{natbib}
\usepackage{placeins}
\usepackage{hyperref}

\theoremstyle{plain}

\newtheorem{theorem}{Theorem}

\newtheorem{lemma}{Lemma}
\newtheorem{corollary}{Corollary}
\theoremstyle{remark}
\newtheorem{remark}{Remark}
\theoremstyle{definition}

\newtheorem*{notation}{Notation}

\begin{document}

\title{Bias robustness of depth estimators in multivariate settings}
\author{Jorge G. Adrover\thanks{%
		Corresponding author, email: jorge.adrover@unc.edu.ar, ORCID: 0000-0003-1018-7895} \\
	FAMAF, National University of C\'ordoba, CIEM and CONICET, Argentina \and Marcelo Ruiz \thanks{%
		Email: mruiz@exa.unrc.edu.ar, ORCID: 0009-0008-1724-8949} \\
	Department of Mathematics, National University of R\'{\i}o Cuarto, Argentina \\
}
\date{}
\maketitle
\date{}
\maketitle
 
\begin{abstract}
The concept of statistical depth extends the notions of the median and 
quantiles to other statistical models. These procedures aim to formalize the
idea of identifying deeply embedded fits to a model that are less influenced
by contamination. In the multivariate case, Tukey's median was a 
groundbreaking concept for multivariate location estimation, and its 
counterpart for scatter matrices has recently attracted considerable 
interest. The breakdown point and the maximum asymptotic bias are key 
concepts used to summarize an estimator's behavior under contamination.  We 
explicitly obtain the maximum bias curve, contamination sensitivity and 
breakdown point of the deepest scatter matrices. In 
the multivariate and regression setting we analyse recently introduced 
error bounds that provide a unified framework for studying 
both the statistical convergence rate and robustness of Tukey's median, 
depth-based scatter matrices and multivariate regression 
estimators. We observe that slight variations in these inequalities allow us
to visualize the maximum bias behavior of the deepest estimators. We also point out that
all the halfspace depths under consideration can be obtained from a unifying concept called residual smallness depth. A numerical study is performed to  compare the finite sample bias performance of several robust estimators in  the multivariate setting. 
\end{abstract}



\section{Introduction}

The concept of asymptotic maximum bias for estimators in a contamination
neighborhood has been treated in the robust statistics literature since the pioneering
paper by \cite{Huber1964} on robust location estimation. The maximum asymptotic bias function $B(\varepsilon)$ provides a more accurate description of the estimator's behavior. It captures global behavior by quantifying how much the estimator deviates over the entire contamination neighborhood, under different levels of contamination $\varepsilon$, from the parameters at the central model. Two closely related
concepts to the maximum bias are contamination sensitivity and breakdown
point. The seminal paper by \cite{Huber1964} highlighted the median as an
estimator with optimal worst-case behavior under contamination since it
minimizes the maximum bias in the class of equivariant location estimators.

Among the many contributions on maximum bias, we mention the following:
\cite{hesimpson1993} obtained a lower bound on the contamination bias of an
estimator that holds for a wide class of parametric families and a second
lower bound which applies to locally linear estimates yielding that such
estimates cannot be bias minimax among all Fisher-consistent estimates in
higher dimensions; \cite{MZ1993A}
established expressions for the maximum asymptotic bias of M-estimates of
scale over the $\varepsilon $-contamination neighborhood and the
corresponding asymptotically minimax bias robust estimates of scale; \cite%
{GRZ2008} considered robust scale estimators in the nonparametric regression
setting and their maximum bias curves; \cite{MYZ1989} found minimax bias
estimators for two different classes of regression estimates: M-estimates
with general scale and GM-estimates in the context of known intercept and
elliptical regressors; \cite{MaronnaYohai1993} delved into a projection
estimator for regression whose maximum bias is twice the lower bound
attainable for maximum bias in regression; \cite{BZ2001} calculated the
maximum bias for several classes of estimators in a general setup of unknown
intercept and nonelliptical regressors; \cite{MaronnaYohai1995} studied the bias
behavior of the Stahel-Donoho estimator for scatter among a bundle of
multivariate dispersion estimators, \cite{A1998} computed the asymptotic
maximum bias for M-estimators for multivariate dispersion showing that Tyler's scatter M-estimator, \cite{Tyler1987}), minimizes the maximum bias in the $\varepsilon $%
-neighborhood in the case of known location; \cite{ChenTyler2002}
investigated various properties of Tukey's median, including its influence
function and maximum contamination bias; \cite{AMY2002} derived the maximum
asymptotic bias for the regression depth estimator and compared its
performance with that of \cite{MaronnaYohai1993} in the context of known
intercept and elliptical regressors; \cite{ZCY2004} computed the maximum
bias and influence functions for projection based estimators for
multivariate location; \cite{ZuoCui2005} derived and examined the influence
function and the maximum bias of the projection depth weighted scatter
estimators.

The concept of statistical depth has received considerable attention as a
way to extend the notions of the median and quantiles to more general
statistical models. Originally introduced by \cite{Tukey1975}, depth was
defined as the minimum proportion of data points lying on either side of a
point, and this idea was later extended to the bivariate case. \cite{DonohoGasko1992} formalized the general notion of the halfspace depth of a
point $\mathbf{z}$ in a $p$-dimensional space with respect to a probability
measure $P$ in which the search is for a location or center of the
probability measure $P$. A key feature of the median is that it is flanked
on each side by half of the data mass, helping to shield it from the effects
of outliers. This robustness---being well-surrounded by data---is echoed in
Tukey's concept of the multivariate median. In this case, a point in a
Euclidean space is considered a median if it maximizes the minimum mass
contained in any closed halfspace whose boundary includes the point. In this
sense, the search for a point or location center that fits deeply within the
data cloud is the idea that pervades some other statistical models.

This pursuit of a point that is "deeply inserted" into the data, to enhance
robustness in univariate and multivariate location models, has been extended
to other statistical frameworks. \cite{RousseeuwHubert1999} defined regression depth, which measures how deeply a
linear fit is embedded within the data. This is evaluated by the smallest
amount of data mass in the two opposing wedges formed by the fit plane and
vertical planes orthogonal to the explanatory subspace. \cite%
{MizeraMuller2004} further explored depth concepts in the context of
location-scale models. In fact, the univariate regression depth treated by 
\cite{RousseeuwHubert1999} is a special case of the multivariate regression
depth defined by \cite{Mizera2002}. A very similar but earlier definition
was proposed in \cite{BernEppstein2000}. \cite{BaiHe1999} derived the
asymptotic distribution of the maximal depth regression estimator and the Tukey's deepest point,
whose limiting distribution is characterized through a max--min operation of a
continuous process. \cite{Nagy2019} discussed the relationships
between halfspace depth and problems in affine and convex geometry, offering
an extensive overview of various depth notions.

\cite{ChenGaoRen2018} broke new ground by introducing a unified way to study
the statistical convergence rate and robustness jointly in order to come up
with multivariate estimators for location and scatter which achieve the
minimax rate $p/n+\varepsilon ^{2}$, with $p$ the dimension and $n$ the
sample size. To achieve this goal, they extended the idea of depth to covariance
matrix estimation by introducing the concept of matrix depth. Their
estimator achieves the optimal rate under Huber's $\varepsilon$-contamination model for estimating covariance/scatter matrices with various
structural assumptions, such as bandedness and sparsity. \cite%
{PaindaveineVanBever2018} also developed halfspace depth concepts for
scatter, concentration and shape matrices. While their concept of scatter
depth coincides with that of \cite{ChenGaoRen2018} and \cite{Zhang2002},
rather than focusing on the deepest scatter matrix, they studied the
properties of the depth function and its associated depth regions. \cite%
{Louvet2024} obtained the influence function and sensitivity curve for the
scatter halfspace depth. \cite{Gao2020} also considered estimators that are
maximizers of multivariate regression depth functions, by studying its
minimax rates in the settings of $\varepsilon $-contamination models for
various regression problems including nonparametric regression, sparse
linear regression, reduced rank regression, etc. 
We provide further insight into the error bounds obtained in \cite{ChenGaoRen2018}, showing that they can be derived in a more transparent manner by examining the effect of the estimator’s maximum bias under an $\varepsilon$-contamination neighborhood.
It is not surprising that an error bound over the $\varepsilon$-contamination neighborhood involves both the convergence rate $\sqrt{p/n}$, which captures stochastic variability, and the maximum bias $B(\varepsilon)$, which characterizes the estimator’s asymptotic behavior under contamination.
The maximum bias $B(\varepsilon)$ diverges to infinity as $\varepsilon$ approaches the breakdown point $b^*$; therefore, it cannot be of order $\varepsilon$ on any interval $(b, b^{\ast })$, with $0 \le b < b^{\ast }$. Since depth-based estimators have bounded contamination sensitivity, $B(\varepsilon)$ is of order $\varepsilon$ in a neighborhood of $0$. However, it is well known in the regression setting that there are many examples of robust estimators for which $B(\varepsilon)$ is of order $\sqrt{\varepsilon}$ near $0$; see \cite{He1991} and \cite{YZ1997}.

Section 2 reviews the concept of halfspace depths in several settings:
multivariate location and scatter, regression and multivariate regression
models. Section 3 introduces the concept of maximum asymptotic bias in the
different settings treated in Section 2. In Section 4 the maximum bias
function is derived for the deepest scatter matrix estimator introduced by 
\cite{ChenGaoRen2018} in the case of known location, showing that this
estimator shares Tukey's median's asymptotic breakdown point of $1/3$.
Section 5 examines the concentration inequalities in \cite{ChenGaoRen2018},
by showing that slight variations in the derivation of these inequalities
allow us to visualize the role of the maximum bias in the behavior of the
deepest estimators and its effect in the error bounds. Section 6 brings the
attention to a unified view of halfspace depths through a residual smallness
concept. Section 7 deals with a Monte Carlo simulation study comparing
several robust proposals under contamination when the location estimator is
included. The impact of including the estimation of the multivariate
location is typically avoided in maxbias derivations for scatter matrices
because performing a theoretical analysis becomes intractable. The proofs of
the results are deferred to the Appendix.

\section{Depth in several statistical models: a review}
Let us consider the multivariate location and scatter model (MLSM), 
$\mathbf{X}=\boldsymbol{\mu}_{0}+V_{0}\mathbf{u}$  
with $\boldsymbol{\mu}_{0}\in \mathbb{R}^{p},$ $V_{0}\in 
\mathbb{R}^{p\times p}$ an invertible matrix and $\mathbf{u}\sim P_{0}$ is an
unobservable random vector. $\boldsymbol{\mu }_{0}$\textbf{\ }and\textbf{\ }$V_{0}V_{0}^{t}$ (except for a constant) are assumed to be identifiable (for
instance this holds if $P_{0}$ is a centrosymmetric distribution around $\mathbf{0}$). Tukey's depth of a vector $\boldsymbol{\theta }\in 
\mathbb{R}^{p}$ is given by%
\begin{equation}
D_{T}\left( \boldsymbol{\theta },P\right) =\inf_{\mathbf{u\in }\mathcal{S}
^{p-1}}P\left( \mathbf{u}^{t}\mathbf{X}\leq \mathbf{u}^{t}\boldsymbol{\theta }
\right) ,  \label{tukey-depth}
\end{equation}%
with $\mathcal{S}^{p-1}{\ =}\left \{ \mathbf{z\in }%
\mathbb{R}^{p}:\left \Vert \mathbf{z}\right \Vert =1\right \}$.
Then, $\boldsymbol{\hat{\theta 
}}\left( P\right) =\arg \max_{\boldsymbol{\theta }\in 
\mathbb{R}^{p}}D_{T}\left(\boldsymbol{\theta}, P\right) $ is taken to be the deepest
estimator.

In the setting of multivariate linear regression, consider the model $Y=B^{t}%
\mathbf{X}+\sigma Z,$ where $B\in\mathbb{R}^{p\times m}$ and $\sigma >0$ are unknown, $Y\in 
\mathbb{R}^{m}$ is the random response vector and $\mathbf{X}\in 
\mathbb{R}^{p}$ is the vector of covariates with $Z\in 
\mathbb{R}^{m}$ an unobservable random vector independent of $\mathbf{X}$. We assume
that $B$ and $\sigma $ are identifiable (for instance, if $Z\sim F_{0}$ is
centrosymmetric around $\mathbf{0}$). The multivariate regression depth of $%
B\in\mathbb{R}^{p\times m}$ is defined as%
\begin{equation}
D_{MR}\left( B,P\right) =\inf_{U\in 
\mathbb{R}^{p\times m}-\left \{ 0\right \} }P\left( \left \langle U^{t}\mathbf{X}
,Y-B^{t}\mathbf{X}\right \rangle \geq 0\right) ,  \label{multregmizera}
\end{equation}%
with $\left \langle \cdot ,\cdot \right \rangle $ being the Euclidean scalar
product in $\mathbb{R}^{m}$. The deepest estimator is defined to be 
$\hat{B}\left( P\right) =\arg \max_{B\in 
\mathbb{R}^{p\times m}}D_{MR}\left( B,P\right)$.
The definition of multivariate regression depth was considered by \cite%
{BernEppstein2000} and \cite{Mizera2002}. When $m=1$, we have the
univariate regression depth, $D_{R}\left( \boldsymbol{\beta },P\right) =\inf_{\mathbf{u}\in 
\mathbb{R}^{p}-\left \{ 0\right \} }P\left( \left( \mathbf{u}^{t}\mathbf{X}\right)
\left( y-\boldsymbol{\beta }^{t}\mathbf{X}\right) \geq 0\right)$, see \cite{RousseeuwHubert1999}.

Interest in multivariate scatter depth has been revitalized over the past
decade, particularly following the seminal papers by \cite{ChenGaoRen2018}
and \cite{PaindaveineVanBever2018}, which explored different aspects of the
concept of depth for multivariate scatter. If we take the set $\mathcal{E}=\left\{ A\in 
\mathbb{R}^{p\times p}:A=A^{t}\text{ and }\mathbf{x}^{t}A\mathbf{x}>0\text{ for all } \mathbf{x\neq 0}\right\}$,  $\mathbf{X}\sim P$,  $\mathbf{X}\in 
\mathbb{R}^{p},$ the depth of \textit{$\Gamma $}$\in \mathcal{E}$ is taken to be 
\begin{equation*}
D_{S}\left( \mathit{\Gamma },P\right) =\inf_{\mathbf{u}\in \mathcal{S}%
^{p-1}}\min \left\{ 
\begin{array}{l}
P\left( \left\vert \mathbf{u}^{t}\left( \mathbf{X}-\mathbf{v}_{0}(P)\right)
\right\vert ^{2}\leq \mathbf{u}^{t}\mathit{\Gamma }\mathbf{u}\right) , \\
P\left( \left\vert \mathbf{u}^{t}\left( \mathbf{X}-\mathbf{v}_{0}(P)\right)
\right\vert ^{2}\geq \mathbf{u}^{t}\mathit{\Gamma }\mathbf{u}\right)%
\end{array}%
\right\}
\end{equation*}%
with $\mathbf{v}_{0}$ a preliminary affine equivariant location functional
used to yield an affine equivariant multivariate scatter functional. If
location is known, we can assume without loss of generality that $\mathbf{v}%
_{0}=\mathbf{0}$. For the known location case, the depth and the deepest
estimator are defined as, 
$D_{S}\left( \mathit{\Gamma },P\right) =\inf_{\mathbf{u}\in \mathcal{S}%
^{p-1}}\min \left\{ P\left( \left\vert \mathbf{u}^{t}\mathbf{X}\right\vert
^{2}\leq \mathbf{u}^{t}\mathit{\Gamma }\mathbf{u}\right) ,P\left( \left\vert 
\mathbf{u}^{t}\mathbf{X}\right\vert ^{2}\geq \mathbf{u}^{t}\mathit{\Gamma }%
\mathbf{u}\right) \right\}$,
$\mathit{\hat{\Gamma  }}\left( P\right) =\arg \underset{\mathit{\Gamma }\in 
\mathcal{E}}{\max }D_{S}\left( \mathit{\Gamma },P\right)$ and $D_M(P)=D_S\left( \hat{\mathit{\Gamma}}\left( P\right),P\right)$. 
Given a random sample either in the multivariate model $\left\{ \mathbf{X}%
_{i}\right\} _{i=1}^{n}$ or in the multivariate regression model $\left\{
\left( \mathbf{Y}_{i},\mathbf{X}_{i}\right) _{i=1}^{n}\right\} $ and the
corresponding empirical distribution function $P_{n}$ based on the sample,
then we set the depth estimators as $\hat{\boldsymbol{\theta}}_{n}=\arg \max_{%
\boldsymbol{\theta \in }%
\mathbb{R}
^{p}}D_{T}\left( \boldsymbol{\theta },P_{n}\right)$, $\hat{\Gamma}%
_{n}=\arg \underset{\mathit{\Gamma }\succeq \mathbf{0}}{\max }D_{S}\left( 
\mathit{\Gamma },P_{n}\right) $ and $\hat{B}_{n}=\arg \underset{B\in 
\mathbb{R}^{p\times m}}{\max }D_{MR}\left( B,P_{n}\right)$.

Our framework assumes that we have observations coming from Huber's $%
\varepsilon$-contamination neighborhood in which we have a majority of
observations coming from a parametric model and a minority coming from an
unknown distribution. Depth estimators are supposed to be much less affected
by the presence of spurious observations as it can be assessed by using
robustness measures. If $\Theta $ denotes the set in which the parameters
are assumed to lie in the central model, the breakdown point of an estimator
is a value that quantifies the level $\varepsilon $ of contamination
required to cause the estimator to move outside any compact subset of $%
\Theta $. The asymptotic breakdown point is an aspect of a more powerful
notion to measure the performance of an estimator, the asymptotic maximum
bias, whose concept is treated in the next section. In spite of the deep
understanding that the maxbias function provides, it is usually neglected
because of the technicalities that its derivation requires.

\section{Maximum bias in different statistical settings}\label{maxbiasdifstat}

\subsection{Maximum bias in the multivariate location and scatter model}

Given an elliptical distribution $\mathbf{X}\sim P_{0}^{E}\left( \cdot
\right) =P_{\mathbf{0}}\left( V_{0}^{-1}\left( \cdot -\boldsymbol{\mu }%
_{0}\right) \right) $ in  the MLSM, 
with $P_{0}$ a
centrosymmetric distribution around $\mathbf{0}$, the $\varepsilon $%
-contamination neighborhood for the multivariate model is given by $\mathcal{%
P}_{\varepsilon }\left( P_{0}^{E}\right) =\left\{ \left( 1-\varepsilon
\right) P_{0}^{E}\left( \cdot \right) +\varepsilon G\left( \cdot \right) ,%
\text{ }G\text{ any distribution on }%
\mathbb{R}^{p}\right\}$, with $\varepsilon \in \lbrack 0,1).$ Set \textit{$\Sigma $}$%
_{0}=V_{0}V_{0}^{t}$ and let $\mathcal{F}$ be a subset of distributions such
that $\mathcal{P}_{\varepsilon }\left( P_{0}^{E}\right) \subset \mathcal{F}$
and it also contains the empirical distribution functions. Observe that, if $\mathcal{L}\left( \mathbf{X}\right) \in \mathcal{P}_{\varepsilon }\left(
P_{0}^{E}\right) $ then $\mathcal{L}\left( \mathbf{X}\right) =\mathcal{L}%
\left( \mathit{\Sigma }_{0}^{1/2}\mathbf{\tilde{X}}+\boldsymbol{\mu }%
_{0}\right) $ with $\mathcal{L}\left( \mathbf{\tilde{X}}\right) \in 
\mathcal{P}_{\varepsilon }\left( P_{0}\right) .$ We say that the functionals 
$\boldsymbol{\ \hat{\mu}}:\mathcal{F}\rightarrow 
\mathbb{R}^{p}$ and $\mathit{\hat{\Gamma  } }:\mathcal{F}\rightarrow \mathcal{E}$ 
are \emph{affine and translation
equivariant} if and only if, for any invertible matrix $A\in 
\mathbb{R}
^{p\times p}$, it holds that $\boldsymbol{\hat{\mu}}\left( \mathcal{L}\left( A%
\mathbf{X+b}\right) \right) =A\boldsymbol{\hat{\mu}}\left( \mathcal{L}\left( 
\mathbf{X}\right) \right) \mathbf{\ +b}$ and $\mathit{\hat{\Gamma  } }\left( 
\mathcal{L}\left( A\mathbf{X+b}\right) \right) =A\mathit{\hat{\Gamma  } }\left( 
\mathcal{L}\left( \mathbf{X}\right) \right) A^{t}$. We say that the
functionals $\boldsymbol{\hat{\mu}}$ and $\mathit{\hat{\Gamma  } }$ are \emph{Fisher
consistent} if and only if $\boldsymbol{\hat{\mu}}\left( P_{0}^{E}\right) =%
\boldsymbol{\mu }_{0}$ and $\mathit{\hat{\Gamma  } }\left( P_{0}^{E}\right) =c\left(
V_{0}V_{0}^{t}\right) ,$ $c>0$ for any $\boldsymbol{\mu }_{0},V_{0}\in 
\mathbb{R}^{p\times p}$ invertible. The effect of the distortion caused by having $%
P\in \mathcal{P}_{\varepsilon }\left( P_{0}^{E}\right) $ can be
measured in the following invariant manner, 
\begin{eqnarray}
b_{L}\left( \boldsymbol{\hat{\mu}},\varepsilon ,P\right) &=&\left( \boldsymbol{\hat{\mu%
}}\left( P\right) -\boldsymbol{\hat{\mu}}\left( P_{0}^{E}\right) \right)
^{t}\left( \mathit{\hat{\Gamma  } }\left( P_{0}^{E}\right) \right) ^{-1}\left( 
\boldsymbol{\hat{\mu}}\left( P\right) -\boldsymbol{\hat{\mu}}\left(
P_{0}^{E}\right) \right) \notag \\
b_{E}\left( \mathit{\hat{\Gamma  } },\varepsilon ,P\right)  &=&\sup_{\mathbf{u}\in 
\mathcal{S}^{p-1}}\frac{\mathbf{u}^{t}\mathit{\hat{\Gamma  } }\left( P\right) 
\mathbf{u}}{\mathbf{u}^{t}\mathit{\hat{\Gamma  } }\left( P_{0}\right) \mathbf{u}}
\label{biasexplo} \\ &=&\lambda _{(1)}\left( \mathit{\hat{\Gamma  } }\left( P_{0}\right) ^{-1/2}\mathit{%
\hat{\Gamma  } }\left( P\right) \mathit{\hat{\Gamma  } }\left( P_{0}\right) ^{-1/2}\right)  
=\lambda _{(1)}\left( \mathit{\hat{\Gamma  } }\left( P\right) \mathit{\hat{\Gamma  } }%
\left( P_{0}\right) ^{-1}\right) ,  \notag 
\end{eqnarray}
\begin{eqnarray}
b_{I}\left( \mathit{\hat{\Gamma  } },\varepsilon ,P\right) &=&\inf_{\mathbf{u}\in 
\mathcal{S}^{p-1}}\frac{\mathbf{u}^{t}\mathit{\hat{\Gamma  } }\left( P\right) 
\mathbf{u}}{\mathbf{u}^{t}\mathit{\hat{\Gamma  } }\left( P_{0}\right) \mathbf{u}}
\label{biasimplo} \\
&=&\lambda _{(p)}\left( \mathit{\hat{\Gamma  } }\left( P_{0}\right) ^{-1/2}\mathit{%
\hat{\Gamma  } }\left( P\right) \mathit{\hat{\Gamma  } }\left( P_{0}\right) ^{-1/2}\right)
=\lambda _{(p)}\left( \mathit{\hat{\Gamma  } }\left( P\right) \mathit{\hat{\Gamma  } }\left(
P_{0}\right) ^{-1}\right)  \notag \\
&=&\sup_{\mathbf{u}\in \mathcal{S}^{p-1}}\frac{\mathbf{u}^{t}\mathit{\hat{\Gamma  } }%
^{-1}\left( P\right) \mathbf{u}}{\mathbf{u}^{t}\mathit{\hat{\Gamma  } }^{-1}\left(
P_{0}\right) \mathbf{u}}=\lambda _{(1)}\left( \mathit{\hat{\Gamma  } }\left(
P_{0}\right) ^{1/2}\mathit{\hat{\Gamma  } }^{-1}\left( P\right) \mathit{\hat{\Gamma  } }%
\left( P_{0}\right) ^{1/2}\right) .  \notag
\end{eqnarray}%
The second and third equalities in (\ref{biasexplo}) and (\ref{biasimplo})
follow from standard arguments in multivariate analysis; see for instance
Section A7, p. 523 of \cite{Seber1984}. (\ref{biasexplo}) refers to the
\textquotedblleft explosion\textquotedblright\ behavior of the functional
over the neighborhood by comparing the two quadratic forms based on the
functional under the true probability $P$ and $P_{0},$ respectively.
Similarly,\ (\ref{biasimplo}) tries to display the "implosion" behavior of
the functional over the neighborhood. Then we can define the asymptotic
maximum biases for the location and scatter functionals as%
\begin{eqnarray*}
B_{L}\left( \boldsymbol{\hat{\mu}},\varepsilon ,P_{0}^{E}\right) &=&\sup_{P\in 
\mathcal{P}_{\varepsilon }}b_{L}\left( \boldsymbol{\hat{\mu}},\varepsilon
,P\right) , \\
B_{E}\left( \mathit{\hat{\Gamma  } },\varepsilon ,P_{0}^{E}\right) &=&\sup_{P\in 
\mathcal{P}_{\varepsilon }}b_{E}\left( \mathit{\hat{\Gamma  } },\varepsilon
,P\right) ,\text{ }B_{I}\left( \mathit{\hat{\Gamma  } },\varepsilon
,P_{0}^{E}\right) =\sup_{P\in \mathcal{P}_{\varepsilon }}b_{I}\left( \mathit{%
\hat{\Gamma  } },\varepsilon ,P\right),\\
B\left( \mathit{\hat{\Gamma  } },\varepsilon ,P_{0}^{E}\right) &=&\max \left\{
B_{E}\left( \mathit{\hat{\Gamma  } },\varepsilon ,P_{0}^{E}\right) ,B_{I}\left( 
\mathit{\hat{\Gamma  } },\varepsilon ,P_{0}^{E}\right) \right\} .
\end{eqnarray*}%
We say that the asymptotic explosion and implosion breakdown points are
given by 
 $\varepsilon _{L}^{\ast } =\inf \left\{ \varepsilon >0:B_{L}\left( \mathit{%
\hat{\Gamma  } },\varepsilon ,P_{0}^{E}\right) =\infty \right\}$, 
$\varepsilon _{E}^{\ast } =\inf \left\{ \varepsilon >0:B_{E}\left( \mathit{%
\hat{\Gamma  } },\varepsilon ,P_{0}^{E}\right) =\infty \right\}$, $\varepsilon
_{I}^{\ast }$ $=$ \newline $\inf \left\{ \varepsilon >0:B_{I}\left( \mathit{\hat{\Gamma  } }%
,\varepsilon ,P_{0}^{E}\right) =\infty \right\}$ and
$\varepsilon ^{\ast } =\min \left( \varepsilon _{E}^{\ast },\varepsilon
_{I}^{\ast }\right)$.

If we consider equivariant and Fisher consistent functionals $\boldsymbol{\hat{%
\mu}}$ and $\mathit{\hat{\Gamma  } }$ for location and scatter, it is easily proved
that 
$B_{L}\left( \boldsymbol{\hat{\mu}},\varepsilon ,P_{0}^{E}\right) = cB_{L}\left( 
\boldsymbol{\hat{\mu}},\varepsilon ,P_{0}\right)$, 
$B_{E}\left( \mathit{\hat{\Gamma  } },\varepsilon ,P_{0}^{E}\right) =B_{E}\left( 
\mathit{\hat{\Gamma  } },\varepsilon ,P_{0}\right)$ and 
$B_{I}\left( \mathit{\hat{\Gamma  } },\varepsilon ,P_{0}^{E}\right) =B_{I}\left( 
\mathit{\hat{\Gamma  } },\varepsilon ,P_{0}\right)$, 
which entails that the maximum bias can be computed using $\boldsymbol{\mu }_{0}=%
\mathbf{0}$ and $\mathit{\Sigma }_{0}=I.$

The \textit{contamination bias of a functional }$T$, see \cite%
{hesimpson1993}, is a local measure closely related to the maximum bias of $%
T $ at the central model $F,$ $B\left( T,\varepsilon ,F\right) .$ It is
defined as 
\begin{equation*}
\gamma \left( T,F\right) =\left. \frac{\partial B\left( T,\varepsilon
,F\right) }{\partial \varepsilon }\right\vert_{\varepsilon =0}  .
\end{equation*}%
For small $\varepsilon $, the maximum bias can be approximated by $%
B(T,\varepsilon ,F)\approx \varepsilon \gamma (T,F)$.

\subsection{Maximum bias in the multivariate regression model}

We will assume in the multivariate regression model (MRM) that the intercept is
known and the covariates have an elliptical distribution $G_{0}$ around $%
\mathbf{0}\in \mathbb{R}^{p}$ with finite second moments; see \cite{Gao2020}. Put $\pi _{\mathbf{x}}=\prod_{j=1}^{p}(-\infty ,x_{j}]$, $\mathbf{x}=\left(
x_{1},\dots ,x_{p}\right) $, then the joint cumulative distribution function is given by
$H_{B,\sigma }\left( \mathbf{y,x}\right) =E_{\mathbf{X}}\left\{ F_{0}\left( 
\frac{\mathbf{y}-B_{0}^{t}\mathbf{w}}{\sigma }\right) I_{\pi _{\mathbf{x}
}}\left( \mathbf{w}\right) \right\}$,  
and 
$\mathcal{P}_{\varepsilon }\left( H_{B,\sigma }\right) =\left \{ \left(
1-\varepsilon \right) H_{B,\sigma }\left( \mathbf{y,x}\right) +\varepsilon
G\left( \mathbf{y,x}\right) ,\text{ }G\text{ any distribution} 
\right \}$
is the $\varepsilon$-contamination neighborhood 
for the MRM with $\varepsilon \in \lbrack 0,1)$. Let $\mathcal{F}$ be a subset of
distributions such that $\mathcal{P}_{\varepsilon }\left( H_{B,\sigma
}\right) \subset \mathcal{F}$ $\ $\ for all $B\in 
\mathbb{R}
^{p\times m},\sigma >0$ and $\varepsilon \leq \varepsilon ^{\prime }$ and it
also contains the empirical distribution functions. Observe that, if $%
\mathcal{L}\left( \mathbf{Y,X}\right) \in \mathcal{P}_{\varepsilon }\left(
H_{B,\sigma }\right)$, then $\mathcal{L}\left( \mathbf{Y,X}\right) =\mathcal{%
\ L}\left( \mathbf{\tilde{Y}+}B^{t}\mathbf{X,X})\right) $ with $\mathcal{L}%
\left( \mathbf{\tilde{Y},X}\right) \in \mathcal{P}_{\varepsilon }\left(
H_{0,1}\right)$. We say that the functional $T:\mathcal{F}\rightarrow 
\mathbb{R}
^{p\times m}$ is affine, regression and scale equivariant if, for any
invertible matrix $A\in 
\mathbb{R}
^{p\times p},$ $C\in 
\mathbb{R}
^{p\times m}$ and $s\in 
\mathbb{R}
^{+},$ it holds that 
$T\left( \mathcal{L}\left( \mathbf{Y,}A\mathbf{X}\right) \right) =\left(
A^{t}\right) ^{-1}T\left( \mathcal{L}\left( \mathbf{Y,X}\right) \right)$,
$T\left( \mathcal{L}\left( \mathbf{Y+}C^{t}\mathbf{X,X}\right) \right)
=T\left( \mathcal{L}\left( \mathbf{Y,X}\right) \right) +C$ and
$T\left( \mathcal{L}\left( s\mathbf{Y,X}\right) \right) =sT\left( \mathcal{%
L }\left( \mathbf{Y,X}\right) \right)$
respectively. Therefore, we may define the bias of a functional $T:\mathcal{P%
}_{\varepsilon }\rightarrow 
\mathbb{R}^{p\times m}$ at a distribution $H$ as 
\begin{equation*}
b_{MR}\left( T,\varepsilon ,H\right) =\left\{ \frac{tr\left( \left( T\left(
H\right) -T\left( H_{B,\sigma }\right) \right) ^{t}S\left( G\right) \left(
T\left( H\right) -T\left( H_{B,\sigma }\right) \right) \right) }{\sigma ^{2}}
\right\}^{1/2},
\end{equation*}%
where $H_{B,\sigma }$ is the distribution under the central model, $S:%
\mathcal{P}_{\varepsilon }\rightarrow 
\mathbb{R}^{p\times p}$ is an affine equivariant estimator for the dispersion matrix
of $\mathbf{X}$. This definition in the multivariate regression model
coincides with that of the univariate regression model; see \cite{AMY2002}.
Then we can define the asymptotic maximum biases for the functionals as 
$B_{MR}\left( \hat{B},\varepsilon ,H_{B,\sigma }\right) =\sup_{H\in \mathcal{%
P }_{\varepsilon }}b_{MR}\left( \hat{B},\varepsilon ,H\right).$
It holds that $B_{MR}\left( \hat{B},\varepsilon ,H_{B,\sigma }\right) =B_{MR}\left( \hat{B}
,\varepsilon ,H_{0,1}\right)$ if we consider equivariant and Fisher consistent functionals $\hat{B}$,
which entails that the maximum bias can be computed using $\sigma
=1 $ and $B=0.$ If $m=1,$ the deepest estimator $\boldsymbol{\hat{\beta}}%
^{t}\in \mathbb{R}
^{p}$ and $B_{MR}\left( \boldsymbol{\hat{\beta}},\varepsilon ,H_{0,1}\right)
=\sup_{H\in \mathcal{P}_{\varepsilon }}\left \Vert \boldsymbol{\hat{\beta}}%
\left( H\right) \right \Vert $, the definition of maximum bias coincides
with the definition given in univariate regression; see, for instance, \cite%
{MYZ1989}, p.1610.

\section{Maximum bias for the deepest scatter matrix}

From now on we will assume that

A1. $P_{0}$ is a multivariate normal around $\mathbf{0}$ and covariance 
matrix $I$.

We will summarize the
steps to get the maximum bias function for the deepest scatter matrix. First
of all we compute the depth of any matrix under the normal model and the
deepest estimator under this distribution, which is given by $\left(\Phi
^{-1}\left( 3/4\right)\right)^2 I$. We next verify that in case of having a sequence
of contaminations in the $\varepsilon $-contamination neighborhood which
yields the deepest estimator to have either the largest eigenvalue going to
infinity or the smallest one going to $0,$ then the level of contamination $%
\varepsilon \geq 1/3,$ which entails that $\varepsilon ^{\ast }\geq 1/3.$
Next, we calculate the depth of any matrix in $\mathcal{E}$ 
under point mass contaminations located along the
direction given by a vector $\mathbf{e\in \mathbb{R}}^{p}$. 
Moreover we consider the set of the matrices whose eigenvector
associated with the largest eigenvalue coincides with the vector $\mathbf{e}$,%
which yields a more specific formula rather than the general one for any
matrix given in the previous step. Then, in this class, we calculate a
deepest scatter matrix (since uniqueness cannot be suspected at all) in
which the largest and smallest eigenvalues attainable for point mass
contaminations coincide with those of the bounds found in the concentration
inequalities. On the one hand we show that in case of having a sequence of
matrices such that the maximum eigenvalue tends to infinity its depth could
only converge to a value less than or equal to $\min \left( \varepsilon
,1-\varepsilon \right) .$ On the other hand, we prove that the deepest
estimator in the case of using point mass contaminations is less than or
equal to $\left( 1-\varepsilon \right) /2$. Therefore we can conclude that
if we have a level of contamination $\varepsilon $ such that the deepest
estimator has smallest and largest eigenvalues bounded above and below for
any contamination, the level of contamination must be less than or equal to $%
1/3$ and we get that the breakdown $\varepsilon ^{\ast }=1/3,$ which
coincides with that of Tukey's median. The reminiscence of Tukey's median is
even emphasized since we can use a similar reasoning to that of \cite%
{ChenTyler2002} to get a bound for the asymptotic maximum bias curve.

Let us state the notation and results briefly described in the previous
paragraph. We know that $\mathit{\hat{\Gamma  } }\left( P_{0}\right) $ is Fisher
consistent up to a constant, that is $\mathit{\hat{\Gamma  } }\left( P_{0}\right) =%
\left\{ \Phi ^{-1}\left( \frac{3}{4}\right)\right\} ^{2} \mathit{\Sigma }$.  Without loss of generality to study the
breakdown point and maximum bias, we can assume that $\mathit{\Sigma }=I.$
If $F$ is a centrosymmetric distribution in $\mathbb{R}
^{p}$ and $\mathbf{X}\sim F\,$with then $\mathbf{w}^{t}\mathbf{X}\sim 
\mathbf{z}^{t}\mathbf{X}$ for all $\mathbf{w,z\in \mathcal{S}}^{p-1}.$ Let
us take a symmetric positive matrix $\mathit{\Gamma },$ with eigenvalues $%
l_{1}\geq l_{2}\geq \dots $ $\geq l_{p}>0$ and eigenvectors $\mathbf{v}%
_{1},\dots ,\mathbf{v}_{p}$ respectively, $\mathit{\Gamma }%
=\sum_{j=1}^{p}l_{j}\mathbf{v}_{j}\mathbf{v}_{j}^{t}.$ 

Let $g : \emph{S}^{p-1}\rightarrow \left[ 0,1\right]$ be the function, %
$g\left( \mathbf{u}\right) =P_{0}\left( -\sqrt{\mathbf{u}^{t}\mathit{\Gamma 
}\mathbf{u}}\leq \mathbf{u}^{t}\mathbf{X}\leq \sqrt{\mathbf{u}^{t}\mathit{%
\Gamma }\mathbf{u}}\right)$.%
Then, we have that the depth under a centrosymmetric model of any positive
definite matrix is given by
\begin{lemma}
\label{depthmodel}$\ D\left( \mathit{\Gamma },P_{0}\right) =\min_{\mathbf{%
u\in }\emph{\ S}^{p-1}}\left( g\left( \mathbf{u}\right) ,1-g\left( \mathbf{u}%
\right) \right) =\min \left( g\left( \mathbf{v}_{p}\right) ,1-g\left( 
\mathbf{v}_{1}\right) \right) .$
\end{lemma}

Moreover, we can easily derive the deepest estimator, which confirms the
Fisher-consistency of the procedure, except for a constant.

\begin{corollary}
\label{coroldepth} $\hat{\Gamma  }(P_0)=\left(\Phi^{-1}(3/4)\right)^2 I$ and $D_S(\hat{\Gamma  }(P_0),P_0)=0.5$. 
\end{corollary}

The following result shows that if a sequence of contaminations
distributions yields a sequence of depth estimators exploding or imploding,
the level of contamination must be greater than $1/3$.

\begin{corollary}
\label{rateexplosion} Let $\left\{ P_{\varepsilon ,n}\right\} _{n=1}^{\infty
}$ be such that $P_{\varepsilon ,n}=\left( 1-\varepsilon \right)
P_{0}+\varepsilon P_{n}.$ If the associated depth estimators $\left\{ 
\mathit{\hat{\Gamma} }_{n}\right\} _{n=1}^{\infty }$, $\mathit{\hat{\Gamma} }%
_{n}=\sum_{j=1}^{p}l_{j}^{(n)}\mathbf{v}_{j}^{(n)} \left(\mathbf{v}_{j}^{(n)}\right)^t%
 $ have either their largest eigenvalue $l_{1}^{(n)}\rightarrow
\infty $ or their smallest eigenvalue $l_{p}^{(n)}\rightarrow 0$ then $%
\varepsilon >1/3$.
\end{corollary}

Given a unit vector $\mathbf{e}$ and $r>0,$ take the point mass
contaminations $\delta _{r}=\delta _{r\mathbf{e}}$. If we have $%
P_{\varepsilon ,r}=\left( 1-\varepsilon \right) P_{0}+\varepsilon \delta
_{r} $ and $\delta $ stands for the Kronecker delta (it is $1$ if the
inequality holds, $0$ otherwise), we consider the function $h^{\mathbf{e}%
}\left( \mathbf{v}\right) =\mathbf{v}^{t}\mathit{\Gamma }\mathbf{v/}\left( 
\mathbf{v}^{t}\mathbf{e}\right) ^{2}$ and $D\left( \mathit{\Gamma }%
,P_{\varepsilon ,r}\right) =\min_{\mathbf{u\in }\mathcal{S}^{p-1}}$ $%
m_{r}\left( \mathbf{u}\right) ,$ with 
$$
m_{r}\left( \mathbf{u}\right) =\min \left\{ \left( 1-\varepsilon \right)
g\left( \mathbf{u}\right) +\varepsilon \delta \left( h^{\mathbf{e}}\left( 
\mathbf{u}\right) \geq r^{2}\right) ,\left( 1-\varepsilon \right) \left(
1-g\left( \mathbf{u}\right) \right) +\varepsilon \delta \left( h^{\mathbf{e}%
}\left( \mathbf{u}\right) \leq r^{2}\right) \right\} .
$$

\begin{notation}
The function $h^{\mathbf{e}}\left( \mathbf{v}\right) $ plays a crucial role
in the derivation of the maximum bias. Take the sets 
$B_{r}^{\mathbf{e}} = \left \{ \mathbf{v}\in S^{p-1\text{ }}\text{and }h^{ 
\mathbf{e}}\left( \mathbf{v}\right) <r^{2}\right \}$,
$F_{r}^{\mathbf{e}} = \left \{ \mathbf{v}\in S^{p-1\text{ }}\text{and }h^{ 
\mathbf{e}}\left( \mathbf{v}\right) =r^{2}\right \}$ and
$A_{r}^{\mathbf{e}} =\left \{ \mathbf{v}\in S^{p-1\text{ }}\text{and }h^{ 
\mathbf{e}}\left( \mathbf{v}\right) >r^{2}\right \}$. 
Then we have the related quantities, 
\begin{equation*}
gb_{r}^{\mathbf{e}}=\inf_{B_{r}^{\mathbf{e}}}g\left( \mathbf{v}\right) , 
\text{ \ \ }Gb_{r}^{\mathbf{e}}=\sup_{B_{r}^{\mathbf{e}}}g\left( \mathbf{v}
\right) ,\text{ \ \ }Ga_{r}^{\mathbf{e}}=\sup_{A_{r}^{\mathbf{e}}}g\left( 
\mathbf{v}\right) ,\text{ \ \ }ga_{r}^{\mathbf{e}}=\inf_{A_{r}^{\mathbf{e}
}}g\left( \mathbf{v}\right) .
\end{equation*}
If either $B_{r}^{\mathbf{e}}=\emptyset $ or $A_{r}^{\mathbf{e}}=\emptyset $
put $gb_{r}^{\mathbf{e}}=\infty $ and $Gb_{r}^{\mathbf{e}}=-\infty $ or $%
ga_{r}^{\mathbf{e}}=\infty $ and $Ga_{r}^{\mathbf{e}}=-\infty .$
\end{notation}

The following lemma gives us the depth of any symmetric positive matrix $%
\mathit{\Gamma }$ under point mass contaminations.

\begin{lemma}
\label{depthpointmass} It holds that 
\begin{equation}
D\left( \mathit{\Gamma },P_{\varepsilon ,r}\right) =\min \left\{ 
\begin{array}{c}
\left( 1-\varepsilon \right) \left( 1-Gb_{r}^{\mathbf{e}}\right)
+\varepsilon ,\left( 1-\varepsilon \right) ga_{r}^{\mathbf{e}}+\varepsilon ,
\\ 
\left( 1-\varepsilon \right) gb_{r}^{\mathbf{e}},\left( 1-\varepsilon
)\left( 1-Ga_{r}^{\mathbf{e}}\right) \right)%
\end{array}%
\right\} .  \label{pointmassdepth}
\end{equation}
\end{lemma}

From this lemma we can derive the depth for a matrix whose eigenvector
associated with the largest eigenvalue coincides with the direction of the
point mass contamination.

\begin{corollary}
\label{deptheigencont} Take $\mathit{\Gamma }$ to be a matrix whose
eigenvector $\mathbf{v}_{1}$ coincides with\textbf{\ }$\mathbf{e.}$ If $%
r>l_{1}^{1/2}$ consider $\mathbf{v}_{m,r}=\arg \min_{B_{r}^{\mathbf{e}}\cup
F_{r}^{\mathbf{e}}}\mathbf{v}^{t}\mathit{\Gamma }\mathbf{v}$ and $\mathbf{v}%
_{M,r}=\arg \max_{A_{r}^{e}\cup F_{r}^{\mathbf{e}}}\mathbf{v}^{t}\mathit{%
\Gamma }\mathbf{v.}$ Then 
\begin{equation*}
D\left( \mathit{\Gamma },P_{\varepsilon ,r}\right) =\left\{ 
\begin{array}{cc}
\min \left\{ 
\left( 1-\varepsilon \right) g(\mathbf{v}_{p})+\varepsilon , 
\left( 1-\varepsilon \right) \left( 1-g\left( \mathbf{v}_{1}\right) \right)%
\right\} & \text{if }r\leq l_{1}^{1/2} \\ 
\min \left\{ 
\begin{array}{c}
\left( 1-\varepsilon \right) \left( 1-g\left( \mathbf{v}_{1}\right) \right)
+\varepsilon ,\left( 1-\varepsilon \right) g\left( \mathbf{v}_{p}\right)
+\varepsilon , \\ 
\left( 1-\varepsilon \right) g\left( \mathbf{v}_{m,r}\right) ,\left(
1-\varepsilon )\left( 1-g\left( \mathbf{v}_{M,r}\right) \right) \right)%
\end{array}%
\right\} & \text{if }r>l_{1}^{1/2}%
\end{array}%
\right. .
\end{equation*}
\end{corollary}

The following lemma restricted to matrices whose eigenvector associated with
the largest eigenvalue shares the same direction as the contamination gives
us which should be the deepest estimator in this class.

\begin{lemma}
\label{deepestrestricted} Let us take $\mathit{\Gamma }=\sum_{j=1}^{p}l_{j}\mathbf{%
v}_{j}\mathbf{v}_{j}^{t}$, $l_{1}\geq ...\geq l_{p},$ $\mathbf{v}_{i}^{t}%
\mathbf{v}_{j}=\delta _{ij}$ and $\mathbf{v}_{1}=\mathbf{e.}$ Given $\gamma =\left\{ \Phi ^{-1}\left( \frac{3-4\varepsilon }{4\left(
	1-\varepsilon \right) }\right) \right\}^{2}$, $l_{1}=\left\{ \Phi ^{-1}\left( 
\frac{3-\varepsilon }{4\left( 1-\varepsilon \right) }\right) \right\}
^{2},l_{p}=\left\{ \Phi ^{-1}\left( \frac{3-5\varepsilon }{4\left(
	1-\varepsilon \right) }\right) \right\}^{2},l_{p}=l_{p-1}=\cdots =l_{2}<\Phi
^{-1}\left( \frac{3}{4}\right)$, then the
deepest matrices in this class are given by $\sqrt{\gamma }I, \text{if }r\leq l_{1}^{-1/2}$, and
$\sum l_{j}\mathbf{v}_{j}\mathbf{v}_{j}^{t}, \text{if }r>l_{1}^{-1/2}$. 
Moreover, $D\left( \mathit{\Gamma },P_{\varepsilon ,r}\right) =\left( 1-\varepsilon \right) /2.$
\end{lemma}

\begin{remark} \label{comentario}
Let $\beta=\left\{ \Phi ^{-1}\left( \frac{3}{4}\right)\right\} ^{2}$. Then $%
\Phi \left( l_{1}^{1/2}\right) $ and $\Phi \left( l_{p}^{1/2}\right)$ will appear
in the error bounds given in Lemma \ref{chenscatter}. Furthermore, $B_{E}\left(
\varepsilon \right) = \frac{1}{\sqrt{\beta }}\Phi ^{-1}\left( \frac{%
3-\varepsilon }{4\left( 1-\varepsilon \right) }\right) -1 $ and $%
B_{I}\left( \varepsilon \right) =1-\frac{1}{\sqrt{\beta }} \Phi ^{-1}\left( 
\frac{3-5\varepsilon }{4\left( 1-\varepsilon \right) } \right)$.
\end{remark}

In case of having point mass contaminations going to infinity and taking
matrices whose largest eigenvalue goes to infinity as well with certain
rate, their depth should converge to the level of contamination as we can
establish in the following result.

\begin{lemma}
\label{behavior} Let us consider the case of having a family of matrices $%
\mathit{\Gamma }_{r}$ with the largest eigenvalue $l_{1}^{(r)}\rightarrow $ $%
\infty $ , $r^{2}l_{1}^{^{(r)}-1/2}\rightarrow 1/2,$ and the other
eigenvalues going to 0. Then, $\lim_{r\rightarrow \infty }D\left( \mathit{%
\Gamma },P_{\varepsilon ,r}\right) =\min \left( \varepsilon ,1-\varepsilon
\right) $.
\end{lemma}

The following lemma will give us a bound for the depth of the deepest
estimator under point mass contaminations.

\begin{lemma}
\label{depthpointmassgeneral} If $\mathit{\hat{\Gamma  } }$ stands for the deepest
estimator, then $D\left( \mathit{\hat{\Gamma  } },P_{\varepsilon ,r}\right) \leq
\left( 1-\varepsilon \right) /2$.
\end{lemma}

Finally we can show the asymptotic breakdown point for the deepest estimator
for multivariate scatter.

\begin{theorem}
\label{breakdownpoint}The asymptotic breakdown point of the deepest
estimator is $1/3$.
\end{theorem}


Next, we follow closely Theorems 4.1 and 4.2 of \cite{ChenTyler2002}. Let us
take the contaminated distribution $P_{\varepsilon ,Q}=\left( 1-\varepsilon
\right) P_{0}+\varepsilon Q$ and $A\subset \cal{E}$, 
given a probability $P$ define $L\left(\eta ,P\right) =\left\{ \mathit{\Gamma }\in {\cal E}:D\left( 
\mathit{\Gamma },P\right) \geq D_{M}\left( P\right) -\eta \right\}$, $\Lambda \left( \varepsilon ,P_{0}\right) =\inf_{Q}D_{M}\left(
P_{\varepsilon ,Q}\right)$, $\delta \left( \varepsilon ,P_{0}\right) =\frac{\Lambda \left( \varepsilon
	,P_{0}\right) -\left( 1-\varepsilon \right) D_{M}\left( P_{0}\right) }{%
	1-\varepsilon }$, $M\left( P\right) =\left\{ \mathit{\Gamma }\in {\cal E}:D\left( \mathit{\Gamma 
	},P\right) =D_{M}\left( P\right) \right\} =\bigcap_{0<\eta
	<D_{M}\left( P\right) }L\left( \eta ,P\right)$ and 
$\left\Vert A\right\Vert =\sup_{\mathit{\Gamma }\in A}\left\{ \frac{%
\left\Vert \mathit{\Gamma }\right\Vert _{op}}{\left\Vert \beta I\right\Vert
_{op}},\left\Vert \beta I\right\Vert _{op}\left\Vert \mathit{\Gamma }%
^{-1}\right\Vert _{op}\right\}$%

\begin{lemma}
\label{maximumbias} Let $\varepsilon <1/3$ and $P_{\varepsilon ,Q}=\left(
1-\varepsilon \right) P_{0}+\varepsilon Q$. It holds that

\begin{enumerate}
	
\item[(i)] $\Lambda \left( \varepsilon ,P_{0}\right) \geq \left( 1-\varepsilon
\right) D_{M}(P_{0})$ and $D_{M}\left( P_{\varepsilon ,Q}\right) \leq \left(
1-\varepsilon \right) D_{M}\left( P_{0}\right) +\varepsilon .$

\item[(ii)] Set $\alpha =\frac{e}{1-\varepsilon }-\delta \left( \varepsilon
,P_{0}\right)$. If $\mathit{\Gamma }\notin L\left( \alpha ,P_0\right) $ then $%
D\left( \mathit{\Gamma },P_{\varepsilon ,Q}\right) <\Lambda \left(
\varepsilon ,P_{0}\right) $ and $\mathit{\Gamma }$ cannot be a deepest
estimator.

\item[(iii)] $B\left( \mathit{\Gamma },\varepsilon ,P_{0}\right) \leq \left\Vert
L\left( \frac{\varepsilon }{1-\varepsilon },P_{0}\right) \right\Vert =\max
\left\{ \frac{1}{\sqrt{\beta }}\Phi ^{-1}\left( \frac{3-\varepsilon }{%
4\left( 1-\varepsilon \right) }\right) ,\frac{\sqrt{\beta }}{\Phi
^{-1}\left( \frac{3-5\varepsilon }{4\left( 1-\varepsilon \right) }\right) }%
\right\} $.
\end{enumerate}
\end{lemma}

Thus, we have the following result.

\begin{theorem}
$B\left( \mathit{\Gamma },\varepsilon ,P_{0}\right) =\max \left\{ \frac{1}{%
\sqrt{\beta }}\Phi ^{-1}\left( \frac{3-\varepsilon }{4\left( 1-\varepsilon
\right) }\right) ,\frac{\sqrt{\beta }}{\Phi ^{-1}\left( \frac{3-5\varepsilon 
}{4\left( 1-\varepsilon \right) }\right) }\right\} .$
\end{theorem}

\begin{remark}\label{comentariosegundo}
The implosion bias actually rules the bias since $B\left( \mathit{\Gamma }%
,\varepsilon ,P_{0}\right) =\sqrt{\beta }/\Phi ^{-1}\left( \frac{%
3-5\varepsilon }{4\left( 1-\varepsilon \right) }\right) .$ A proof of this
fact is available at the Appendix. Truth to be said, this fact
is also observed in the simulation study for almost all the estimators under
consideration, since the plots for the empirical bias show that  point mass
contaminations at $K=0$ or $K=1$ seem to provoke the worst bias situation.
This entails that the contamination bias for the deepest
estimator is given by
%
$\gamma \left( \mathit{\Gamma },\varepsilon ,P_{0}\right) =1/\left(2\sqrt{%
	\beta }\varphi \left( \sqrt{\beta }\right) \right).$
\end{remark}

\section{Maximum bias and concentration inequalities for depth estimators}

\cite{ChenGaoRen2018} broke new ground by introducing a unified way to study
the statistical convergence rate and robustness jointly. 
Let us state some notation used throughout this section. Given $\delta \in
(0,1/2),$ put $\alpha =1-2\delta $. Let $\mathcal{P}_{\varepsilon }\left(
P_{0}^{E}\right) $ be the $\varepsilon $-contamination neighborhood with $%
P_{0}^{E}=N(\boldsymbol{\theta },\mathit{\Sigma })$. Set $\mathcal{F}\left(
M\right) $ as the set of symmetric and definite positive matrices $\mathit{%
\Sigma }$ such that the largest eigenvalue $\lambda _{1}(\mathit{\Sigma })$
is less than a constant $M>0.$ Take $\varepsilon^{\prime }<1/3$.
Theorem 2.1 of \cite{ChenGaoRen2018} derived
that, for $\varepsilon \in \left[ 0,\varepsilon ^{\prime }\right] ,$ 
and $\left( p+\log \left( 1/\delta \right)
\right) /n$ sufficiently small, there exists a constant $C>0$ (depending on $%
\varepsilon ^{\prime }$ but independent of $p,n,\varepsilon),$ such that 
\begin{equation}
\inf_{\boldsymbol{\theta },\mathit{\Sigma }\in \mathcal{F}\left( M\right)
,P\in \mathcal{P}_{\varepsilon }\left( P_{0}^{E}\right) }P\left( \left\Vert 
\boldsymbol{\hat{\theta}}_{n}-\boldsymbol{\theta }\right\Vert ^{2}\leq
C\left( \max \left\{ \frac{p}{n},\varepsilon ^{2}\right\} +\frac{\log \left(
1/\delta \right) }{n}\right) \right) \geq \alpha .  \label{tukeymedian}
\end{equation}

The constant $C$ in the error bound (\ref{tukeymedian}) is actually affected
by the asymptotic maximum bias of Tukey's median. \cite{ChenTyler2002}
derived the asymptotic maximum bias for Tukey's median $\boldsymbol{\hat{\theta}}%
_{T}$ for $p\geq 2,$ which turns out to be $B_{L}\left( \hat{\boldsymbol{%
		\theta }}_{T},\varepsilon ,\Phi \right) =\Phi ^{-1}\left( \frac{%
	1+\varepsilon }{2\left( 1-\varepsilon \right) }\right) $. The heuristics
behind expecting the asymptotic maximum bias function to appear in the
concentration inequality is as follows: As the sample size tends to
infinity, the estimator is expected to converge to the functional value $%
\hat{\boldsymbol{\theta }}_{T}\left( P\right) .$ The quantity $\left\Vert 
\hat{\boldsymbol{\theta }}_{T}\left( P\right) -\boldsymbol{\theta }%
\right\Vert $ remains within a range bounded above by the maximum bias
corresponding to the given level of contamination, since the distributions
vary over the entire $\varepsilon $-contamination neighborhood. For $p=1$
the rationale is completely similar, although the bound $B_{L}\left( 
\boldsymbol{\hat{\theta}}_{T},\varepsilon ,\Phi \right) $ is too
large since the maximum bias for the univariate median $\boldsymbol{\hat{\theta}}%
_{M}$ is $B_{L}\left( \mathbf{\hat{\theta}}_{M}\mathbf{,}\varepsilon 
\mathbf{,}\Phi \right) =\Phi ^{-1}\left( 1/(2(1-\varepsilon ))\right)$.
Anyway, for $p\geq 2$, the bound (\ref{tukeymedian}) can be 
derived in a more illuminating manner by explicitly incorporating the 
maximum bias, as the maximum bias governs the behavior of the estimator when
the sample size is sufficiently large, which turns out to be a more
informative inequality without enlarging significantly the error bound in (%
\ref{tukeymedian}). It is well known  the Tukey's median has bounded
contamination sensitivity (\cite{ChenTyler2002}) and therefore it has order $%
\varepsilon $ for $\varepsilon $ near $0$, although $B_{L}\left( \hat{%
\boldsymbol{\theta }}_{T},\varepsilon ,\Phi \right) $ is not of order $%
\varepsilon $ in $\left( 0,1/3\right) $. Table 1 depicts more accurately the incremental quotient $B_{L}\left( \varepsilon \right) /\varepsilon $ as $\varepsilon$ moves in $(0,1/3)$. 
\begin{table}[!hbt]
	\centering
\begin{tabular}{llllllllll}
	\hline
$\varepsilon $ & $0.01$ & $0.05$ & $0.10$ & $0.15$ & $0.20$ & $0.25$ & $0.30$
& $0.33$ & $1/3$ \\ 
$B_{L}\left( \varepsilon \right) $ & $0.03$ & $0.13$ & $0.28$ & $0.46$ & $%
0.67$ & $0.98$ & $1.47$ & $2.43$ & $\infty $ \\ 
$B_{L}\left( \varepsilon \right) /\varepsilon $ & $2.53$ & $2.65$ & $2.82$ & 
$3.05$ & $3.37$ & $3.87$ & $4.88$ & $7.38$ & $\infty $%
\\ \hline
\end{tabular}%
\caption{Behavior of maximum bias $B_L(\varepsilon)$ vs. level of contamination $\varepsilon$}
\end{table}%

Before getting into the details of the proofs, a  simple calculation
displays the effect of the maximum bias in the error bound. We focus on
Tukey's median. Suppose that we have obtained the concentration inequality%
\begin{equation*}
P\left( \left\Vert \boldsymbol{\hat{\theta}}_{n}-\boldsymbol{\theta }%
\right\Vert \leq \Phi ^{-1}\left( \frac{1+\varepsilon }{2\left(
1-\varepsilon )\right) }+40\sqrt{\frac{6e\pi }{1-e^{-1}}}\sqrt{\frac{p+1}{n}}%
+\frac{7}{2}\sqrt{\frac{\log \left( 1/\delta \right) }{n}}\right) \right)
\geq \alpha .
\end{equation*}%
Call $x=b\left( p,n\right) =40\sqrt{\frac{6e\pi }{1-e^{-1}}}\sqrt{\frac{p+1}{%
n}}+\frac{7}{2}\sqrt{\frac{\log \left( 1/\delta \right) }{n}}$. Since $\frac{%
1+\varepsilon }{2\left( 1-\varepsilon )\right) }+x\ $ must be taken less than 
$1$, then $\varepsilon <$ $\frac{1-2x}{3-2x}=\bar{\varepsilon}(x)$. Thus, the
level of contamination can vary up to a certain value in accordance to $x$	
to make a sensible upper bound. Take the level of contamination $%
\varepsilon \left( x\right) =\frac{1-3x}{3-2x}<\bar{\varepsilon}(x),$ $%
B\left( x\right) =B\left( \varepsilon \left( x\right) \right) =\Phi
^{-1}\left( \frac{4-5x}{2\left( 2+x\right) }\right) $ and the error bound
turns out to be $\Phi ^{-1}\left( \frac{4-5x}{2\left( 2+x\right) }+x\right)
=\Phi ^{-1}\left( \frac{4-5x}{2\left( 2+x\right) }+x\right) -\Phi
^{-1}\left( \frac{4-5x}{2\left( 2+x\right) }\right) +\Phi ^{-1}\left( \frac{%
4-5x}{2\left( 2+x\right) }\right) =V\left( x\right) +B\left( x\right)$. By
these means, the error bound comprises two terms exhibiting the usual
trade-off dispersion-bias, $V(x)$ accounting for dispersion through the effect
of the stabilizing rate $\sqrt{p/n}$ and $%
B\left( x\right) $ which measures bias. Table \ref{effectdispersionbias} illustrates
the interaction of dispersion and bias.  
\
\begin{table}[h]
	\centering
	\begin{tabular}{clllllllll}
		\hline
		$x$ & 0.33 & 0.30 & 0.25 & 0.20 & 0.15 & 0.10 & 0.05 & 0.01 & 0.00 \\ 
		$\varepsilon (x)$ & 0.00 & 0.04 & 0.10 & 0.15 & 0.20 & 0.25 & 0.29 & 
		0.33 & $1/3$ \\ 
		$\bar{\varepsilon}(x)$ & 0.15 & 0.17 & 0.20 & 0.23 & 0.26 & 0.29 & 0.31
		& 0.33 & $1/3$ \\ 
		$B(x)$ & 0.01 & 0.11 & 0.28 & 0.47 & 0.69 & 0.97 & 1.37 & 2.11 & $%
		\infty $ \\ 
		$V(x)$ & 0.96 & 0.90 & 0.80 & 0.71 & 0.62 & 0.54 & 0.44 & 0.33 & NaN
		\\ 
		$B\left( x\right) +V\left( x\right) $ & 0.97 & 1.01 & 1.08 & 1.18 & 1.32
		& 1.50 & 1.81 & 2.44 & $\infty $%
		\\ \hline
	\end{tabular}%
	\caption{Interaction between dispersion and bias in the error bound for Tukey's median}
	\label{effectdispersionbias}
\end{table}

Consequently, we can next state an analogous result to that of (\ref{tukeymedian}) but incorporating the
maximum bias of the estimator.

\begin{lemma}
\label{chenlocation}
For $\varepsilon \in \left[ 0,\varepsilon^{\prime }\right]$ 
and $\left( p+\log \left( 1/\delta \right) \right) /n$
sufficiently small, there exists a constant $\tilde{C}>0$ (depending on $\varepsilon
^{\prime }$ but independent of $p,n,\varepsilon ),$ such that 
\begin{equation}
\inf_{\substack{\boldsymbol{\theta },\mathit{\Sigma }\in \mathcal{F}\left( M\right)
,\\ P\in \mathcal{P}_{\varepsilon }\left( P_{0}^{E}\right)} }P\left( \left\Vert 
\hat{\boldsymbol{\theta }}_{T}-\boldsymbol{\theta }\right\Vert ^{2}\leq 
\tilde{C}\left( \max \left\{ \frac{p}{n},B_{L}^{2}\left(\boldsymbol{\hat{\theta}}%
_{T},\varepsilon ,\Phi \right) \right\} +\frac{\log \left( 1/\delta \right) 
}{n}\right) \right) \geq \alpha .  \label{modifiedchen1}
\end{equation}
\end{lemma}
On the other hand, Theorem 3.1 of \cite{ChenGaoRen2018} also derived an
error bound for the deepest estimator for the dispersion
matrix. They showed that, with probability at least $\alpha,$
\begin{equation}
\inf_{\mathit{\Sigma }\in \mathcal{F}\left( M\right) ,P\in \mathcal{P}%
_{\varepsilon }\left( P_{0}^{E}\right) }P\left( \left\Vert \hat{\mathit{\Sigma }}-%
\mathit{\Sigma }\right\Vert _{op}^{2}\leq C\left( \max \left\{ \frac{p}{n}%
,\varepsilon ^{2}\right\} +\frac{\log \left( 1/\delta \right) }{n}\right)
\right) \geq \alpha .  \label{depthinequality}
\end{equation}%
With a similar reasoning to that of Lemma \ref{chenlocation} we can obtain
the following result.

\begin{lemma}
\label{chenscatter} 
Let $\hat{\mathit{\Sigma }}=\beta^{-1}\hat{\mathit{\Gamma }}$. For $\varepsilon \in \left[ 0,\varepsilon^{\prime }\right]$ and 
$\left( p+\log \left( 1/\delta \right)
\right) /n$ sufficiently small, there exists a constant $C^{\ast }>0$ (depending on $%
\varepsilon ^{\prime }$ but independent of $p,n,\varepsilon ),$ such that
\begin{equation}
\inf_{\mathit{\Sigma }\in \mathcal{F}\left( M\right) ,P\in \mathcal{P}%
_{\varepsilon }\left( P_{0}^{E}\right) }P\left( \left\Vert \hat{\mathit{\Sigma }}-%
\mathit{\Sigma }\right\Vert _{op}^{2}\leq C^{\ast }\left( \max \left\{ \frac{%
p}{n},B_{E}^{2}(\varepsilon )\right\} +\frac{\log \left( 1/\delta \right) }{n%
}\right) \right) \geq \alpha ,  \label{modifiedchenscatter}
\end{equation}%
with $B_{E}\left( \varepsilon \right) =\left\{ \frac{1}{\sqrt{\beta }}\Phi
^{-1}\left( a\left( \varepsilon \right) \right) -1\right\}$.
\end{lemma}

\begin{remark}
By these means, the concentration inequality is also able to uncover the
likely maximum bias of the deepest one for $p\geq 2$ since it was derived in Section 4.
\end{remark}

Similarly, in the multivariate regression model, under the assumptions: (a) $\mathbf{X\sim }N_{p}\left( \mathbf{0},\mathit{\Sigma }\right) $, and 
(b) $\left. \mathbf{Y}\right \vert \mathbf{X\sim }N_{m}\left( B^{t}\mathbf{X,%
}\sigma ^{2}I_{m}\right)$, 
the depth estimator $\hat{B}$ for multivariate regression verifies that
there exists a universal constant $C>0$ such that $tr\left( \left( \hat{B}%
_{n}-B\right) ^{t}\mathit{\Sigma }\left( \hat{B}_{n}-B\right) \right) \leq
C\sigma ^{2}\left( \frac{pm}{n}\vee \varepsilon ^{2}\right) ,$ with high
probability and uniformly over the $\varepsilon $-contamination neighborhood
and all $B\in \mathbb{R}
^{p\times m},$ see Theorem 4.1 of \cite{Gao2020}.

We can reformulate this result similarly to those of (\ref{modifiedchen1})
and (\ref{modifiedchenscatter}). Take the functions $h:[0,\infty )\times %
\left[ 0,1\right] ^{p}\rightarrow \left[ 0.5,1\right] $ and $g:[0,\infty
)\rightarrow \left[ 0.5,1\right] $ defined as $h\left( t,\theta _{1},\dots
,\theta _{p}\right) =\Phi \left( t\sqrt{\sum_{i=1}^{p}\theta _{i}Z_{i}^{2}}%
\right) $ and $g\left( t\right) =h\left( t,1,0,...,0\right) .$ Now, note
that $g^{-1}\left( \frac{1+\varepsilon }{2\left( 1-\varepsilon \right) }%
\right) :\left[ 0,1/3\right) \rightarrow \lbrack 0,\infty )$ and put $%
b_{MR}^{2}\left( \varepsilon \right) =g^{-1}\left( \frac{1+\varepsilon }{
2\left( 1-\varepsilon \right) }\right) .$ Then, the following result follows.

\begin{lemma}
\label{gaomultregrmodified} 
Suppose that (a) and (b) hold. Let $\mathcal{P}%
_{\varepsilon }\left( H_{B,\sigma }\right) $ be the $\varepsilon $
-contamination neighborhood with $H_{B,\sigma }$ as in (\ref%
{gaomultregrmodified}). For $\varepsilon \in \left[ 0,\varepsilon ^{\prime }%
\right]$, 
and $\left( p+\log \left( 1/\delta
\right) \right) /n$ sufficiently small, there exists a constant $C>0$
(depending on $\varepsilon ^{\prime }$ but independent of $p,n,\varepsilon
)$, such that
\begin{equation*}
\inf_{\substack{ \boldsymbol{B\in 
\mathbb{R}}^{p\times m},  \\ P\in \mathcal{P}_{\varepsilon }\left( H_{B,\sigma
}\right) }}P\left[ tr\left( \left( \hat{B}_{n}-B\right) ^{t}\mathit{\Sigma }%
\left( \hat{B}_{n}-B\right) \right) \leq C\sigma ^{2}\left( \frac{pm}{n}\vee
b_{MR}^{2}\left( \varepsilon \right) \right) +\frac{\log \left( 1/\delta
\right) }{n}\right] \geq \alpha .
\end{equation*}
\end{lemma}

\cite{AMY2002} derived the maximum bias of depth estimators in case of the
univariate regression setting with known intercept and covariables with
spherical distribution. They proved that given a bivariate normal vector $%
\left( U,V\right) $ with zero mean, $Var\left( U\right) =Var\left( V\right)
=1$ and $\rho =Corr\left( U,V\right) ,$ the function $h\left( \rho \right)
=P\left( sg\left( U\right) =sg\left( V\right) \right) $ determines the
maximum bias $b$, since it solves the equation %
$h^{-1}\left( 1+\varepsilon\right) / \left( 2\left(
1-\varepsilon \right)  \right) = b/(\sqrt{1+b^{2}})$.
The following lemma shows that such a $b$ solving the equation should be $%
g^{-1} \left\{(1+\varepsilon) /(2\left( 1-\varepsilon \right))\right\} $.

\begin{lemma}
\label{biasunivariateregression} Let $m=1$ and suppose that (a) and (b) hold.
Then $B_{MR}\left( \boldsymbol{\hat{\beta}},\varepsilon ,H_{0,1}\right)$
$=$ $g^{-1}\left( \frac{1+\varepsilon }{2\left( 1-\varepsilon \right) }\right) .$
\end{lemma}

\begin{remark}
The concentration inequality given by Gao (2020) contains the information
regarding  the maximum bias.
\end{remark}

\section{Depth as a residual smallness concept}

So far we have shed light on the relationship between concentration
inequalities and the concept of asymptotic maximum bias for depth estimators
in the multivariate and regression setting. We next want to state a unified
view for all the definitions of depth considered in Section 2. \cite%
{Carrizosa1996} and \cite{AMY2002}\emph{\ }independently came up with a
residual smallness concept which comprises the notion of depth given in the
univariate and multivariate model as well as the univariate regression
model. More precisely, the depth of $\boldsymbol{\theta }\in \mathbf{\mathbb{R}}^{p}$ can be defined as 
\begin{equation}
\mathcal{D}_{T}^{E}\left( \boldsymbol{\theta },P\right) {\ =}%
\inf_{\left\Vert \boldsymbol{\lambda }\right\Vert =1,\boldsymbol{\gamma }%
\mathbf{\in \mathbb{R}
}^{p}}P\left( \left\vert \boldsymbol{\lambda }^{t}(\mathbf{x-}\boldsymbol{\
\theta }\mathbf{)}\right\vert \leq \left\vert \boldsymbol{\lambda }^{t}(%
\mathbf{x-}\boldsymbol{\gamma }\mathbf{)}\right\vert \right)
\label{residualsmallness}
\end{equation}%
It is proved that $\mathcal{D}_{T}^{E}\left( \boldsymbol{\theta },P\right) =%
\mathcal{D}_{T}\left( \boldsymbol{\theta },P\right) .$ The idea behind that
definition is that the depth of a fit $\boldsymbol{\theta }$ is determined
by the bad performance displayed by the residuals $\left\vert \boldsymbol{\
\lambda }^{t}(\mathbf{x-}\boldsymbol{\theta }\mathbf{)}\right\vert $
compared to the best competitor $\boldsymbol{\gamma ,}$ whose residuals $%
\left\vert \boldsymbol{\lambda }^{t}(\mathbf{x-}\boldsymbol{\theta }\mathbf{)%
}\right\vert $ have the minimum probability of being worse than those of $%
\left\vert \boldsymbol{\lambda }^{t}(\mathbf{x-}\boldsymbol{\theta }\mathbf{)%
}\right\vert $, therefore the $\boldsymbol{\theta }$ with the best worst
performance is singled out.
In the regression setting, (\ref{residualsmallness}) takes the form 
\begin{equation}
\mathcal{D}_{R}^{E}\left( \boldsymbol{\theta },P\right) {=}\inf_{\boldsymbol{%
\gamma }\mathbf{\in 
\mathbb{R}}^{p}}P\left( \left \vert y-\boldsymbol{\theta }^{t}\mathbf{x}\right \vert
\leq \left \vert y-\boldsymbol{\gamma }^{t}\mathbf{x}\right \vert \right)
\label{residualsmallness2}
\end{equation}%
and the deepest regression estimator is taken to be 
$\hat{\boldsymbol{\theta}}_{R}\left( P\right) =\arg \sup_{\boldsymbol{\theta } }
\mathcal{D}_{R}^{E}\left( \boldsymbol{\theta },P\right)$ .

The equivalence between 
(\ref{residualsmallness2}) and the regression depth 
can be found in \cite{AMY2002}. In the
multivariate regression model we can also generalize the concept of depth
given in (\ref{residualsmallness2}) by defining,%
$\mathcal{D}_{MR}^{E}\left( B,P\right) {\ =}\inf_{U\mathbf{\in 
\mathbb{R}
}^{p\times m}}P\left( \left \Vert Y-B^{t}\mathbf{X}\right \Vert \leq \left
\Vert Y-U^{t}\mathbf{X}\right \Vert \right). \ \
$
 Then, it can be easily seen that this approach coincides with that of (\ref%
{multregmizera}).

\begin{lemma}
\label{multregrequivalence}$\mathcal{D}_{MR}^{E}\left( B,P\right) =\mathcal{%
D }_{MR}\left( B,P\right) .$
\end{lemma}

The residual smallness concept in (\ref{residualsmallness}) can be easily
adapted for joint estimation of location and scale as follows. In the
univariate case the MLSM 
switches to the usual
location-scale model $Y=\mu _{0}+\sigma _{0}U$, $\sigma _{0}>0$. Then, the
depth of $\left( \mu ,\sigma \right) \in 
\mathbb{R}\times (0,\infty )$ is taken to be,%
\begin{eqnarray}
\mathcal{D}_{LS}\left( \mu ,\sigma ,P\right) &=&\min \left\{ 
\begin{array}{c}
\inf_{\lambda \in \mathbb{R}
}P\left(  \left\vert Y-\mu \right\vert \leq \left\vert Y-\lambda
\right\vert \right) , \\ 
\inf_{\gamma >0}P\left( \left\vert \left\vert \frac{Y-\mu }{\sigma }%
\right\vert -1\right\vert \leq \left\vert \left\vert \frac{Y-\mu }{\gamma }%
\right\vert -1\right\vert \right)%
\end{array}%
\right\}  \label{locscaledepthone} \\
\left( \hat{\mu}_{1},\hat{\sigma}_{1}\right) &=&\arg \max_{\mu ,\sigma
}D_{LS}^{1}\left( \mu ,\sigma ,P\right) .  \notag
\end{eqnarray}%
With this definition of depth, we obtain very well known functionals for
location and scale, as it is stated in the following lemma.

\begin{lemma}
\label{residsmallocscale} If $Y\sim P,$ then $\hat{\mu}_{1}=med_{P}\left(
Y\right) $ and $\hat{\sigma}_{1}=med_{P}\left( \left\vert Y-med_{P}\left(
Y\right) \right\vert \right) $ and $D_{LS}\left( \hat{\mu}_{1},\hat{\sigma}%
_{1},P\right) \geq 0.5.$
\end{lemma}

If $P$ is taken to be the empirical distribution function we come up with the usual
median and median absolute deviation around the median (MADM)as location and scale
estimators. \cite{MZ1993A} showed that a scaled version of the MADM is approximately minimax bias-robust within the class of Huber's Proposal 2 joint estimates of location and scale. 

In the multivariate setting, we can adjust similarly the definition given en
(\ref{locscaledepthone}) and we take%
\begin{eqnarray*}
\mathcal{D}_{LS}^{E}\left( \boldsymbol{\mu },\mathit{\Gamma },P\right) &=&\min
\left\{ 
\begin{array}{c}
\inf_{\mathbf{u,}\boldsymbol{\lambda }\in \mathbb{R}
^{p}}P\left( \left\vert \mathbf{u}^{t}\left( \mathbf{X}-\boldsymbol{%
\mu }\right) \right\vert \leq \left\vert \mathbf{u}^{t}\left( \mathbf{X}-%
\boldsymbol{\lambda }\right) \right\vert \right) , \\ 
\inf_{\mathit{\Theta }\in \mathcal{E},\mathbf{u}\in \mathbb{R}
^{p}}P\left( \left\vert \left\vert \frac{\mathbf{u}^{t}\left( \mathbf{%
X}-\boldsymbol{\mu }\right) }{\sqrt{\mathbf{u}^{t}\mathit{\Gamma }\mathbf{u}}%
}\right\vert -1\right\vert \leq \left\vert \left\vert \frac{\mathbf{u}%
^{t}\left( \mathbf{X}-\boldsymbol{\mu }\right) }{\sqrt{\mathbf{u}^{t}\mathit{%
\Theta }\mathbf{u}}}\right\vert -1\right\vert \right)%
\end{array}%
\right\} \\
\left( \boldsymbol{\hat{\mu}},\mathit{\Gamma }\right) &=&\arg \max_{\boldsymbol{%
\mu },\mathit{\Gamma }}\mathcal{D}_{LS}^{E}\left( \boldsymbol{\mu },\mathit{%
\Gamma },P\right) .
\end{eqnarray*}%
If the location is known, the depth of $\mathit{\Gamma }\in \mathcal{E}$ is
taken to be 
\begin{equation*}
\mathcal{D}_{LS}^{E}\left( \mathit{\Gamma },P\right) =\inf_{%
\mathit{\Theta }\in \mathcal{E},\mathbf{u}\in 
\mathbb{R}^{p}}P\left( \left\vert \left\vert \frac{\mathbf{u}^{t}\mathbf{X}}{%
\sqrt{\mathbf{u}^{t}\mathit{\Gamma }\mathbf{u}}}\right\vert -1\right\vert
\leq \left\vert \left\vert \frac{\mathbf{u}^{t}\mathbf{X}}{\sqrt{\mathbf{u}%
^{t}\mathit{\Theta }\mathbf{u}}}\right\vert -1\right\vert  \right)
\end{equation*}%
From (\ref{locscaledepthone}) we obtain that $\mathcal{D}%
_{LS}^{E}\left( \mathit{\Gamma },P\right) =\mathcal{D}_{LS}\left( \mathit{%
\Gamma },P\right) ~$. Therefore the depths functions considered for
multivariate scatter and regression can be embodied into this framework of
residual smallness, as (\ref{residualsmallness}) and (\ref%
{residualsmallness2}).

\FloatBarrier

\section{Numerical study}

We derived the maximum bias curve for the deepest scatter matrix as well as
its breakdown point. Only a few maxbias curves have been derived for
multivariate dispersion measures, despite their being a very informative
overall measure of robustness. Therefore, we cannot rely on comparisons
among many theoretical maxbias curves to understand the behaviour of
dispersion estimators.

On the one hand, maxbias curves are population measures; that is, they
describe the behavior of estimators as functionals on contamination
neighborhoods. On the other hand, finite-sample effects should also be taken
into account. We next report the results of a Monte Carlo simulation study
to investigate these effects, explicitly including the estimation of the
unknown location parameter, whose impact is typically avoided in the maxbias
derivation for scatter matrices because of the intractability to perform a
theoretical analysis.

In \cite{MaronnaYohai2017} and \cite{MaronnaMartinYohaiSalibian2018}, the
behavior of scatter matrix estimators is studied under the assumption of
known location, using the Kullback-Leibler divergence as a performance
measure. \cite{ChenGaoRen2018} addresses the performance of estimators under
unknown location by using the operator norm. \cite{HRK2014} compares the
worst-case bias of several prominent robust multivariate estimators by means
of simulation by using the condition number as bias measure under unknown
location. We consider the empirical version of the bias given in formulae \eqref{biasexplo} and  \eqref{biasimplo} as performance measures (see Section \ref{maxbiasdifstat}) and the condition number
of an estimator as a measure of asymptotic bias).

In the following subsections, we describe the estimators and contamination
scenarios, explain how the estimated bias measures are empirically computed
from replicated samples of the contaminated model, and then report and
discuss the main results.

\subsection{The estimators}

We include several estimators in our simulation study. To compute them, we
use functions from R packages available on the Comprehensive R Archive
Network (CRAN), with default argument values in all cases. We selected
packages whose implementations ensure Fisher consistency of the estimators.
The robust estimators for multivariate scatter 
under consideration are the sample covariance matrix (\textsc{SCOV)},
the minimum volume ellipsoid estimator (\textsc{MVE}), the minimum covariance determinant estimator (\textsc{MCD}), S-estimator for multivariate location and scatter (\textsc{SE)},  S-estimators with non-monotonic weight functions (\textsc{Rocke}), MM-estimators (\textsc{MM)}, the Stahel-Donoho location scatter estimator (\textsc{SD}), and the Deepest Estimator (\textsc{MDepth}). A brief description of them is available at the Appendix. We do not compute Tyler's M-estimator of scatter, since we are only using routines available at the CRAN Project and the R package ICSNP takes the
sample mean as multivariate location estimator by default, which makes the
comparison with other robust proposals unreliable.

\subsection{Contamination scenarios and empirical bias}

Let $\mathcal{P}_{\varepsilon }\left( P_{0}\right) $ be the $\varepsilon $%
-contamination neighborhood with $P_{0}=N(\boldsymbol{\theta },\mathit{%
\Sigma })$, $\boldsymbol{\theta }=\mathbf{0}$ and $\mathit{\Sigma }=I_{p}$,
the $p\times p$ identity matrix. We assume that the contaminating
distribution is $G_{k}=\delta _{x_{k}}$, a point mass at $\mathbf{Z}%
_{k}=(k,\ldots ,k)^{t}$, with $k\in \mathbb{N}\cup \{0\}$. Consider $%
\boldsymbol{X}^{0}=(X_{1}^{0},\ldots ,X_{p}^{0})^{\prime }\sim P_{0}$ and $%
B\sim Ber\left( \varepsilon \right) .$ If $\boldsymbol{X}^{0}$ and $B$ are
independent then $
\boldsymbol{X}=(1-B)\boldsymbol{X}^{0}+B\boldsymbol{Z}_{k}\sim P_{k}$,
with $P_{k}=(1-\varepsilon )P_{\mathbf{0}}\left( \cdot \right) +\varepsilon
G_{k}(\cdot )\in \mathcal{P}_{\varepsilon }\left( P_{0}\right)$. For a
contamination rate $\varepsilon $ and a fixed constant $k$, 
we generate a sample $\{\mathbf{x}_{1},\ldots ,\mathbf{x}_{n}\}$
from the random vector $\mathbf{X}\sim P_{k}$.

We will consider the following scenarios: (a) Contamination proportions $\varepsilon = 0.1$, $0.2$; (b) $k$ $\in $ $\left \{ 0,1,5,10,15,20,25\right \} $; (c) Dimensions $p=2,5$, $10$, $15$, and (d) Sample size: $n=10p$ (see \cite{MaronnaYohai2017}, Table 3), $n=40p$, $n=500p$. 
For each combination of $p$, $\varepsilon$, and $k$, we generate $R = 50$
independent data sets, denoted as $\mathbb{X}_r = \{ \mathbf{x}_1^{(r)},
\ldots, \mathbf{x}_n^{(r)} \}$, for $r = 1, \ldots, R$.


The effect of the distortion caused by having $P\in \mathcal{P}_{\varepsilon
}\left( P_{0}\right) $ was measured through formulas \ref{biasexplo} and \ref{biasimplo}. Following 
\cite{MaronnaYohai1995}, we also use the condition number (CN) of an
estimator $\mathit{\Gamma }$ as an additional measure of bias, abbreviated
as $B_{\text{CN}}$, and defined as the supremum of $\text{CN associated with 
}\mathit{\Gamma }$ over the contamination neighbourhood. 
Since the conclusions regarding to $B_{\text{CN}}$ are similar to those of using \ref{biasexplo} and \ref{biasimplo}, tables and plots related to $B_{\text{CN}}$ are deferred to the Appendix. 

To obtain empirical versions of $B$ and $B_{\text{CN}}$, we proceed as
follows:

\begin{itemize}
\item \textbf{Step 1.} For every $k=0,1,\ldots ,25$, generate $R=50$
independent data sets, denoted by 
$\mathbb{X}_{k}^{(r)}=\{\mathbf{x}_{1}^{(r)},\ldots ,\mathbf{x}%
_{n}^{(r)}\},\quad r=1,\ldots ,R$, 
where each $\mathbb{X}_{k}^{(r)}$ is a sample from the distribution $P_{k}$.
Let $\mathit{\Gamma }_{k}^{(r)}$ denote the corresponding estimated scatter
matrix. For every $\mathit{\Gamma }_{k}^{(r)}$, compute $\hat{b}_{E,k}^{(r)}$
and $\hat{b}_{I,k}^{(r)}$, $r=1,\ldots ,R$, by using the empirical version
of $P_{k}$ corresponding to $\mathbb{X}_{k}^{(r)}$, denoted by $%
P_{n,k}^{(r)} $. Then, define the bias of the scatter matrix $\mathit{\Gamma 
}_{k}^{(r)}$, based on the $r$-th sample from $P_{k}$, as 
$\hat{b}_{k}^{(r)}=\max \left\{ \hat{b}_{I,k}^{(r)},\hat{b}%
_{E,k}^{(r)}\right\} =\max \left\{ \lambda _{(1)}^{(r)},\lambda
_{(p)}^{(r)-1}\right\}$, 
where $\lambda _{(1)}^{(r)}$ and $\lambda _{(p)}^{(r)}$ denote the largest
and the smallest eigenvalues of $\mathit{\Gamma }_{k}^{(r)}$. For $1\leq
r\leq R,$ let CN$_{k}^{(r)}=$CN$\left( \mathit{\Gamma }_{k}^{(r)}\right)
=\lambda _{(1)}^{(r)}/\lambda _{(p)}^{(r)}$ be the condition number of the
estimated scatter matrix $\mathit{\Gamma }_{k}^{(r)}.$ If $l$ is a location
estimator, then $\hat{b}_{k}$ and CN$_{k}$ stand for 
\begin{equation}
\hat{b}_{k}=l\left( \{\hat{b}_{k}^{(r)}\}\right) \text{ and } CN_{k}=l\{CN_{k}^{(r)}\}. \label{bk}
\end{equation}%
\item \textbf{Step 2.} The measures 
$\hat{B}=\max_{0\leq k\leq 25}\{\hat{b}_{i,k}\}$ and $\hat{B}_{%
\text{CN}}=\max_{0\leq k\leq 25}\{CN_{k}\}$,  
are empirical approximations to the maximum biases $B$ and $B_{\text{CN}}$,
respectively.

\end{itemize}

In order to place $\hat{b}_{I,k}^{(r)},\hat{b}_{E,k}^{(r)}$ and 
$CN_{k}^{(r)}$ on a more comparable scale, we consider their logarithms. When
sampling from a contaminated distribution $P_{k},$ the empirical
distributions of $\log \left( \hat{b}_{I,k}^{(r)}\right)$, $\log \left( \hat{%
	b}_{E,k}^{(r)}\right) $ and $\log \left( CN_{k}^{(r)}\right) $ exhibit
asymmetry and the presence of outliers, as shown by the boxplots in Figures 
\ref{box_B_k1_p_2_n_20} and \ref{box_B_k0_p_15_n_150}. Therefore, it is more
appropriate to use the median as the location measure $l$ in (\ref{bk}) 
rather than the mean. For the sake of simplicity, we refer to $\hat{B}$ and $\hat{B}_{\mathrm{CN}}$ as the “maximum medians” of $B$ and $CN$, respectively, in the tables and figures. Further details are provided in Subsection \ref{rc}.

\subsection{Efficiency}

When the data come from the symmetric central model, the empirical
distributions of $\log \left( CN_{k}^{(r)}\right) $ and $\log \left(				
b_{k}^{(r)}\right) $\ tend to be more symmetric and the mean is therefore
more representative. Thus, given a scatter estimator, we can take $\hat{B}$\
and $\hat{B}_{\mathrm{CN}}$\ as defined in Step 2., 
but using the
mean rather than the median in (\ref{bk}), 
which will be called the (empirical) mean absolute error, abbreviated as MAE. 
Furthermore, if $\mathit{\Gamma }$\ stands for a scatter estimator and $S_{n}$ the
sample covariance matrix, the ratio of the MAE for $S_{n}$\ to the
corresponding value for $\mathit{\Gamma }$, 
$\text{Eff}=\mathrm{MAE}(S_{n})/\mathrm{MAE}(\mathit{\Gamma })$,
may be considered a measure of efficiency at the normal central model $P_{0}$%
; see \cite{MaronnaYohai1995}.

\subsection{Results and conclusions}\label{rc}

We discuss the results presented in this subsection together with those reported in Appendix B.   In  Tables \ref{eff_n_50_B_mean} and \ref{eff_n_200_B_mean} the
efficiencies of scatter estimators are computed for different combinations
of $n$ and $p$. For most estimators, efficiency remains stable or even
slightly increases with growing $p$, suggesting favorable finite-sample
behavior when the ratio $p/n$ is moderate. In particular, the MM estimator
maintains high efficiency across all considered dimensions and sample sizes.

For the SE, for small dimensions the efficiency is low, but it
increases considerably as $p$ grows. This is a phenomenon described in the
literature; see, for instance, Section 6.4.4. of \cite%
{MaronnaMartinYohaiSalibian2018}. As discussed therein, when $p$ is large
enough, almost all observations receive similar weights, except for
observations far from the bulk of the data, yielding an estimator which
closely approximates the sample covariance matrix. But increasing efficiency
usually entails a loss of robustness when the dimension grows.

\begin{table}[h]
\centering
\begin{tabular}{rrrrrrrrr}
\hline
$p$ & MVE & MCD & SE & ROCKE & MM & SD & MDEPTH &  \\ \hline
$2$ & 0.62 & 0.69 & 0.58 & 0.43 & 0.81 & 0.62 & 0.28 &  \\ 
$5$ & 0.54 & 0.69 & 0.80 & 0.52 & 0.92 & 0.56 & 0.40 &  \\ 
$10$ & 0.54 & 0.68 & 0.88 & 0.44 & 0.95 & 0.58 & 0.47 &  \\ \hline
\end{tabular}%
\caption{$\text{Eff}$ in log-scale based on $B$ (means) for $n=50$ over
dimensions.}
\label{eff_n_50_B_mean}
\end{table}

\begin{table}[ht]
\centering
\begin{tabular}{rrrrrrrrr}
\hline
$p$ & MVE & MCD & SE & ROCKE & MM & SD & MDEPTH &  \\ \hline
$2$ & 0.704 & 0.67 & 0.58 & 0.45 & 0.87 & 0.64 & 0.15 &  \\ 
$5$ & 0.72 & 0.80 & 0.85 & 0.62 & 0.90 & 0.72 & 0.25 &  \\ 
$10$ & 0.77 & 0.81 & 0.94 & 0.65 & 0.93 & 0.75 & 0.33 &  \\ 
$15$ & 0.76 & 0.82 & 0.97 & 0.64 & 0.92 & 0.71 & 0.38 &  \\ \hline
\end{tabular}%
\caption{$\text{Eff}$ in log-scale based on $B$ (means) for $n=200$ over
dimensions.}
\label{eff_n_200_B_mean}
\end{table}

\begin{figure}[h]
	\centering 
	\includegraphics[width=\textwidth]{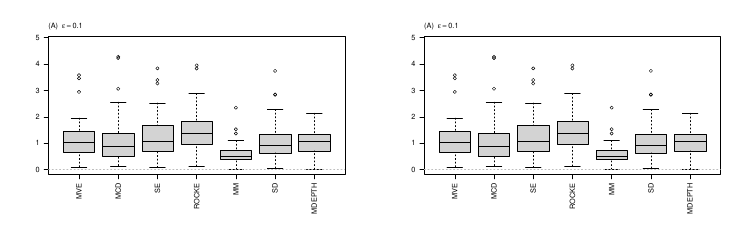}
	\caption{Boxplots of the biases $\{ \hat{b}_{k}^{(r)}\}_{1\leq r\leq R}$
		(log-scale) through the $R=50$ replicates for each scatter estimator with $%
		p=2$, $n=20$ and $k=1$. Panels correspond to (A) $\varepsilon=0.1$ and 
		(B) $\varepsilon=0.2$.}
	\textbf{Alt text:} Two side-by-side boxplots comparing bias distributions across several scatter estimators. Panel A shows results for contamination level $\varepsilon=0.1$, and panel B for $\varepsilon=0.2$. Bias distributions become more asymmetric under higher contamination.	
	\label{box_B_k1_p_2_n_20}
\end{figure}

\begin{figure}[h]
	\centering 
	\includegraphics[width=\textwidth]{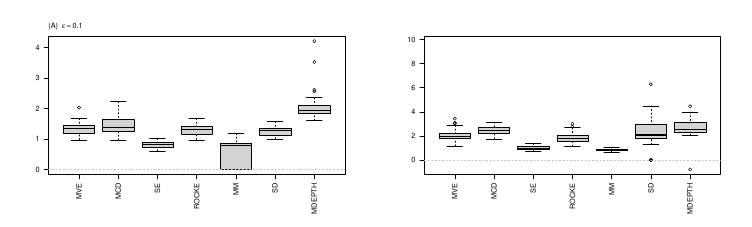}
	\caption{Boxplots of the biases $\{ \hat{b}_{k}^{(r)}\}_{1\leq r\leq R}$
		(log-scale) through the $R=50$ replicates for each scatter estimator with $
		p=15$, $n=150$ and $k=0$ }
	\textbf{Alt text:} Two side-by-side boxplots comparing bias distributions across several scatter estimators. Panel A shows results for contamination level $\varepsilon=0.1$, and panel B for $\varepsilon=0.2$. Bias distributions become more asymmetric under higher contamination.	
	\label{box_B_k0_p_15_n_150}
\end{figure}

Figures \ref{box_B_k1_p_2_n_20} and \ref{box_B_k0_p_15_n_150} show boxplots
of the biases $\{ \hat{b}_{k}^{(r)}\}$, on a log scale, across the $R$
replicates for each scatter estimator, under selected scenarios
(combinations of $p$, $n$, and $k$). The boxplots highlight the need for
using medians rather than means as representative measures of the behavior
of the bias measures. Note that most empirical distributions exhibit heavy
tails, skewness, and outliers, especially when the contamination level is $%
\varepsilon =0.2$.

\begin{figure}[htbp]
	\centering
	
	\includegraphics[width=\textwidth]{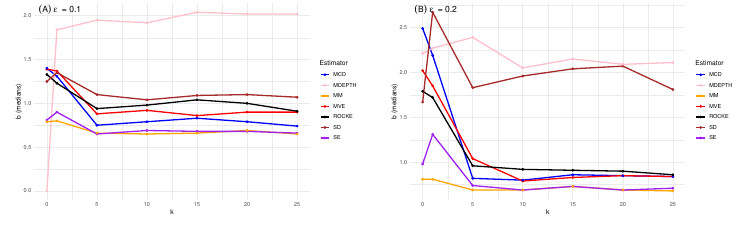}
	
	\caption{$\hat{b}_{k}=\mathrm{median}_{1\leq r\leq R}\{ \hat{b}_{k}^{(r)} \}$ (log-scale) versus $k$ for each scatter estimator, under contamination levels $\varepsilon$. Dimension $p=10$ and $n=100$.}
	
	\textbf{Alt text:} Two side-by-side line plots showing the median values of the bias as a function of $k$ for several scatter estimators on a log scale. Panel A corresponds to a lower contamination level and panel B to a higher contamination level.
	\label{b_K_p_10_n_100}
	
\end{figure}

\begin{figure}[htbp]
	\centering
	\includegraphics[width=\textwidth]{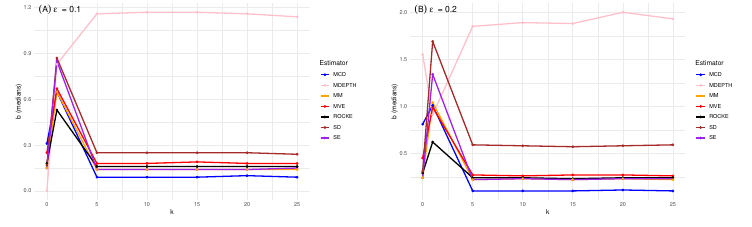}
	
	\caption{$\hat{b}_{k}=\mathrm{median}_{1\leq r\leq R}\{ \hat{b}_{k}^{(r)} \}$ (log-scale) versus $k$ for each scatter estimator, under contamination levels $\varepsilon$. Dimension $p=10$ and $n=5000$.}
	\textbf{Alt text:} Two side-by-side line plots showing the median values of the bias as a function of $k$ for several scatter estimators on a log scale. Panel A corresponds to a lower contamination level and panel B to a higher contamination level.
	\label{b_K_p_10_n_5000}

\end{figure}

Figures \ref{b_K_p_10_n_100} and \ref{b_K_p_10_n_5000} display the behavior of $%
\hat{b}_{k}=\text{median}_{1\leq r\leq R}\{ \hat{b}_{k}^{(r)}\}$ (see (\ref{bk})), for each scatter estimator, as functions of $k$,
for selected dimensions $p$ and sample sizes $n$, under both contamination
levels. Tables \ref{b_p_2} to \ref{b_p_15}  
show the maximum medians of 
$B$ for each scatter estimator across sample sizes $n$ and contamination
levels $\varepsilon $ for dimension $p=2,\ldots ,15$. From the tables we can
conclude that MM tends to give the smallest maximum medians of $B$ for small
and moderate sample sizes and for lower dimensions, showing the best overall
performance. ROCKE often outperforms MM for larger sample sizes and in
higher dimensions ($p\geq 10$). The figures suggest that, in low-dimensional
settings ($p=2$), MM and MCD perform best. For larger sample sizes, MM
remains the best-performing estimator, followed by SE and ROCKE, while MCD
also shows competitive performance.
Generally speaking, the bias curves decrease as $k$ grows,
reflecting the robustness of these estimators to extreme outliers.

\FloatBarrier
\begin{table}[ht]
\centering
\begin{tabular}{ccrrrrrrrr}
\hline
$n$ & $\varepsilon$ & SCOV & MVE & MCD & SE & ROCKE & MM & SD & MDEPTH \\ 
\hline
20 & 0.10 & 4.77 & 1.06 & 1.37 & 1.48 & 1.65 & 0.78 & 1.30 & 1.65 \\ 
20 & 0.20 & 5.35 & 1.73 & 1.67 & 1.90 & 2.20 & 0.82 & 1.60 & 1.82 \\ 
80 & 0.10 & 4.74 & 0.50 & 0.48 & 0.51 & 0.65 & 0.34 & 0.48 & 1.18 \\ 
80 & 0.20 & 5.32 & 0.77 & 1.14 & 0.99 & 1.14 & 0.52 & 0.93 & 1.62 \\ 
1000 & 0.10 & 4.76 & 0.32 & 0.48 & 0.41 & 0.51 & 0.21 & 0.40 & 1.14 \\ 
1000 & 0.20 & 5.31 & 0.54 & 1.18 & 0.80 & 0.90 & 0.52 & 0.77 & 1.52 \\ \hline
\end{tabular}%
\caption{Empirical maximum bias $\hat{B}$ (log scale) for each scatter
estimator across sample sizes $n$ and contamination levels $\protect%
\varepsilon$ for dimension $p=2$.}
\label{b_p_2}
\end{table}

\begin{table}[ht]
\centering
\begin{tabular}{llrrrrrrrr}
\hline
$n$ & $\varepsilon$ & SCOV & MVE & MCD & SE & ROCKE & MM & SD & MDEPTH \\ 
\hline
50 & 0.10 & 5.75 & 1.41 & 1.24 & 0.89 & 1.55 & 0.69 & 1.34 & 1.73 \\ 
50 & 0.20 & 6.24 & 2.49 & 2.35 & 1.42 & 2.26 & 0.78 & 2.09 & 2.32 \\ 
200 & 0.10 & 5.69 & 0.58 & 0.67 & 0.55 & 0.66 & 0.41 & 0.70 & 1.37 \\ 
200 & 0.20 & 6.22 & 0.80 & 1.34 & 0.87 & 0.96 & 0.63 & 1.23 & 1.82 \\ 
2500 & 0.10 & 5.65 & 0.42 & 0.36 & 0.49 & 0.46 & 0.29 & 0.53 & 1.15 \\ 
2500 & 0.20 & 6.22 & 0.61 & 0.80 & 0.83 & 0.60 & 0.47 & 1.04 & 1.54 \\ \hline
\end{tabular}%
\caption{Empirical maximum bias $\hat{B}$ (log scale) for each scatter
estimator under contamination levels $\protect\varepsilon$ and different
sample sizes $n$ for dimension $p=5$.}
\label{b_p_5}
\end{table}

\begin{table}[ht]
\centering
\begin{tabular}{ccrrrrrrrr}
\hline
$n$ & $\varepsilon$ & SCOV & MVE & MCD & SE & ROCKE & MM & SD & MDEPTH \\ 
\hline
100 & 0.10 & 6.34 & 1.39 & 1.40 & 0.90 & 1.33 & 0.77 & 1.35 & 2.04 \\ 
100 & 0.20 & 6.88 & 2.02 & 2.49 & 1.31 & 1.79 & 0.81 & 2.67 & 2.39 \\ 
400 & 0.10 & 6.30 & 0.71 & 0.69 & 0.86 & 0.61 & 0.66 & 0.97 & 1.53 \\ 
400 & 0.20 & 6.89 & 1.02 & 1.50 & 1.31 & 0.77 & 1.01 & 1.85 & 2.32 \\ 
5000 & 0.10 & 6.33 & 0.67 & 0.65 & 0.85 & 0.53 & 0.65 & 0.87 & 1.17 \\ 
5000 & 0.20 & 6.91 & 1.40 & 1.01 & 1.34 & 0.62 & 1.04 & 1.69 & 2.00 \\ \hline
\end{tabular}%
\caption{Empirical maximum bias $\hat{B}$ (log scale) for each scatter
estimator under different sample sizes $n$ and contamination levels $\protect%
\varepsilon$ for dimension $p=10$.}
\label{b_p_10}
\end{table}

\begin{table}[th]
\centering
\begin{tabular}{rrrrrrrrrr}
\hline
$n$ & $\varepsilon$ & SCOV & MVE & MCD & SE & ROCKE & MM & SD & MDEPTH \\ 
\hline
150 & 0.10 & 6.68 & 1.35 & 1.37 & 1.10 & 1.32 & 0.87 & 1.79 & 2.37 \\ 
150 & 0.20 & 7.28 & 1.94 & 2.43 & 1.69 & 1.77 & 1.05 & 3.47 & 2.62 \\ 
600 & 0.10 & 6.73 & 1.24 & 0.92 & 1.15 & 0.56 & 0.93 & 1.31 & 2.23 \\ 
600 & 0.20 & 7.31 & 1.32 & 1.56 & 1.65 & 0.71 & 1.35 & 2.36 & 1.90 \\ 
7500 & 0.10 & 6.73 & 0.89 & 0.88 & 1.12 & 0.50 & 0.92 & 0.89 & 1.23 \\ 
7500 & 0.20 & 7.31 & 1.27 & 1.35 & 1.66 & 0.55 & 1.39 & 1.46 & 2.02 \\ \hline
\end{tabular}%
\caption{Empirical maximum bias $\hat{B}$ (log scale) for each scatter
estimator under different sample sizes $n$ and contamination levels $\protect%
\varepsilon $ for dimension   $p=15$.}
\label{b_p_15}
\end{table}

\FloatBarrier

\FloatBarrier

\section{Concluding remarks}
It is well known in the robust statistics literature that the maximum bias of an estimator plays a fundamental role in its mean squared error for moderate to large sample sizes. Consequently, the derivation of the maximum bias curve provides substantial insight into the robustness properties of an estimator, including two closely related concepts: contamination sensitivity and breakdown point.
Since the introduction of Tukey’s median and halfspace depth, several notions of depth have been proposed across different statistical models, leading to suitable extensions for multivariate scatter matrices. In this work, we derive the maximum bias function of the deepest scatter estimator under the $\varepsilon$-contamination neighborhood, as well as its breakdown point and contamination sensitivity, assuming a normal distribution for the central model. 
The work of \cite{ChenGaoRen2018} establishes error bounds governed by the rate $\sqrt{p/n}+\varepsilon$ for Tukey’s median and deepest scatter estimators in several settings. We emphasize that such bounds should reflect the maximum bias rather than $\varepsilon$, which only approximates the maximum bias for small values of $\varepsilon$. In particular, we analyze error bounds for Tukey’s median, the deepest scatter estimator, and the deepest multivariate regression estimator, highlighting the role of maximum bias in each case. Furthermore, the paper by \cite{hesimpson1993} shows that lower bounds for the maximum bias of equivariant estimators are governed by the so-called variation gauge, rather than directly by the contamination level $\varepsilon$. The variation gauge also plays a key role in the derivation of lower error bounds by \cite{ChenGaoRen2018}.  This observation challenges the interpretation of $\sqrt{p/n}+\varepsilon$ as a minimax rate, since the optimal behavior is determined by the variation gauge rather than the contamination level itself. 
We also provide insights that allow the depth notions considered in this work to be embedded into a unified framework, which we refer to as residual smallness depth. Finally, since the theoretical derivation of maximum bias is carried out under the assumption of known location—due to the intractability of the unknown-location case—we complement our analysis with a numerical study illustrating the effect of contamination on widely used robust estimators of multivariate scatter.

\paragraph{Acknowledgments.}

The authors thank Stanislav Nagy for insightful comments regarding the
numerical study.

\section*{Appendix}

\setcounter{section}{0}
\setcounter{subsection}{0}

\renewcommand{\thesection}{Appendix \Alph{section}}
\renewcommand{\thesubsection}{\Alph{section}\arabic{subsection}}

\renewcommand{\theHsection}{\Alph{section}}
\renewcommand{\theHsubsection}{\Alph{section}\arabic{subsection}}
\section{Proofs}

\subsection{Proofs in Section 4} 

\noindent \textbf{Proof of Lemma 1. } Given $g\left( \mathbf{u%
}\right) =\Phi \left( \sqrt{\mathbf{u}^{t}\mathit{\Gamma} \mathbf{u}}\right) -\Phi
\left( -\sqrt{\mathbf{u}^{t}\mathit{\Gamma} \mathbf{u}}\right) ,$ the Lagrangian is
given by 
\begin{equation*}
	h\left( \mathbf{u,}\lambda \right) =g\left( \mathbf{u}\right) +\lambda
	\left( \mathbf{u}^{t}\mathbf{u-}1\right) .
\end{equation*}%
By differentiating the function $h$ we get 
\begin{eqnarray*}
	\mathbf{0}\mathbf{=}\frac{\partial h\left( \mathbf{u,}\lambda \right) }{%
		\partial \mathbf{u}}=\varphi \left( \sqrt{\mathbf{u}^{t}\mathit{\Gamma} \mathbf{u}}%
	\right) \left[ \frac{\mathit{\Gamma} \mathbf{u}}{\sqrt{\mathbf{u}^{t}\mathit{\Gamma} \mathbf{u}%
	}}\right] &&-\varphi \left( -\sqrt{\mathbf{u}^{t}\mathit{\Gamma} \mathbf{u}}\right) %
	\left[ -\frac{\mathit{\Gamma} \mathbf{u}}{\sqrt{\mathbf{u}^{t}\mathit{\Gamma} \mathbf{u}}}%
	\right] +2\lambda \mathbf{u} \\
	0 &\mathbf{=}&\frac{\partial h\left( \mathbf{u,}\lambda \right) }{\partial
		\lambda }=\mathbf{u}^{t}\mathbf{u-}1.
\end{eqnarray*}%
Therefore, we obtain that%
$\mathbf{0}\mathbf{=}\frac{\partial h\left( \mathbf{u,\lambda }\right) }{%
	\partial \mathbf{u}}=\left[ 2\varphi \left( \sqrt{\mathbf{u}^{t}\mathit{\Gamma} 
	\mathbf{u}}\right) \right] \frac{\mathit{\Gamma} \mathbf{u}}{\sqrt{\mathbf{u}%
		^{t}\mathit{\Gamma} \mathbf{u}}}+2\lambda \mathbf{u}$.

Call $b\left( \mathbf{u}\right) =-\left[ 2\varphi \left( \sqrt{\mathbf{u}%
	^{t}\mathit{\Gamma} \mathbf{u}}\right) \right] /\sqrt{\mathbf{u}^{t}\mathit{\Gamma} \mathbf{u}}%
<0,$ which yields $b\left( \mathbf{u}\right) \mathit{\Gamma} \mathbf{u}\mathbf{=}%
2\lambda \mathbf{u}$ where $\lambda =\frac{1}{2}b\left( \mathbf{u}\right) 
\mathbf{u}^{t}\mathit{\Gamma} \mathbf{u}$ and consequently, $\left( \mathit{\Gamma} \mathbf{-}%
\left[ \mathbf{u}^{t}\mathit{\Gamma} \mathbf{u}\right] I\right) \mathbf{u}\mathbf{=}%
\mathbf{0.}$ Thus, this entails that the critical points of the Lagrangian
are the eigenvectors of $\mathit{\Gamma} .$ Therefore, the function, 
\begin{equation*}
	g\left( \mathbf{u}\right) =\Phi \left( \sqrt{\mathbf{u}^{t}\mathit{\Gamma} \mathbf{u}}%
	\right) -\Phi \left( -\sqrt{\mathbf{u}^{t}\mathit{\Gamma} \mathbf{u}}\right)
\end{equation*}%
has the following Lagrangian critical points we have 
\begin{eqnarray*}
	g\left( \mathbf{v}_{j}\right) &=&\Phi \left( l_{j}^{1/2}\right) -\Phi \left(
	-l_{j}^{1/2}\right) =2\Phi \left( l_{j}^{1/2}\right) -1,\text{ }j=1,...,p \\
	g\left( \mathbf{v}_{p}\right) &\leq &g\left( \mathbf{v}\right) \leq g\left( 
	\mathbf{v}_{1}\right) \text{ for all }\mathbf{v\in }\emph{S}^{p-1} \\
	1-g\left( \mathbf{v}_{j}\right) &=&1-\left[ \Phi \left( l_{j}^{1/2}\right)
	-\Phi \left( -l_{j}^{1/2}\right) \right] =2\left[ 1-\Phi \left(
	l_{j}^{1/2}\right) \right] ,\text{ }j=1,...,p \\
	1-g\left( \mathbf{v}_{p}\right) &\geq &g\left( \mathbf{v}\right) \geq
	1-g\left( \mathbf{v}_{1}\right) \text{ for all }\mathbf{v\in }\emph{S}^{p-1}
\end{eqnarray*}%
and the minimum is given by $m=\min \left( g\left( \mathbf{v}_{p}\right)
,1-g\left( \mathbf{v}_{1}\right) \right) .$ $\hfill \square $

\noindent \textbf{Proof of Corollary 1. } Since $g(\mathbf{v}%
_{p})\leq g(\mathbf{v}_{1})$ and $g$ is an increasing function of $\sqrt{%
	\mathbf{v}^{\top }\mathit{\Gamma} \mathbf{v}}$, $\min \bigl(g(\mathbf{v}_{p}),\,1-g(%
\mathbf{v}_{1})\bigr)$ is maximized when $g(\mathbf{v}_{p})=1-g(\mathbf{v}%
_{1})=0.5$, which implies $\ell _{p}^{1/2}=\ell _{1}^{1/2}=\Phi ^{-1}(3/4)$.
Hence, $\hat{\mathit{\Gamma}}=\left[ \Phi ^{-1}(3/4)\right] ^{2}I_{p},$ and the
corresponding maximum depth equals to $1/2$. $\hfill \square $

\noindent \textbf{Proof of Corollary 2.} For every vector $%
\mathbf{w}\in \mathbb{R}^{p}$, define 
\begin{equation*}
	Z_{n}(\mathbf{w})=\frac{\mathbf{w}^{t}\mathit{\hat{\Gamma}} _{n}^{-1/2}\mathbf{X}}{%
		\left \Vert \mathbf{w}^{t}\mathit{\hat{\Gamma}} _{n}^{-1/2}\right \Vert }\in \mathbb{R}.
\end{equation*}%
Let 
\begin{equation*}
	c_{n}(\mathbf{w})=\frac{1}{\displaystyle \sum_{j=1}^{p}l_{j}^{(n)-1}\left( 
		\mathbf{w}^{t}\mathbf{e}_{j}\right) ^{2}}.
\end{equation*}%
and consider the probabilities 
\begin{eqnarray*}
	a_{\varepsilon ,n}(\mathbf{w}) &=&(1-\varepsilon )\,P_{0}\! \left( Z_{n}(%
	\mathbf{w})^{2}\leq c_{n}(\mathbf{w})\right) +\varepsilon \,P_{n}\! \left(
	Z_{n}(\mathbf{w}_{r})^{2}\leq c_{n}(\mathbf{w})\right) , \\[1ex]
	b_{\varepsilon ,n}(\mathbf{w}) &=&(1-\varepsilon )\,P_{0}\! \left( Z_{n}(%
	\mathbf{w}_{r})^{2}\geq c_{n}(\mathbf{w})\right) +\varepsilon
	\,P_{n}\! \left( Z_{n}(\mathbf{w}_{r})^{2}\geq c_{n}(\mathbf{w})\right) .
\end{eqnarray*}%
Let $\mathbf{w}_{n\text{ }}$ be such that $D\left( \mathit{\hat{\Gamma}} _{n},P_{n}\right)
=\min \left( a_{\varepsilon ,n}\left( \mathbf{w}_{n}\right) ,b_{\varepsilon
	,n}\left( \mathbf{w}_{n}\right) \right) .$ We next analyze two cases.

\textbf{Case i) }$l_{1}^{(n)}$\textbf{$\rightarrow \infty $.}

By definition of the depth, for every
vector $\mathbf{w}$ we have that $$D(\mathit{\hat{\Gamma}} _{n},P_{n})=\min \{a_{\varepsilon ,n}(%
\mathbf{w}_{n}),b_{\varepsilon ,n}(\mathbf{w}_{n})\} \leq \min
\{a_{\varepsilon ,n}(\mathbf{w}),b_{\varepsilon ,n}(\mathbf{w})\}.$$ Since%
\begin{equation*}
	D(\mathit{\hat{\Gamma}} _{n},P_{n})=\min \{a_{\varepsilon ,n}(\mathbf{w}_{n}),b_{%
		\varepsilon ,n}(\mathbf{w}_{n})\} \leq \min \{a_{\varepsilon ,n}(\mathbf{v}%
	_{1}^{(n)}),b_{\varepsilon ,n}(\mathbf{v}_{1}^{(n)})\},
\end{equation*}%
it holds that $a_{\varepsilon ,n}(\mathbf{v}_{1}^{(n)})\geq D(\mathit{\hat{\Gamma}}
_{n},P_{n})\quad b_{\varepsilon ,n}(\mathbf{v}_{1}^{(n)})\geq D(\mathit{\hat{\Gamma}}
_{n},P_{n})$. 

By considering that $\{ \mathbf{v}_{j}^{(n)}\}_{j=1}^{p}$ is an
orthonormal basis, we have that 
\begin{equation*}
	\frac{1}{\sum_{j=1}^{p}l_{j}^{(n)-1}\bigl(\mathbf{v}_{1}^{(n)T}\mathbf{e}_{j}%
		\bigr)^{2}}=l_{1}^{(n)}.
\end{equation*}%
and 
\begin{eqnarray*}
	b_{\varepsilon ,n}\! \left( \mathbf{v}_{1}^{(n)}\right) &=&(1-\varepsilon
	)\,P_{0}\! \left( Z_{n}\! \left( \mathbf{v}_{1}^{(n)}\right) ^{2}\geq
	l_{1}^{(n)}\right) +\varepsilon \,P_{n}\! \left( Z_{n}\! \left( \mathbf{v}%
	_{1}^{(n)}\right) ^{2}\geq l_{1}^{(n)}\right) \\
	&\geq &D(\mathit{\hat{\Gamma}} _{n},P_{n}).
\end{eqnarray*}%
Since $P_{0}\left( Z_{n}\! \left( \mathbf{v}_{1}^{(n)}\right) ^{2}\geq
l_{1}^{(n)}\right) \rightarrow 0$ when $l_{1}^{(n)}\rightarrow \infty $, we
conclude that 

$$\varepsilon \geq \lim_{n\rightarrow \infty }D\left(
\mathit{\hat{\Gamma}} _{n},P_{n}\right) \geq  \lim_{n\rightarrow \infty }D\left(
I,P_{n}\right) \geq \left( 1-\varepsilon \right) D\left( I,P_{0}\right)
=\left( 1-\varepsilon \right) \frac{1}{2},$$
which says that $\varepsilon \geq
1/3.$

\textbf{Case ii) }$l_{p}^{(n)}\rightarrow 0$. With a similar analysis to the
case i) we can conclude that%
\begin{equation*}
	b_{\varepsilon ,n}\left( \mathbf{v}_{p}^{(n)}\right) \geq D\left( \mathit{\hat{\Gamma}}
	_{n},P_{n}\right) \text{ and }a_{\varepsilon ,n}\left( \mathbf{v}%
	_{p}^{(n)}\right) \geq D\left( \mathit{\hat{\Gamma}} _{n},P_{n}\right) .
\end{equation*}%
Since%
\begin{equation*}
	a_{\varepsilon ,n}\left( \mathbf{v}_{p}^{(n)}\right) =(1-\varepsilon
	)\,P_{0}\! \left( Z_{n}\! \left( \mathbf{v}_{p}^{(n)}\right) ^{2}\leq
	l_{p}^{(n)}\right) +\varepsilon \,P_{n}\! \left( Z_{n}\! \left( \mathbf{v}%
	_{p}^{(n)}\right) ^{2}\leq l_{p}^{(n)}\right) \geq D\! \left( \mathit{\hat{\Gamma}}
	_{n},P_{n}\right),
\end{equation*}%
and considering that $P_{0}\left( Z_{n}\! \left( \mathbf{v}_{p}^{(n)}\right)
^{2}\leq l_{p}^{(n)}\right) \rightarrow 0$ when $l_{p}^{(n)}\rightarrow 0$,
we get 
\begin{equation*}
	\varepsilon \geq \lim_{n\rightarrow \infty }D\left( \mathit{\hat{\Gamma}} _{n},P_{n}\right)
	\geq \lim_{n\rightarrow \infty }D\left( I,P_{n}\right) \geq \left(
	1-\varepsilon \right) D\left( I,P_{0}\right) =\left( 1-\varepsilon \right) 
	\frac{1}{2}
\end{equation*}%
which says that $\varepsilon \geq 1/3.$ $\hfill \square $

\noindent \textbf{Proof of Lemma 2. }  If $h^{\mathbf{e}%
}\left( \mathbf{v}\right) <r^{2},$ then $\left( 1-\varepsilon \right)
g\left( \mathbf{v}\right) +\varepsilon \delta \left( h^{\mathbf{e}}\left( 
\mathbf{v}\right) \geq r^{2}\right) \geq \left( 1-\varepsilon \right)
gb_{r}^{\mathbf{e}}$ and $\left( 1-\varepsilon \right) (1-g\left( \mathbf{v}%
\right) )+\varepsilon \delta \left( h^{\mathbf{e}}\left( \mathbf{v}\right)
\leq r^{2}\right) \geq \left( 1-\varepsilon \right) \left( 1-Gb_{r}^{\mathbf{%
		e}}\right) +\varepsilon .$ If\textbf{\ }$h^{\mathbf{e}}\left( \mathbf{v}%
\right) >r^{2},$ we conclude that  $\left( 1-\varepsilon \right) g\left( \mathbf{v}%
\right) +\varepsilon \delta \left( h^{\mathbf{e}}\left( \mathbf{v}\right)
\geq r^{2}\right) \geq \left( 1-\varepsilon \right) ga_{r}^{\mathbf{e}%
}+\varepsilon $ and $\left( 1-\varepsilon \right) (1-g\left( \mathbf{v}%
\right) )+\varepsilon \delta \left( h^{\mathbf{e}}\left( \mathbf{v}\right)
\leq r^{2}\right) \geq \left( 1-\varepsilon \right) \left( 1-Ga_{r}^{\mathbf{%
		e}}\right) $ and the lemma follows. $\hfill \square $

\begin{lemma}
	\label{karushkuhntucker}The critical points of $\Phi \left( \sqrt{\mathbf{v}%
		^{t}\mathit{\Gamma} \mathbf{v}}\right) $ (respectively $1-\Phi \left( \sqrt{\mathbf{v}%
		^{t}\mathit{\Gamma} \mathbf{v}}\right) )$ subject to $\mathbf{v}^{t}\mathit{\Gamma} \mathbf{v-}%
	r^{2}\left( \mathbf{v}^{t}\mathbf{e}\right) ^{2}\geq 0$ (idem $\mathbf{v}%
	^{t}\mathit{\Gamma} \mathbf{v-}r^{2}\left( \mathbf{v}^{t}\mathbf{e}\right) ^{2}\leq
	0) $ are either $\mathbf{v}_{p},$ $\mathbf{v}_{1}$ or occurs at $F_{r}^{%
		\mathbf{e}}.$
\end{lemma}

\noindent \textbf{Proof. }  Take the Lagrangian
\[
\Phi \! \left(\sqrt{\mathbf{v}^t \mathit{\Gamma} \mathbf{v}}\right)
+ \lambda \! \left(\mathbf{v}^t \mathit{\Gamma} \mathbf{v}
- r^2(\mathbf{v}^t \mathbf{e})^2\right)
+ \mathit{\Gamma} (\mathbf{v}^t \mathbf{v} - 1).
\]

The Karush--Kuhn--Tucker conditions say that if a local minimum
occurs at $\tilde{\mathbf v}$, then either
\[
\tilde{\mathbf v}^t \mathit{\Gamma} \tilde{\mathbf v}
- r^2(\tilde{\mathbf v}^t \mathbf e)^2 = 0,
\]
or we have to consider the Lagrangian
\[
\Phi \! \left(\sqrt{\mathbf{v}^t \mathit{\Gamma} \mathbf{v}}\right)
+ \mathit{\Gamma} (\mathbf{v}^t \mathbf{v}-1).
\]

In the last case, it was shown that the
critical points are the eigenvectors of $\mathit{\Gamma} .$ Suppose that $\mathbf{v}%
_{j_{0}\text{,}}$ $j_{0}<p,$is the eigenvector corresponding to the minimum $%
l_{j_{0}}$ such that $\mathbf{v}_{j_{0}\text{ }}^{t}\mathit{\Gamma} \mathbf{v}_{j_{0}%
	\text{ }}>r^{2}\left( \mathbf{v}_{j_{0}\text{ }}^{t}\mathbf{e}\right) ^{2}$,
by the continuity of the functions we could find another $\mathbf{\bar{v}}$
such that $\mathbf{v}_{j_{0}\text{ }}^{t}\mathit{\Gamma} \mathbf{v}_{j_{0}\text{ }}>%
\mathbf{\bar{v}}_{\text{ }}^{t}\mathit{\Gamma} \mathbf{\bar{v}}>r^{2}\left( \mathbf{%
	\bar{v}}^{t}\mathbf{e}\right) ^{2}>r^{2}\left( \mathbf{v}_{j_{0}\text{ }}^{t}%
\mathbf{e}\right) ^{2}$ \ which is a contradiction. Therefore the critical
points occur at $\mathbf{\tilde{v}=v}_{p}$ $\left( \mathbf{v}_{1}\right) $
or $\mathbf{\tilde{v}}^{t}\mathit{\Gamma} \mathbf{\tilde{v}-}r^{2}\left( \mathbf{%
	\tilde{v}}^{t}\mathbf{e}\right) ^{2}=0.$

\noindent \textbf{Proof of Corollary 3.}  Set $q_{j}=\frac{%
	l_{j}^{1/2}}{\mathbf{v}_{j}^{t}\mathbf{e}},$ then $q_{1}=l_{1}^{1/2}$ and $%
q_{j}=\infty ,$ $j=2...,p.$ If $r<l_{1}^{1/2}$ then

\begin{equation*}
	B_{r}^{\mathbf{e}}=\emptyset ,Gb_{r}^{\mathbf{e}}=-\infty \text{ and }%
	gb_{r}^{\mathbf{e}}=\infty .
\end{equation*}%
\begin{equation*}
	\left \{ \mathbf{v}_{1},\mathbf{v}_{2},\dots ,\mathbf{v}_{p}\right \}
	\subseteq A_{r}^{\mathbf{e}}\text{ then }Ga_{r}^{\mathbf{e}}=g\left( \mathbf{%
		v}_{1}\right) \text{ and }ga_{r}^{\mathbf{e}}=g\left( \mathbf{v}_{p}\right) .
\end{equation*}%
Thus,

\begin{equation*}
	D\left( \mathit{\Gamma} ,P_{\varepsilon ,r}\right) =\min \left \{ \left( 1-\varepsilon
	\right) g(\mathbf{v}_{p})+\varepsilon ,\left( 1-\varepsilon \right) \left(
	1-g\left( \mathbf{v}_{1}\right) \right) \right \} .
\end{equation*}%
If $r>l_{1}^{1/2}$ then 
\begin{eqnarray*}
	\left \{ \mathbf{v}_{2},...,\mathbf{v}_{p}\right \} &\subseteq &\left( B_{r}^{%
		\mathbf{e}}\right) ^{c},\mathbf{v}_{1}\in B_{r}^{\mathbf{e}},\text{ }Gb_{r}^{%
		\mathbf{e}}=g\left( \mathbf{v}_{1}\right) , \\
	\left \{ \mathbf{v}_{2},...,\mathbf{v}_{p}\right \} &\subseteq &A_{r}^{\mathbf{%
			e}},\text{ }\mathbf{v}_{1}\notin A_{r}^{\mathbf{e}},ga_{r}^{\mathbf{e}%
	}=g\left( \mathbf{v}_{p}\right) .
\end{eqnarray*}%
We need to determine $gb_{r}^{\mathbf{e}}$ and $Ga_{r}^{\mathbf{e}}.$ Then
the function $g\left( \mathbf{v}\right) $ has a minimum in $B_{r}^{\mathbf{e}%
}\cup F_{r}^{\mathbf{e}}$ (a compact set) and the minimum has to occur in $%
F_{r}^{\mathbf{e}}.$ Likewise, the function $g\left( \mathbf{v}\right) $ has
a maximum in $A_{r}^{\mathbf{e}}\cup F_{r}^{\mathbf{e}}$ which must occur at 
$F_{r}^{\mathbf{e}}.$ In both cases, Lemma \ref{karushkuhntucker} supports
that statement. We use spherical coordinates,%
\begin{eqnarray*}
	\mathbf{v}^{t}\mathbf{v}_{1} &=&\cos \eta _{0},\text{ } \\
	\mathbf{v}^{t}\mathbf{v}_{2} &=&\sin \eta _{0}\cos \eta _{1}, \\
	\mathbf{v}^{t}\mathbf{v}_{3} &=&\sin \eta _{0}\sin \eta _{1}\cos \eta _{2},%
	\text{ }0\leq \eta _{0},..\eta _{p-3}\leq \pi \\
	&&\vdots \\
	\mathbf{v}^{t}\mathbf{v}_{p-1} &=&\sin \eta _{0}\sin \eta _{1}\sin \eta
	_{2}\dots \cos \eta _{p-2},0\leq \eta _{p-2}\leq 2\pi \\
	\mathbf{v}^{t}\mathbf{v}_{p} &=&\sin \eta _{0}\sin \eta _{1}\sin \eta
	_{2}\dots \sin \eta _{p-2},
\end{eqnarray*}%
($\eta _{j}$ is the angle between $\mathbf{v}_{j+1}$ and $\mathbf{v,}$ $%
j=0,...,p-1\mathbf{).}$ Then%
\begin{equation*}
	\mathbf{v}^{t}\mathit{\Gamma} \mathbf{v=}r^{2}\left( \mathbf{v}^{t}\mathbf{v}%
	_{1}\right) ^{2}
\end{equation*}%
is equivalent to 
\begin{equation*}
	\begin{aligned}
		l_{1}\cos ^{2}\eta _{0}
		&+\sin ^{2}\eta _{0}\Big[
		l_{2}\cos ^{2}\eta _{1}
		+l_{3}\sin ^{2}\eta _{1}\cos ^{2}\eta _{2} \\
		&\qquad
		+\dots
		+l_{p}\sin ^{2}\eta _{1}\sin ^{2}\eta _{2}\dots
		\sin ^{2}\eta _{p-2}
		\Big]  \\
		&= r^{2}\cos ^{2}\eta _{0}
	\end{aligned}
\end{equation*}

\begin{equation*}
	tg^{2}\eta _{0}=\frac{r^{2}-l_{1}}{d}
\end{equation*}
with
\[
d = l_{2}\cos ^{2}\eta _{1}+\sin ^{2}\eta _{1}\bigl(\dots
(l_{p-2}\cos ^{2}\eta _{p-3}
+\sin ^{2}\eta _{p-3}(l_{p-1}\cos ^{2}\eta _{p-2}
+l_{p}\sin ^{2}\eta _{p-2}))\dots \bigr).
\]

Observe that the denominator is  a nested convex combination. Thus, we can
ensure that $\frac{r^{2}-l_{1}}{l_{2}} \leq tg^{2}\eta _{0}\leq \frac{r^{2}-l_{1}}{l_{p}}$ if and only if either $\sqrt{\frac{r^{2}-l_{1}}{l_{2}}} \leq tg\left( \eta
_{0}\right) \leq \sqrt{\frac{r^{2}-l_{1}}{l_{p}}}$ or $-\sqrt{\frac{%
		r^{2}-l_{1}}{l_{p}}}\leq tg\left( \eta _{0}\right) \leq -\sqrt{\frac{%
		r^{2}-l_{1}}{l_{2}}}$ if and only if either $0 \leq \arctan \left( \sqrt{\frac{r^{2}-l_{1}}{l_{2}}}%
\right) \leq \eta _{0}\leq \arctan \left( \sqrt{\frac{r^{2}-l_{1}}{l_{p}}}%
\right) \leq \pi /2$ or $\text{or }\pi /2 \leq \arctan \left( -\sqrt{\frac{r^{2}-l_{1}}{l_{p}}}%
\right) \leq \eta _{0}\leq \arctan \left( -\sqrt{\frac{r^{2}-l_{1}}{l_{2}}}%
\right) \leq \pi.$
Then, by similarity of triangles we can say that, for either $l=l_{2}$ or $%
l_{p},$ 
\begin{equation*}
	\frac{\sqrt{1+\frac{r^{2}-l_{1}}{l}}}{1}=\frac{1}{\cos \left( \boldsymbol{\theta} \right) 
	},\text{ \  \  \ }\frac{\sqrt{1+\frac{r^{2}-l_{1}}{l}}}{1}=\frac{\sqrt{\frac{%
				r^{2}-l_{1}}{l}}}{\sin \left( \boldsymbol{\theta} \right) }.
\end{equation*}%
Therefore, we obtain that 
\begin{eqnarray*}
	\max_{A_{r}^{\mathbf{e}}\cup F_{r}^{\mathbf{e}}}\mathbf{v}^{t}\mathit{\Gamma} \mathbf{%
		v} &\mathbf{=}&\max_{F_{r}^{\mathbf{e}}}\mathbf{v}^{t}\mathit{\Gamma} \mathbf{v=}%
	r^{2}\cos ^{2}\left[ \arctan \left( \sqrt{\frac{r^{2}-l_{1}}{l_{2}}}\right) %
	\right] =r^{2}\frac{l_{2}}{r^{2}+l_{2}-l_{1}}, \\
	\mathbf{v}_{M,r} &=&\arg \max_{F_{r}^{\mathbf{e}}}\mathbf{v}^{t}\mathit{\Gamma} 
	\mathbf{v} \\
	&\mathbf{=}&\left( \cos \arctan \left( \sqrt{\frac{r^{2}-l_{1}}{l_{2}}}%
	\right) ,\sin \arctan \left( \sqrt{\frac{r^{2}-l_{1}}{l_{2}}}\right)
	,0,\dots ,0)\right) \\
	&=&\left( \sqrt{\frac{l_{2}}{r^{2}+l_{2}-l_{1}}},\sqrt{\frac{r^{2}-l_{1}}{%
			r^{2}+l_{2}-l_{1}}},0,\dots ,0\right), \\
	\min_{B_{r}^{\mathbf{e}}\cup F_{r}^{\mathbf{e}}}\mathbf{v}^{t}\mathit{\Gamma} \mathbf{%
		v} &\mathbf{=}&r^{2}\cos ^{2}\left[ \arctan \left( \sqrt{\frac{r^{2}-l_{1}}{%
			l_{p}}}\right) \right] =r^{2}\frac{l_{p}}{r^{2}+l_{p}-l_{1}},\\
	\mathbf{v}_{m,r} &=&\arg \min_{B_{r}^{\mathbf{e}}\cup F_{r}^{\mathbf{e}}}%
	\mathbf{v}^{t}\mathit{\Gamma} \mathbf{v} \\
	&\mathbf{=}&\left( \cos \arctan \left( \sqrt{\frac{r^{2}-l_{1}}{l_{p}}}%
	\right) ,0,\dots ,0,\sin \arctan \left( \sqrt{\frac{r^{2}-l_{1}}{l_{p}}}%
	\right) \right) \\
	&=&\left( \sqrt{\frac{l_{p}}{r^{2}+l_{p}-l_{1}}},0,\dots ,0,\sqrt{\frac{%
			r^{2}-l_{1}}{r^{2}+l_{p}-l_{1}}}\right),
\end{eqnarray*}%
which means that%
\begin{eqnarray*}
	gb_{r}^{\mathbf{e}} &=&g\left( \mathbf{v}_{m,r}\right) =2\Phi \left( r\cos %
	\left[ \arctan \left( \sqrt{\frac{r^{2}-l_{1}}{l_{p}}}\right) \right]
	\right) -1, \\
	Ga_{r}^{\mathbf{e}} &=&g\left( \mathbf{v}_{M,r}\right) =2\Phi \left( r\cos %
	\left[ \arctan \left( \sqrt{\frac{r^{2}-l_{1}}{l_{2}}}\right) \right]
	\right) -1.
\end{eqnarray*}%
Since $\mathbf{v}_{M,r}$ is a maximum over $A_{r}^{\mathbf{e}}\cup F_{r}^{%
	\mathbf{e}},$ $g\left( \mathbf{v}_{2}\right) \leq g\left( \mathbf{v}%
_{M,r}\right) \leq g\left( \mathbf{v}_{1}\right) .$ Then, we have that 
\begin{equation}
	D\left( \mathit{\Gamma} ,P_{\varepsilon ,r}\right) =\left \{ 
	\begin{array}{cc}
		\min \left \{ 
		\begin{array}{c}
			\left( 1-\varepsilon \right) g(\mathbf{v}_{p})+\varepsilon , \\ 
			\left( 1-\varepsilon \right) \left( 1-g\left( \mathbf{v}_{1}\right) \right)%
		\end{array}%
		\right \} & \text{if }r\leq l_{1}^{1/2} \\ 
		\min \left \{ 
		\begin{array}{c}
			\left( 1-\varepsilon \right) \left( 1-g\left( \mathbf{v}_{1}\right) \right)
			+\varepsilon ,\left( 1-\varepsilon \right) g\left( \mathbf{v}_{p}\right)
			+\varepsilon , \\ 
			\left( 1-\varepsilon \right) g\left( \mathbf{v}_{m,r}\right) ,\left(
			1-\varepsilon )\left( 1-g\left( \mathbf{v}_{M,r}\right) \right) \right)%
		\end{array}%
		\right \} & \text{if }r>l_{1}^{1/2}%
	\end{array}%
	\right. .  \label{maxdepth1}
\end{equation}%
$\hfill \square $

\noindent \textbf{Proof of Lemma 3. } If $r<l_{1}^{1/2}$
then the deepest matrix should verify 
\begin{eqnarray*}
	\left( 1-\varepsilon \right) g\left( \mathbf{v}_{p}\right) +\varepsilon
	&=&\left( 1-\varepsilon \right) \left( 1-g\left( \mathbf{v}_{1}\right)
	\right) \\
	g\left( \mathbf{v}_{p}\right) +g\left( \mathbf{v}_{1}\right) &=&\frac{%
		1-2\varepsilon }{1-\varepsilon }.
\end{eqnarray*}%
To increase the depth, $g\left( \mathbf{v}_{p}\right) $ should increase and $%
g\left( \mathbf{v}_{1}\right) $ should decrease, and then they would reach 
\begin{equation*}
	2g\left( \mathbf{v}_{p}\right) =2g\left( \mathbf{v}_{1}\right) =\frac{%
		1-2\varepsilon }{1-\varepsilon }.
\end{equation*}%
This implies that 
\begin{eqnarray*}
	g\left( \mathbf{v}_{p}\right) &=&g\left( \mathbf{v}_{1}\right) =\frac{%
		1-2\varepsilon }{2(1-\varepsilon )} \\
	l_{p}^{1/2} &=&l_{1}^{1/2}=\Phi ^{-1}\left( \frac{3-4\varepsilon }{%
		4(1-\varepsilon )}\right)
\end{eqnarray*}%
and the deepest matrix is a multiple of the identity matrix.

Let us take $l_{1}^{1/2}<r.$ Thus we get%
\begin{equation*}
	D\left( \mathit{\Gamma} ,P_{\varepsilon ,r}\right) =\min \left \{ 
	\begin{array}{c}
		\left( 1-\varepsilon \right) \left( 1-g\left( \mathbf{v}_{1}\right) \right)
		+\varepsilon ,\left( 1-\varepsilon \right) g\left( \mathbf{v}_{p}\right)
		+\varepsilon , \\ 
		\left( 1-\varepsilon \right) g\left( \mathbf{v}_{m,r}\right) ,\left(
		1-\varepsilon )\left( 1-g\left( \mathbf{v}_{M,r}\right) \right) \right)%
	\end{array}%
	\right \} .
\end{equation*}%
with $g(\mathbf{v}_{p})\leq ...\leq g(\mathbf{v}_{2})\leq g\left( \mathbf{v}%
_{M,r}\right) \leq g\left( \mathbf{v}_{1}\right) .$ The best estimator
should be selected so that $g\left( \mathbf{v}_{m,r}\right) =g\left( \mathbf{%
	v}_{M,r}\right) $ (by making $l_{2}=\dots =l_{p},$ $g\left( \mathbf{v}%
_{m,r}\right) $ increases) and $\mathbf{v}_{m,r}=\mathbf{v}_{M,r}.$ Call $%
\mathbf{v}_{0,r}=\mathbf{v}_{m,r}$ and $\left( 1-\varepsilon \right) g\left( 
\mathbf{v}_{0,r}\right) =(1-\varepsilon )\left( 1-g\left( \mathbf{v}%
_{0,r}\right) \right) ,$ which entails that $g\left( \mathbf{v}_{0,r}\right)
=1/2$ or $\sqrt{\mathbf{v}_{0,r}^{t}\mathit{\Gamma} \mathbf{v}_{0,r}}=\mathbf{\Phi }%
^{-1}\left( 3/4\right) .$ $\ $Then we take matrices $\mathit{\Gamma} $ such that $%
\Phi ^{-1}\left( 3/4\right) =0.6745\leq l_{1}^{1/2}$ $.$ All the four
quantities involved in $D\left( \mathit{\Gamma} ,P_{\varepsilon ,r}\right) $ must be
equal in order to get the best estimator In this case we get that 
\begin{eqnarray*}
	\left( 1-\varepsilon \right) \left( 1-g_{1}\right) +\varepsilon &=&\frac{%
		\left( 1-\varepsilon \right) }{2}=\left( 1-\varepsilon \right)
	g_{p}+\varepsilon \\
	\left( 1-g_{1}\right) +\frac{\varepsilon }{1-\varepsilon } &=&\frac{1}{2}%
	\text{ \ };\text{ \ }g_{p}+\frac{\varepsilon }{1-\varepsilon }=\frac{1}{2} \\
	g_{1} &=&\frac{1}{2}+\frac{\varepsilon }{1-\varepsilon }\text{ \  \ ; \  \ }%
	g_{p}=\frac{1}{2}-\frac{\varepsilon }{1-\varepsilon } \\
	g_{1} &=&\frac{1+\varepsilon }{2\left( 1-\varepsilon \right) }\text{ \  \ ; \
		\ }g_{p}=\frac{1-2\varepsilon }{2\left( 1-\varepsilon \right) },
\end{eqnarray*}%
which implies that 
\begin{eqnarray*}
	g_{1} &=&2\Phi \left( l_{1}^{1/2}\right) -1=\frac{1}{2}+\frac{\varepsilon }{%
		1-\varepsilon } \\
	\Phi \left( l_{1}^{1/2}\right) &=&\frac{3}{4}+\frac{\varepsilon }{2\left(
		1-\varepsilon \right) }=\frac{3-3\varepsilon +2\varepsilon }{4\left(
		1-\varepsilon \right) }=\frac{3-\varepsilon }{4\left( 1-\varepsilon \right) }%
	. \\
	g_{p} &=&2\Phi \left( l_{p}^{1/2}\right) -1=\frac{1}{2}-\frac{\varepsilon }{%
		1-\varepsilon } \\
	\Phi \left( l_{p}^{1/2}\right) &=&\frac{3}{4}-\frac{\varepsilon }{2\left(
		1-\varepsilon \right) }=\frac{3-3\varepsilon -2\varepsilon }{4\left(
		1-\varepsilon \right) }=\frac{3-5\varepsilon }{4\left( 1-\varepsilon \right) 
	}.
\end{eqnarray*}%
$\hfill \square $

\noindent \textbf{Proof of Lemma 4. } We have to analyze,

\begin{equation*}
	\min \left \{ 
	\begin{array}{c}
		\left( 1-\varepsilon \right) P_{0}\left( \frac{\mathbf{w}_{r}^{t}\mathit{\Gamma} _{r}%
		}{\left \Vert \mathbf{w}_{r}^{t}\mathit{\Gamma} _{r}^{-1/2}\right \Vert }\mathbf{XX}^{t}%
		\frac{\mathit{\Gamma} _{r}^{-1/2}\mathbf{w}_{r}}{\left \Vert \mathbf{w}_{r}^{t}\mathit{\Gamma}
			_{r}^{-1/2}\right \Vert }\leq \frac{1}{\sum_{j=2}^{p}l_{j}^{^{(r)}-1}\left( 
			\mathbf{w}_{r}^{t}\mathbf{e}_{j}\right) ^{2}}\right) + \\ 
		\varepsilon \delta \left( r^{2}l_{1}^{^{(r)}-1/2}\mathbf{w}_{r}^{t}\mathbf{ee%
		}^{t}\mathbf{w}_{r}\leq 1\right) , \\ 
		\left( 1-\varepsilon \right) P_{0}\left( \frac{\mathbf{w}_{r}^{t}\mathit{\Gamma}
			_{r}^{-1/2}}{\left \Vert \mathbf{w}_{r}^{t}\mathit{\Gamma} _{r}^{-1/2}\right \Vert }%
		\mathbf{XX}^{t}\frac{\mathit{\Gamma} _{r}^{-1/2}\mathbf{w}_{r}}{\left \Vert \mathbf{w}%
			_{r}^{t}\mathit{\Gamma} _{r}^{-1/2}\right \Vert }\geq \frac{1}{%
			\sum_{j=2}^{p}l_{j}^{^{(r)}-1}\left( \mathbf{w}_{r}^{t}\mathbf{e}_{j}\right)
			^{2}}\right) + \\ 
		\varepsilon \delta \left( r^{2}l_{1}^{^{(r)}-1/2}\mathbf{w}_{r}^{t}\mathbf{ee%
		}^{t}\mathbf{w}_{r}\geq 1\right)%
	\end{array}%
	\right \} ,
\end{equation*}%
where $\mathbf{w}_{r}$ is a vector which yields the minimum.

(i) If $\  \mathbf{w}_{r}$ $\in L\left( \mathbf{e}_{2},\dots ,\mathbf{e}%
_{p}\right) $ we get for $r$ sufficiently large, $\min \left \{ \varepsilon
,\left( 1-\varepsilon \right) \right \} .$

(ii) If $\mathbf{w}_{r}$ $\in \left[ L\left( \mathbf{e}_{2},\dots ,\mathbf{e}%
_{p}\right) \right] ^{c},$ we get for $r$ sufficiently large, $\min \left \{
\varepsilon ,\left( 1-\varepsilon \right) \right \} .$ Then, for $r$ large
enough, since $\varepsilon \leq 1/2,$ 
\begin{equation}
	D\left( \mathit{\Gamma} _{r},P_{\varepsilon ,r}\right) =\varepsilon ,  \label{depth2}
\end{equation}%
and we conclude the statement of the lemma. $\hfill \square $

\noindent \textbf{Proof of Lemma 5. } Take $\mathit{\Gamma}
=\sum_{j=1}^{p}l_{j}\mathbf{v}_{j}\mathbf{v}_{j}^{t}$. The depth of a matrix
is given by, if $r>l_{1}^{1/2},$ 
\begin{equation*}
	\min \left \{ \left( 1-\varepsilon \right) \left( 1-Gb_{r}^{\mathbf{e}%
	}\right) +\varepsilon ,\left( 1-\varepsilon \right) ga_{r}^{\mathbf{e}%
	}+\varepsilon ,\left( 1-\varepsilon \right) gb_{r}^{\mathbf{e}},\left(
	1-\varepsilon )\left( 1-Ga_{r}^{\mathbf{e}}\right) \right) \right \} \text{.}
\end{equation*}%
Call $g_{m}=\min_{\mathbf{v\in }F_{r}^{\mathbf{e}}}g\left( \mathbf{v}\right) 
$ and $g_{M}=\arg \max_{\mathbf{v\in }F_{r}^{\mathbf{e}}}g\left( \mathbf{v}%
\right) $. The different configurations that a matrix $\mathit{\Gamma} $ might have
are given in the following table depending on where $\mathbf{v}_{1}$ and $%
\mathbf{v}_{p}$ belong to,%
\begin{equation*}
	\begin{tabular}{|l|l|l|l|l|l|l|l|l|l|}
		\hline
		$A_{r}^{\mathbf{e}}$ & $\mathbf{v}_{1},\mathbf{v}_{p}$ &  &  & $\mathbf{v}%
		_{1}$ & $\mathbf{v}_{p}$ & $\mathbf{v}_{1}$ & $\mathbf{v}_{p}$ &  &  \\ 
		\hline
		$B_{r}^{\mathbf{e}}$ &  & $\mathbf{v}_{1},\mathbf{v}_{p}$ &  & $\mathbf{v}%
		_{p}$ & $\mathbf{v}_{1}$ &  &  & $\mathbf{v}_{p}$ & $\mathbf{v}_{1}$ \\ 
		\hline
		$F_{r}^{\mathbf{e}}$ &  &  & $\mathbf{v}_{1},\mathbf{v}_{p}$ &  &  & $%
		\mathbf{v}_{p}$ & $\mathbf{v}_{1}$ & $\mathbf{v}_{1}$ & $\mathbf{v}_{p}$ \\ 
		\hline
	\end{tabular}%
\end{equation*}

(i) $\mathbf{v}_{1},\mathbf{v}_{p}\in A_{r}^{\mathbf{e}}$: the depth of $%
\mathit{\Gamma} $ should be 
\begin{equation*}
	\min \left( 1-\varepsilon \right) \left( 1-g_{M}\right) +\varepsilon ,\left(
	1-\varepsilon \right) g_{p}+\varepsilon ,\left( 1-\varepsilon \right)
	g_{m},\left( 1-\varepsilon )\left( 1-g_{1}\right) \right)
\end{equation*}

(ii) $\mathbf{v}_{1},\mathbf{v}_{p}\in B_{r}^{\mathbf{e}}:$%
\begin{equation*}
	\min \left \{ \left( 1-\varepsilon \right) \left( 1-g_{1}\right)
	+\varepsilon ,\left( 1-\varepsilon \right) g_{m}+\varepsilon ,\left(
	1-\varepsilon \right) g_{p},\left( 1-\varepsilon )\left( 1-g_{M}\right)
	\right) \right \}
\end{equation*}

(iii) $\mathbf{v}_{1},\mathbf{v}_{p}\in F_{r}^{\mathbf{e}}:$%
\begin{equation*}
	\min \left \{ \left( 1-\varepsilon \right) \left( 1-g_{1}\right)
	+\varepsilon ,\left( 1-\varepsilon \right) g_{p}+\varepsilon ,\left(
	1-\varepsilon \right) g_{p},\left( 1-\varepsilon )\left( 1-g_{1}\right)
	\right) \right \}
\end{equation*}

(iv) $\mathbf{v}_{1}\in A_{r}^{\mathbf{e}}$, $\mathbf{v}_{p}\in B_{r}^{%
	\mathbf{e}},$ 
\begin{equation*}
	\min \left \{ \left( 1-\varepsilon \right) \left( 1-g_{M}\right)
	+\varepsilon ,\left( 1-\varepsilon \right) g_{m}+\varepsilon ,\left(
	1-\varepsilon \right) g_{p},\left( 1-\varepsilon )\left( 1-g_{1}\right)
	\right) \right \}
\end{equation*}

(v) $\mathbf{v}_{p}\in A_{r}^{\mathbf{e}}$, $\mathbf{v}_{1}\in B_{r}^{%
	\mathbf{e}},$ 
\begin{equation*}
	\min \left \{ \left( 1-\varepsilon \right) \left( 1-g_{1}\right)
	+\varepsilon ,\left( 1-\varepsilon \right) g_{p}+\varepsilon ,\left(
	1-\varepsilon \right) g_{m},\left( 1-\varepsilon )\left( 1-g_{M}\right)
	\right) \right \}
\end{equation*}

(vi) $\mathbf{v}_{1}\in A_{r}^{\mathbf{e}}$, $\mathbf{v}_{p}\in F_{r}^{%
	\mathbf{e}},$ 
\begin{equation*}
	\min \left \{ \left( 1-\varepsilon \right) \left( 1-g_{M}\right)
	+\varepsilon ,\left( 1-\varepsilon \right) g_{p}+\varepsilon ,\left(
	1-\varepsilon \right) g_{p},\left( 1-\varepsilon )\left( 1-g_{1}\right)
	\right) \right \}
\end{equation*}

(vii) $\mathbf{v}_{1}\in B_{r}^{\mathbf{e}}$, $\mathbf{v}_{p}\in F_{r}^{%
	\mathbf{e}},$ 
\begin{equation*}
	\min \left \{ \left( 1-\varepsilon \right) \left( 1-g_{1}\right)
	+\varepsilon ,\left( 1-\varepsilon \right) g_{p}+\varepsilon ,\left(
	1-\varepsilon \right) g_{p},\left( 1-\varepsilon )\left( 1-g_{M}\right)
	\right) \right \}
\end{equation*}

(viii) $\mathbf{v}_{1}\in A_{r}^{\mathbf{e}}$, $\mathbf{v}_{p}\in F_{r}^{%
	\mathbf{e}}$ 
\begin{equation*}
	\min \left \{ \left( 1-\varepsilon \right) \left( 1-g_{M}\right)
	+\varepsilon ,\left( 1-\varepsilon \right) g_{p}+\varepsilon ,\left(
	1-\varepsilon \right) g_{p},\left( 1-\varepsilon )\left( 1-g_{1}\right)
	\right) \right \}
\end{equation*}

(ix) $\mathbf{v}_{1},\mathbf{v}_{p}\in F_{r}^{\mathbf{e}}$ 
\begin{equation*}
	\min \left \{ \left( 1-\varepsilon \right) \left( 1-g_{1}\right)
	+\varepsilon ,\left( 1-\varepsilon \right) g_{p}+\varepsilon ,\left(
	1-\varepsilon \right) g_{p},\left( 1-\varepsilon )\left( 1-g_{1}\right)
	\right) \right \}
\end{equation*}%
In any case, we have that $D\left( \mathit{\Gamma} ,P_{\varepsilon ,r}\right) \leq
\min \left \{ \left( 1-\varepsilon \right) g_{s},\left( 1-\varepsilon
\right) \left( 1-g_{t}\right) \right \} ,s\in \left \{ p,m\right \} ,$ $t\in
\left \{ 1,M\right \} ,$ $g_{s}\leq 1-g_{t}$ and the statement follows. $%
\hfill \square $

\noindent \textbf{Proof of Theorem 1. }(i) Lemma above
says that if the estimators moves to the boundary of the domain of the
eigenvalues, then $\varepsilon \geq 1/3$ \ and $\varepsilon ^{\ast }\geq
1/3. $

(ii) If the depth estimators $\left \{ \hat{\mathit{\Gamma}}_{r}\right \} $ remain
bounded in the $\varepsilon $-neighborhood, we have that for point mass
contaminations we have that $D\left( \mathit{\Gamma} _{r},P_{\varepsilon ,r}\right)
\leq D\left( \hat{\mathit{\Gamma}}_{r},P_{\varepsilon ,r}\right) $ for $r$ large
enough. By the preceding lemmas, we have that $\varepsilon \leq \left(
1-\varepsilon \right) /2,$ and $\varepsilon \leq 1/3.$ Then, if the depth
estimator remains bounded in the $\varepsilon $-contamination neighborhood
then the level of contamination $\varepsilon $ is less than $1/3.$
Equivalently, if $\varepsilon >1/3$ then the depth estimator becomes
unbounded in the $\varepsilon $-contamination neighborhood. Therefore $%
\varepsilon ^{\ast }\leq 1/3.$ $\hfill \square $

\noindent \textbf{Proof of Lemma 6. } Since $P_{0}=\Phi $
then $D_{M}\left( \Phi \right) =1/2$ and as  in Remark 1 of the paper let $\beta=\left[ \Phi ^{-1}\left( \frac{3}{4}\right)\right] ^{2}$.

Observe that $\delta \geq 0$ since given $Q$ any distribution on $\mathbb{R}^{p},$

\begin{eqnarray*}
	P_{\varepsilon ,Q}\left( \left \vert \mathbf{u}^{t}\mathbf{X}\right \vert
	^{2}\leq \mathbf{u}^{t}\mathit{\Gamma} \mathbf{u}\right) &\geq &\left( 1-\varepsilon
	\right) P_{0}\left( \left \vert \mathbf{u}^{t}\mathbf{X}\right \vert
	^{2}\leq \mathbf{u}^{t}\mathit{\Gamma} \mathbf{u}\right) \text{ \ },\text{ \  \ } \\
	P_{\varepsilon ,Q}\left( \left \vert \mathbf{u}^{t}\mathbf{X}\right \vert
	^{2}\geq \mathbf{u}^{t}\mathit{\Gamma} \mathbf{u}\right) &\geq &\left( 1-\varepsilon
	\right) P_{0}\left( \left \vert \mathbf{u}^{t}\mathbf{X}\right \vert
	^{2}\geq \mathbf{u}^{t}\mathit{\Gamma} \mathbf{u}\right)
\end{eqnarray*}

\begin{eqnarray*}
	&&\min \left \{ P_{\varepsilon ,Q}\left( \left \vert \mathbf{u}^{t}\mathbf{X}%
	\right \vert ^{2}\leq \mathbf{u}^{t}\mathit{\Gamma} \mathbf{u}\right) ,P_{\varepsilon
		,Q}\left( \left \vert \mathbf{u}^{t}\mathbf{X}\right \vert ^{2}\geq \mathbf{u%
	}^{t}\mathit{\Gamma} \mathbf{u}\right) \right \} \\
	&\geq &\left( 1-\varepsilon \right) \min \left \{ P_{0}\left( \left \vert 
	\mathbf{u}^{t}\mathbf{X}\right \vert ^{2}\leq \mathbf{u}^{t}\mathit{\Gamma} \mathbf{u}%
	\right) ,P_{0}\left( \left \vert \mathbf{u}^{t}\mathbf{X}\right \vert
	^{2}\geq \mathbf{u}^{t}\mathit{\Gamma} \mathbf{u}\right) \right \} \\
	&\geq &\left( 1-\varepsilon \right) D\left( \mathit{\Gamma} ,P_{0}\right) \\
	D\left( \mathit{\Gamma} ,P_{\varepsilon ,Q}\right) &\geq &\left( 1-\varepsilon
	\right) D\left( \mathit{\Gamma} ,P_{0}\right) \\
	D\left( \mathit{\Gamma} ,P_{\varepsilon ,Q}\right) &=&\inf_{\mathbf{u\in }\mathcal{S}%
		^{p-1}}\min \left \{ P_{\varepsilon ,Q}\left( \left \vert \mathbf{u}^{t}%
	\mathbf{X}\right \vert ^{2}\leq \mathbf{u}^{t}\mathit{\Gamma} \mathbf{u}\right)
	,P_{\varepsilon ,Q}\left( \left \vert \mathbf{u}^{t}\mathbf{X}\right \vert
	^{2}\geq \mathbf{u}^{t}\mathit{\Gamma} \mathbf{u}\right) \right \}
\end{eqnarray*}%
\begin{eqnarray*}
	D_{M}\left( P_{\varepsilon ,Q}\right) &\geq &D\left( \sqrt{\beta}I,P_{\varepsilon
		,Q}\right) \geq \left( 1-\varepsilon \right) D_{M}(P_{0}) \\
	\Lambda \left( \varepsilon ,P_{0}\right) &=&\inf_{Q}D_{M}\left(
	P_{\varepsilon ,Q}\right) \geq \left( 1-\varepsilon \right) D_{M}(P_{0})
\end{eqnarray*}%
Since 
\begin{eqnarray*}
	\min \left \{ 
	\begin{array}{c}
		P_{\varepsilon ,Q}\left( \left \vert \mathbf{u}^{t}\mathbf{X}\right \vert
		^{2}\leq \mathbf{u}^{t}\mathit{\Gamma} \mathbf{u}\right) , \\ 
		P_{\varepsilon ,Q}\left( \left \vert \mathbf{u}^{t}\mathbf{X}\right \vert
		^{2}\geq \mathbf{u}^{t}\mathit{\Gamma} \mathbf{u}\right)%
	\end{array}%
	\right \} &\leq &P_{\varepsilon ,Q}\left( \left \vert \mathbf{u}^{t}\mathbf{X%
	}\right \vert ^{2}\leq \mathbf{u}^{t}\mathit{\Gamma} \mathbf{u}\right) \\
	&\leq &\left( 1-\varepsilon \right) P_{0}\left( \left \vert \mathbf{u}^{t}%
	\mathbf{X}\right \vert ^{2}\leq \mathbf{u}^{t}\mathit{\Gamma} \mathbf{u}\right)
	+\varepsilon \\
	\min \left \{ 
	\begin{array}{c}
		P_{\varepsilon ,Q}\left( \left \vert \mathbf{u}^{t}\mathbf{X}\right \vert
		^{2}\leq \mathbf{u}^{t}\mathit{\Gamma} \mathbf{u}\right) , \\ 
		P_{\varepsilon ,Q}\left( \left \vert \mathbf{u}^{t}\mathbf{X}\right \vert
		^{2}\geq \mathbf{u}^{t}\mathit{\Gamma} \mathbf{u}\right)%
	\end{array}%
	\right \} &\leq &P_{\varepsilon ,Q}\left( \left \vert \mathbf{u}^{t}\mathbf{X%
	}\right \vert ^{2}\geq \mathbf{u}^{t}\mathit{\Gamma} \mathbf{u}\right) \\
	&\leq &\left( 1-\varepsilon \right) P_{0}\left( \left \vert \mathbf{u}^{t}%
	\mathbf{X}\right \vert ^{2}\geq \mathbf{u}^{t}\mathit{\Gamma} \mathbf{u}\right)
	+\varepsilon
\end{eqnarray*}%
Then,%
\begin{eqnarray*}
	&&\min \left \{ P_{\varepsilon ,Q}\left( \left \vert \mathbf{u}^{t}\mathbf{X}%
	\right \vert ^{2}\leq \mathbf{u}^{t}\mathit{\Gamma} \mathbf{u}\right) ,P_{\varepsilon
		,Q}\left( \left \vert \mathbf{u}^{t}\mathbf{X}\right \vert ^{2}\geq \mathbf{u%
	}^{t}\mathit{\Gamma} \mathbf{u}\right) \right \} \\
	&\leq &\left( 1-\varepsilon \right) \min \left \{ P_{0}\left( \left \vert 
	\mathbf{u}^{t}\mathbf{X}\right \vert ^{2}\leq \mathbf{u}^{t}\mathit{\Gamma} \mathbf{u}%
	\right) ,P_{0}\left( \left \vert \mathbf{u}^{t}\mathbf{X}\right \vert
	^{2}\geq \mathbf{u}^{t}\mathit{\Gamma} \mathbf{u}\right) \right \} +\varepsilon
\end{eqnarray*}%
and%
\begin{equation*}
	D\left( \mathit{\Gamma} ,P_{\varepsilon ,Q}\right) \leq \left( 1-\varepsilon \right)
	D\left( \mathit{\Gamma} ,P_{0}\right) +\varepsilon .
\end{equation*}%
Consequently, (i) follows.

(ii) Set $\alpha =\frac{e}{1-\varepsilon }-\delta \left( \varepsilon
,P_{0}\right) .$ Then $0\leq \alpha <1/2.$ \ If $\mathit{\Gamma} \notin L\left(
\alpha .P_{0}\right) $ we have that $D\left( \mathit{\Gamma} ,P_{0}\right)
<D_{M}\left( P_{0}\right) -\alpha .$ Therefore, 
\begin{eqnarray*}
	D\left( \mathit{\Gamma} ,P_{\varepsilon ,Q}\right) &\leq &\left( 1-\varepsilon
	\right) D\left( \mathit{\Gamma} ,P_{0}\right) +\varepsilon <\left( 1-\varepsilon
	\right) \left[ D_{M}\left( P_{0}\right) -\frac{\varepsilon }{1-\varepsilon}%
	+\delta \left( \varepsilon ,P_{0}\right) \right] +\varepsilon \\
	&\leq &\left( 1-\varepsilon \right) \left[ D_{M}\left( P_{0}\right) -\frac{%
		\varepsilon }{1-\varepsilon }+\frac{\Lambda \left( \varepsilon ,P_{0}\right)
		-\left( 1-\varepsilon \right) D_{M}\left( P_{0}\right) }{1-\varepsilon }%
	\right] +\varepsilon =\Lambda \left( \varepsilon ,P_{0}\right)
\end{eqnarray*}%
and (ii) follows. 

(iii)   $L\left( \alpha ,P_{0}\right) =L\left( \varepsilon /\left(
1-\varepsilon \right) -\delta ,P_{0}\right) \subset L\left( \varepsilon
/\left( 1-\varepsilon \right) ,P_{0}\right) $ Then, $\mathit{\Gamma} \in \left(
L\left( \alpha ,P_{0}\right) \right) ^{c}$ implies, by \ (ii) that $\mathit{\Gamma}
\in \left( M\left( P_{\varepsilon ,Q}\right) \right) ^{c}$ or $\left(
L\left( \alpha ,P_{0}\right) \right) ^{c}\subset \left( M\left(
P_{\varepsilon ,Q}\right) \right) ^{c}$ which says that.  for all
distribution $Q$,  $M\left(
P_{\varepsilon ,Q}\right) \subset L\left( \alpha ,P_{0}\right) \subset
L\left( \varepsilon /\left( 1-\varepsilon \right) ,P_{0}\right) $.

Since $D\left( \mathit{\Gamma} ,P_{0}\right) =\min \left( g\left( 
\mathbf{v}_{p}\right) ,1-g\left( \mathbf{v}_{1}\right) \right) $ and $%
D_{M}\left( P_{0}\right) =1/2,$ we have that%
\begin{eqnarray*}
	L\left( \varepsilon /\left( 1-\varepsilon \right) ,P_{0}\right) &=&\left \{
	\mathit{\Gamma} \succeq 0:\min \left( g\left( \mathbf{v}_{p}\right) ,1-g\left( 
	\mathbf{v}_{1}\right) \right) \geq \frac{1}{2}-\varepsilon /\left(
	1-\varepsilon \right) \right \} \\
	&=&\left \{ \mathit{\Gamma} \succeq 0:1-g\left( \mathbf{v}_{1}\right) \geq g\left( 
	\mathbf{v}_{p}\right) \geq \frac{1}{2}-\varepsilon /\left( 1-\varepsilon
	\right) \right \} \\
	&&\cup \left \{ \mathit{\Gamma} \succeq 0:g\left( \mathbf{v}_{p}\right) \geq 1-g\left( 
	\mathbf{v}_{1}\right) \geq \frac{1}{2}-\varepsilon /\left( 1-\varepsilon
	\right) \right \} \\
	&=&\left \{ \mathit{\Gamma} \succeq 0:\frac{1}{2}+\varepsilon /\left( 1-\varepsilon
	\right) \geq g\left( \mathbf{v}_{1}\right) \geq g\left( \mathbf{v}%
	_{p}\right) \geq \frac{1}{2}-\varepsilon /\left( 1-\varepsilon \right)
	\right \} .
\end{eqnarray*}%
Then, $\left \Vert \mathit{\Gamma} \right \Vert _{op}\leq \Phi ^{-1}\left( \frac{%
	3-\varepsilon }{4\left( 1-\varepsilon \right) }\right) ,$ $\left \Vert \mathit{\Gamma}
^{-1}\right \Vert _{op}\leq 1/\Phi ^{-1}\left( \frac{3-5\varepsilon }{4\left(
	1-\varepsilon \right) }\right) $ if $\mathit{\Gamma} \in L\left( \varepsilon /\left(
1-\varepsilon \right) ,P\right) $ and 
\begin{eqnarray*}
	\left \Vert L\left( \varepsilon /\left( 1-\varepsilon \right) ,P_{0}\right)
	\right \Vert &=&\sup_{\mathit{\Gamma} \in L\left( \varepsilon /\left( 1-\varepsilon
		\right) ,P_{0}\right) }\max \left \{ \frac{\left \Vert \mathit{\Gamma} \right \Vert _{op}%
	}{\sqrt{\beta} },\frac{\sqrt{\beta} }{\left \Vert \mathit{\Gamma} ^{-1}\right \Vert _{op}}\right \} \\
	&=&\max \left \{ \frac{1}{\sqrt{\beta} }\Phi ^{-1}\left( \frac{3-\varepsilon }{%
		4\left( 1-\varepsilon \right) }\right) ,\frac{\sqrt{\beta} }{\Phi ^{-1}\left( \frac{%
			3-5\varepsilon }{4\left( 1-\varepsilon \right) }\right) }\right \}
\end{eqnarray*}%
The bias functions turn out to be 
\begin{eqnarray*}
	b_{S}\left( \hat{\mathit{\Gamma}},\varepsilon ,P\right) &=&\lambda _{(1)}\left( \hat{%
		\mathit{\Gamma}}\left( P_{0}\right) ^{-1/2}\hat{\mathit{\Gamma}}\left( P\right) \hat{\mathit{\Gamma}}%
	\left( P_{0}\right) ^{-1/2}\right) =\sqrt{\beta} ^{-1}\lambda _{(1)}\left( \hat{%
		\mathit{\Gamma}}\left( P\right) \right) =\sqrt{\beta} ^{-1}\left \Vert \hat{\mathit{\Gamma}}\left(
	P\right) \right \Vert _{op} \\
	b_{I}\left( \hat{\mathit{\Gamma}},\varepsilon ,P\right) &=&\lambda _{(1)}\left( \hat{%
		\mathit{\Gamma}}\left( P_{0}\right) ^{1/2}\hat{\mathit{\Gamma}}^{-1}\left( P\right) \hat{\mathit{\Gamma}%
	}\left( P_{0}\right) ^{1/2}\right) =\sqrt{\beta} \lambda _{(1)}\left( \hat{\mathit{\Gamma}}%
	^{-1}\left( P\right) \right) =\sqrt{\beta} \left \Vert \hat{\mathit{\Gamma}}^{-1}\left(
	P\right) \right \Vert _{op}.
\end{eqnarray*}%
\begin{equation*}
	B\left( \hat{\mathit{\Gamma}},\varepsilon ,P_{0}\right) =\max \left \{ \sqrt{\beta}
	^{-1}\sup_{P\in \mathcal{P}_{\varepsilon }}\left \Vert \hat{\mathit{\Gamma}}\left(
	P\right) \right \Vert _{op},\sqrt{\beta} \sup_{P\in \mathcal{P}_{\varepsilon
	}}\left \Vert \hat{\mathit{\Gamma}}^{-1}\left( P\right) \right \Vert _{op}\right \}
\end{equation*}%
Since $M\left( P_{\varepsilon ,Q}\right) \subset L\left( \varepsilon /\left(
1-\varepsilon \right) ,P_{0}\right) $ then $B\left( \hat{\mathit{\Gamma}},\varepsilon
,P_{0}\right) \leq \left \Vert L\left( \varepsilon /\left( 1-\varepsilon
\right) ,P_{0}\right) \right \Vert .$ $\hfill \square $

\subsection{Proofs in Section 5.}


\noindent \textbf{Proof of Lemma 7. }\cite{ChenGaoRen2018},
in the proof of Theorem 2.1, showed that, for $\mathbf{\boldsymbol{\theta} }=\mathbf{0}$ 
\textbf{\  \ }%
\begin{eqnarray}
	\Phi \left( \left \Vert \hat{\boldsymbol{\theta}}_T\right \Vert \right) &\leq &\frac{1}{2}+%
	\frac{\varepsilon }{1-\varepsilon }+40\sqrt{\frac{6e\pi }{1-e^{-1}}}\sqrt{%
		\frac{p+1}{n}}+\frac{7}{2}\sqrt{\frac{\log \left( 1/\delta \right) }{n}}
	\label{cgr2.1} \\
	&=&\frac{1+\varepsilon }{2\left( 1-\varepsilon \right) }+40\sqrt{\frac{6e\pi 
		}{1-e^{-1}}}\sqrt{\frac{p+1}{n}}+\frac{7}{2}\sqrt{\frac{\log \left( 1/\delta
			\right) }{n}}  \notag
\end{eqnarray}%
The upper bound becomes useless for $\varepsilon \geq 1/3$ since it is
greater than 1. Therefore, we need to consider $\varepsilon <1/3-c.$ Put $%
b\left( p,n\right) =40\sqrt{\frac{6e\pi }{1-e^{-1}}}\sqrt{\frac{p+1}{n}}+%
\frac{7}{2}\sqrt{\frac{\log \left( 1/\delta \right) }{n}}$ and $a\left(
\varepsilon \right) =\frac{1}{2}+\frac{\varepsilon }{1-\varepsilon }=\frac{%
	1+\varepsilon }{2\left( 1-\varepsilon \right) }.$ Thus, (\ref{cgr2.1}) is
equivalent to 
\begin{eqnarray*}
	\left \Vert \hat{ \boldsymbol{\theta}_{T}}\right \Vert &\leq &\Phi ^{-1}\left( \frac{1}{2}+\frac{%
		\varepsilon }{1-\varepsilon }+b\left( p,n\right) \right) -\Phi ^{-1}\left( 
	\frac{1}{2}\right) +\Phi ^{-1}\left( \frac{1}{2}\right) \\
	&=&\frac{1}{\phi \left( \Phi ^{-1}\left( \mathit{\Gamma} _{C}\right) \right) }\left( 
	\frac{\varepsilon }{1-\varepsilon }+b\left( p,n\right) \right) \leq C\left( 
	\sqrt{\frac{p}{n}}\vee \varepsilon +\sqrt{\frac{\log \left( 1/\delta \right) 
		}{n}}\right) ,\text{ } \\
	\text{with }\mathit{\Gamma} _{C} &\in &\left( \frac{1}{2},\frac{1}{2}+\frac{%
		\varepsilon }{1-\varepsilon }+b\left( p,n\right) \right) .
\end{eqnarray*}%
Therefore, consider $\Phi ^{-1}\left( a\left( \varepsilon \right) \right) $
rather than $\Phi ^{-1}\left( \frac{1}{2}\right) $ to center $\Phi
^{-1}\left( \frac{1}{2}+\frac{\varepsilon }{1-\varepsilon }+b\left(
p,n\right) \right) $ and we get 
\begin{eqnarray}
	\left \Vert \hat{ \boldsymbol{\theta}_{T}}\right \Vert &\leq &\Phi ^{-1}\left( a\left(
	\varepsilon \right) +b\left( p,n\right) \right) -\Phi ^{-1}\left( a\left(
	\varepsilon \right) \right) +\Phi ^{-1}\left( a\left( \varepsilon \right)
	\right)  \label{(cotamedtukey)} \\
	&=&\frac{1}{\phi \left( \Phi ^{-1}\left( \eta _{B}\right) \right) }b\left(
	p,n\right) +\Phi ^{-1}\left( a\left( \varepsilon \right) \right) ,\text{ \ }%
	\eta _{B}\in \left( a\left( \varepsilon \right) ,a\left( \varepsilon \right)
	+b\left( p,n\right) \right)  \notag \\
	&\leq &\frac{1}{\phi \left( \Phi ^{-1}\left( a\left( \varepsilon \right)
		+b\left( p,n\right) \right) \right) }b\left( p,n\right) +\Phi ^{-1}\left(
	a\left( \varepsilon \right) \right)  \notag
\end{eqnarray}%
If $\varepsilon <1/3-c$, where $c$ is a positive constant, $a\left(
\varepsilon \right) $ is increasing on $[0,1/3-c)$, let $a(m_{c})$ be the
maximum value. Take $d$ such that $1-a\left( m_{c}\right) >d>0$ and $\left(
p,n\right) $ $\in $ $A_{c,d}$ $=$ $\left \{ \left( p,n\right) :b\left(
p,n\right) <1-a\left( m_{c}\right) -d\right \}$. Consequently, there exists
a constant 
\begin{equation*}
	C_{c,d}=\sup_{\left( p,n\right) \in A_{c,d}}\left[ \phi \left( \Phi
	^{-1}\left( a\left( m_{c}\right) +b\left( p,n\right) \right) \right) \right]
	^{-1}
\end{equation*}%
such that 
\begin{equation*}
	\left \Vert \hat{ \boldsymbol{\theta}_{T}}\right \Vert \leq \tilde{C}_{c,d}\left( \sqrt{\frac{p}{%
			n}}\vee B_{L}\left( \mathbf{\hat{ \boldsymbol{\theta}}}_{T},\varepsilon ,\Phi \right) +%
	\sqrt{\frac{\log \left( 1/\delta \right) }{n}}\right) ,
\end{equation*}%
and the lemma follows. $\hfill \square $

\noindent \textbf{Proof of Lemma 8. } 
From \cite{ChenGaoRen2018} (p. 1955) with probability greater than $1-2\delta$, the following inequality holds 
\begin{equation}
	\sup_{\mathbf{u\in }\emph{S}^{p-1}}\left \vert \Phi \left( \sqrt{\beta }%
	\right) -\Phi \left( 	\sqrt{\beta }x_{\mathbf{u}}\right) \right \vert \leq \frac{\varepsilon }{%
		2(1-\varepsilon )}+a_{n,\delta },
	\label{cheninequality}
\end{equation}%
where $a_{n,\delta }=40\sqrt{\frac{6e\pi }{1-e^{-1}}}\sqrt{\frac{3+2p}{n%
}}+\frac{7}{4}\sqrt{\frac{\log (1/\delta )}{n}}$ and $x_{\mathbf{u}}=\sqrt{%
	\frac{\mathbf{u}^{T}\left( \hat{\mathit{\Gamma}}/\beta \right) u}{\mathbf{u}%
		^{T}\Sigma \mathbf{u}}}$.

Note now that 
\begin{eqnarray}
	&&\sup_{\mathbf{u}\in S^{p-1}}\left \vert \Phi (\sqrt{\beta })-\Phi \left( 
	\sqrt{\beta }x_{\mathbf{u}}\right) \right \vert  \label{f1} \\
	&=&\max \left( \Phi \left( \sqrt{\beta }\sup_{\mathbf{u}\in S^{p-1}}x_{%
		\mathbf{u}}\right) -\Phi (\sqrt{\beta }),\Phi (\sqrt{\beta })-\Phi \left( 
	\sqrt{\beta }\inf_{\mathbf{u}\in S^{p-1}}x_{\mathbf{u}}\right) \right)
\end{eqnarray}%
If $\inf_{\mathbf{u\in }S^{p-1}}x_{\mathbf{u}}\leq \sup_{\mathbf{u}\in
	S^{p-1}}x_{\mathbf{u}}\leq 1$, then 
\begin{equation*}
	\sup_{\mathbf{u}\in S^{p-1}}\left \vert \Phi (\sqrt{\beta })-\Phi \left( 
	\sqrt{\beta }x_{\mathbf{u}}\right) \right \vert =\Phi (\sqrt{\beta })-\Phi
	\left( \sqrt{\beta }\inf_{\mathbf{u}\in S^{p-1}}x_{\mathbf{u}}\right) .
\end{equation*}
If $1\leq \inf_{\mathbf{u\in }S^{p-1}}x_{\mathbf{u}}\leq \sup_{\mathbf{u}\in
	S^{p-1}}x_{\mathbf{u}},$ then 
\begin{equation*}
	\sup_{\mathbf{u}\in S^{p-1}}\left \vert \Phi (\sqrt{\beta })-\Phi \left( 
	\sqrt{\beta }x_{\mathbf{u}}\right) \right \vert =\Phi \left( \sqrt{\beta }%
	\sup_{\mathbf{u}\in S^{p-1}}x_{\mathbf{u}}\right) -\Phi (\sqrt{\beta }).
\end{equation*}

If $\inf_{\mathbf{u\in }S^{p-1}}x_{\mathbf{u}}\leq 1\leq \sup_{\mathbf{u}\in
	S^{p-1}}x_{\mathbf{u}}$ we have to analyse both cases. If the maximum on the
right-hand side of (\ref{f1}) occurs in $\Phi \left( \sqrt{\beta }\sup_{%
	\mathbf{u}\in S^{p-1}}x_{\mathbf{u}}\right) -\Phi (\sqrt{\beta })$ then $%
\sup_{\mathbf{u}\in S^{p-1}}x_{\mathbf{u}}\geq 1,$ and we have, 
\begin{equation*}
	\Phi \left( \sqrt{\beta }\inf_{\mathbf{u}\in S^{p-1}}x_{\mathbf{u}}\right)
	\leq \Phi (\sqrt{\beta })+\frac{\varepsilon }{2\left( 1-\varepsilon \right) }%
	+\frac{1}{2}a_{n,\delta }.
\end{equation*}%
If we denote $a\left( \varepsilon \right) =\frac{3}{4}+\frac{\varepsilon }{%
	2\left( 1-\varepsilon \right) }=\frac{3-\varepsilon }{4\left( 1-\varepsilon
	\right) }$ 
then 
\begin{eqnarray*}
	\sup_{\mathbf{u}\in S^{p-1}}x_{\mathbf{u}} &\leq &\frac{1}{\sqrt{\beta }}\Phi ^{-1}\left(
	a\left( \varepsilon \right) +\frac{1}{2}a_{n,\delta }\right) -\frac{1}{\sqrt{%
			\beta }}\Phi ^{-1}\left( a\left( \varepsilon \right) \right) +1 \\
	&&+\left[ \frac{1}{\sqrt{\beta }}\Phi ^{-1}\left( a\left( \varepsilon
	\right) \right) -1\right]
\end{eqnarray*}%

Since $a(\varepsilon) > \tfrac{1}{2}$, there exists $\eta \in \left[a(\varepsilon),\, a(\varepsilon)+\tfrac{1}{2}a_{n,\delta}\right]$ such that
\begin{equation*}
	\sup_{\mathbf{u}\in S^{p-1}} x_{\mathbf{u}} - 1
	\le
	\frac{1}{2\sqrt{\beta}}
	\frac{1}{\phi\!\left(\Phi^{-1}(\eta)\right)}\, a_{n,\delta}
	+
	\left[\frac{1}{\sqrt{\beta}}\Phi^{-1}\!\left(a(\varepsilon)\right)-1\right].
\end{equation*}
Moreover, since $\phi(\Phi^{-1}(t))$ is decreasing for $t > \tfrac{1}{2}$, we obtain
\begin{equation*}
	0 \le \sup_{\mathbf{u}\in S^{p-1}} x_{\mathbf{u}} - 1
	\le
	\frac{1}{2\sqrt{\beta}}
	\frac{1}{\phi\!\left(\Phi^{-1}\!\left(a(\varepsilon)+\tfrac{1}{2}a_{n,\delta}\right)\right)}\, a_{n,\delta}
	+
	\left[\frac{1}{\sqrt{\beta}}\Phi^{-1}\!\left(a(\varepsilon)\right)-1\right]
	= g_S(\varepsilon,n,\delta,p).
\end{equation*}

Call 
\begin{equation*}
	B_{E}\left( \varepsilon \right) =\left[ \frac{1}{\sqrt{\beta }}\Phi
	^{-1}\left( a\left( \varepsilon \right) \right) -1\right] .\text{ }
\end{equation*}

Denote $b\left(
\varepsilon \right) =\frac{3}{4}-\frac{\varepsilon }{2\left( 1-\varepsilon
	\right) }=\frac{3-5\varepsilon }{4\left( 1-\varepsilon \right) }$ and note that $\Phi(\sqrt{\beta})=3/4$. If we assume now that the maximum in (\ref{f1}) occurs at $\Phi (\sqrt{\beta })-\Phi \left( 
\sqrt{\beta }\inf_{\mathbf{u}\in S^{p-1}}x_{\mathbf{u}}\right) ,$ and then $%
\inf_{\mathbf{u}\in S^{p-1}}x_{\mathbf{u}}\leq 1$, we get that 
\begin{eqnarray*}
	0 &\leq &\Phi (\sqrt{\beta })-\Phi \left( \sqrt{\beta }\inf_{u\in
		S^{p-1}}x_{u}\right) \leq \frac{\varepsilon }{2\left( 1-\varepsilon \right) }%
	+\frac{1}{2}a_{n,\delta } \\
	\Phi \left( \sqrt{\beta }\inf_{u\in S^{p-1}}x_{u}\right) &\geq &\Phi (\sqrt{%
		\beta })-\frac{\varepsilon }{2\left( 1-\varepsilon \right) }-\frac{1}{2}%
	a_{n,\delta } \\
	\inf_{u\in S^{p-1}}x_{u} &\geq &\frac{1}{\sqrt{\beta }}\Phi ^{-1}\left(
	b\left( \varepsilon \right) \right) -1+1-\frac{1}{\sqrt{\beta }}\Phi
	^{-1}\left( b\left( \varepsilon \right) \right) \\
	&&+\frac{1}{\sqrt{\beta }}\Phi ^{-1}\left( b\left( \varepsilon \right) -%
	\frac{1}{2}a_{n,\delta }\right)
\end{eqnarray*}%

Then, there exists $\eta \in \left( b\left( \varepsilon \right) -\frac{1}{2}%
a_{n,\delta },b\left( \varepsilon \right) \right) $ such that 
\begin{eqnarray*}
	\inf_{\mathbf{u}\in S^{p-1}}x_{\mathbf{u}}-1 &\geq &-\frac{1}{2\sqrt{\beta }}%
	\frac{1}{\phi \left( \Phi ^{-1}\left( \eta \right) \right) }a_{n,\delta }+%
	\frac{1}{\sqrt{\beta }}\Phi ^{-1}\left( b\left( \varepsilon \right) \right)
	-1,\text{ } \\
	&\geq &-\frac{1}{2\sqrt{\beta }}\frac{1}{\phi \left( \Phi ^{-1}\left(
		b\left( \varepsilon \right) \right) \right) }a_{n,\delta }+\frac{1}{\sqrt{%
			\beta }}\Phi ^{-1}\left( b\left( \varepsilon \right) \right) -1.
\end{eqnarray*}%
Hence,%
\begin{equation*}
	0\leq 1-\inf_{\mathbf{u}\in S^{p-1}}x_{\mathbf{u}}\leq \frac{1}{2\sqrt{\beta 
	}}\frac{1}{\phi \left( \Phi ^{-1}\left( b\left( \varepsilon \right) \right)
		\right) }a_{n,\delta }+1-\frac{1}{\sqrt{\beta }}\Phi ^{-1}\left( b\left(
	\varepsilon \right) \right) =g_{I}\left( \varepsilon ,n,\delta ,p\right) .
\end{equation*}%
where we have used that  $b(\varepsilon)>1 / 2$ and 
$\phi\left(\Phi^{-1}(t)\right)$ is decreasing in $t>1 / 2$.

Call 
\begin{equation*}
	B_{I}\left( \varepsilon \right) =1-\frac{1}{\sqrt{\beta }}\Phi ^{-1}\left( 
	\frac{3-5\varepsilon }{4\left( 1-\varepsilon \right) }\right) .
\end{equation*}%

Since $a(\varepsilon) \geq b(\varepsilon)$ and $\Phi^{-1}$ is increasing, we have $\Phi^{-1}(a(\varepsilon)) \geq \Phi^{-1}(b(\varepsilon))$, which implies $B_{I}(\varepsilon) \leq B_{E}(\varepsilon)$. Observe that there is an $\eta \in \left[ \beta \inf_{%
	\mathbf{u}\in S^{p-1}}x_{\mathbf{u}},\beta \right] \cup \left[ \beta ,\beta
\sup_{\mathbf{u\in }S^{p-1}}x_{\mathbf{u}}\right] $ such that 
\begin{eqnarray*}
	\sup_{\mathbf{u}\in S^{p-1}}\left \vert \Phi (\sqrt{\beta })-\Phi \left( 
	\sqrt{\beta x_{\mathbf{u}}}\right) \right \vert &=&\sup_{\mathbf{u}\in
		S^{p-1}}\left \vert \left( \beta -\beta x_{\mathbf{u}}\right) \frac{\varphi
		\left( \eta \right) }{\sqrt{\eta }}\right \vert \\
	&\geq &\frac{\varphi \left( 1+g_{S}\left( \varepsilon ,n,p,\delta \right)
		\right) }{\sqrt{1+g_{S}\left( \varepsilon ,n,p,\delta \right) }}\sup_{%
		\mathbf{u}\in S^{p-1}}\left \vert \left( \beta -\beta x_{\mathbf{u}}\right)
	\right \vert .
\end{eqnarray*}%
Then, we can conclude that%
\begin{equation*}
	\inf_{\Sigma \in \mathcal{F}\left( M\right) ,P\in \mathcal{P}_{\varepsilon
		}\left( P_{0}\right) }P\left( \left \Vert \hat{\Sigma}-\Sigma \right \Vert
	_{op}^{2}\leq C^{\ast }\left( \max \left \{ \frac{p}{n},B_{E}^{2}(\varepsilon
	)\right \} +\frac{\log \left( 1/\delta \right) }{n}\right) \right) \geq
	\alpha .
\end{equation*}

$\hfill \square $

\noindent \textbf{Proof of Lemma 9. }Since we have
regression through the origin model (a zero intercept model) and the
covariables have an ellipsoidal distribution and we are dealing with affine
equivariant estimators we can assume that $B=0,$ $\Sigma =I$ and $\sigma =1$%
, Then, $P_{B},$ the joint distribution of $\left( \mathbf{X},\mathbf{Y}%
\right) \in 
\mathbb{R}
^{p}\times 
\mathbb{R}
^{m},$ is given by%
\begin{equation*}
	\left( 
	\begin{array}{c}
		\mathbf{Y} \\ 
		\mathbf{X}%
	\end{array}%
	\right) \sim N_{p+m}\left( \left( 
	\begin{array}{c}
		0_{m\times 1} \\ 
		0_{p\times 1}%
	\end{array}%
	\right) ,\left( 
	\begin{array}{cc}
		I_{m\times m} & 0_{m\times p} \\ 
		0_{p\times m} & I_{p\times p}%
	\end{array}%
	\right) \right)
\end{equation*}%
Then, if $A,U\in 
\mathbb{R}
^{p\times m},$ it holds that 
\begin{equation*}
	\left( 
	\begin{tabular}{l}
		$\mathbf{Y}-A^{t}\mathbf{X}$ \\ 
		$U^{t}\mathbf{X}$%
	\end{tabular}%
	\right) \sim N_{2m}\left( \left( 
	\begin{array}{c}
		0_{m\times 1} \\ 
		0_{m\times 1}%
	\end{array}%
	\right) ,\left( 
	\begin{array}{cc}
		I_{m\times m}+B^{t}B & -B^{t}U \\ 
		-U^{t}B & U^{t}U%
	\end{array}%
	\right) \right) .
\end{equation*}%
Given $D\left( A,P_{B}\right) =\inf_{U\in 
	\mathbb{R}
	^{p\times m}}P_{B}\left( (U^{t}X,Y-A^{t}X\right) \geq 0),$ from the proof of
Theorem 4.1, p.1167 of \cite{Gao2020}, we have the following inequality

\begin{equation}
	D\left( \hat{A},P_{B}\right) \geq \frac{1}{2}-\frac{\varepsilon }{%
		1-\varepsilon }-C_{1}\left( \sqrt{\frac{mp+p}{n}}+\sqrt{\frac{\log \left(
			1/\delta \right) }{2n}}\right) ,  \label{bounddepthMR}
\end{equation}%
with $C_{1}$ an absolute constant (independent of $n,p,m,\varepsilon )$).

Note that $X^{t} U Y \mid X \sim N\left(0,\left \Vert X^{t}U\right \Vert^{2}\right)$, then we have that 

\begin{eqnarray*}
	P_{B}\left( (U^{t}X,Y-A^{t}X\right)\geq 0) &=& 
	E_{X}\left[ P\left( \left. X^{t}UA^{t}X\leq X^{t}UY\right \vert X\right) %
	\right] \text{ \ } \\
	&=&E\left[ 1-\Phi \left( \frac{X^{t}UA^{t}X}{\left \Vert X^{t}U\right \Vert }%
	\right) \right] =1-E\Phi \left( \left( \frac{U^{t}X}{\left \Vert
		U^{t}X\right \Vert },A^{t}X\right) \right) .
\end{eqnarray*}%
Therefore, 	
\begin{eqnarray*}
	\inf_{U\in 
		\mathbb{R}
		^{p\times m}}P_{B}\left( (U^{t}X,Y-A^{t}X\right) &\geq &0)=1-\sup_{U\in 
		\mathbb{R}
		^{p\times m}}E\Phi \left( \left( \frac{U^{t}X}{\left \Vert U^{t}X\right \Vert }%
	,A^{t}X\right) \right) \\
	&=&1-E\Phi \left( \left \Vert A^{t}X\right \Vert \right) .
\end{eqnarray*}%
Take the spectral decomposition of the non-negative matrix $AA^{t}\in 
\mathbb{R}
^{p\times p},$ $AA^{t}=PDP^{t}$ $P\in 
\mathbb{R}
^{p\times p}$ orthonormal and $D=diag\left( \lambda _{1},...,\lambda
_{p}\right) $ with $0\leq \lambda _{0}\leq \lambda _{1}\leq \dots \leq
\lambda _{p}$, then $Z=PX\sim N_{p}\left( 0,I\right) .$ Therefore, 
\begin{equation*}
	E\Phi \left( \sqrt{X^{t}AA^{t}X}\right) =E\Phi \left( \sqrt{Z^{t}DZ}\right)
	=E\Phi \left( t\sqrt{Z^{t}\frac{D}{tr\left( D\right) }Z}\right) =E\Phi
	\left( t\sqrt{\sum_{i=1}^{p}\boldsymbol{\theta} _{i}Z_{i}^{2}}\right) ,
\end{equation*}%
with $t=\sqrt{tr\left( D\right) },$ $\left \{ \boldsymbol{\theta} _{i}\right \}
_{i=1}^{p}\subset \left[ 0,1\right] $ and $\sum_{i=1}^{p}\boldsymbol{\theta} _{i}=1$. The
square root is a concave function, and we can obtain that 
\begin{equation*}
	\sqrt{\sum_{i=1}^{p}\boldsymbol{\theta} _{i}Z_{i}^{2}}\geq \sum_{i=1}^{p}\boldsymbol{\theta}
	_{i}\left \vert Z_{i}\right \vert .
\end{equation*}%
$\Phi $ is also a concave function on $[0,\infty )$, then take the functions 
$h:[0,\infty )\times \left[ 0,1\right] ^{p}\rightarrow \left[ 0.5,1\right] $
and $g:[0,\infty )\rightarrow \left[ 0.5,1\right] $ as $h\left( t,\boldsymbol{\theta}
_{1},\dots ,\boldsymbol{\theta} _{p}\right) =\Phi \left( t\sqrt{\sum_{i=1}^{p}\boldsymbol{\theta}
	_{i}Z_{i}^{2}}\right) $ and $g\left( t\right) =h\left( t,1,0,...,0\right) ,$ 
\begin{eqnarray*}
	\Phi \left( t\sqrt{\sum_{i=1}^{p}\boldsymbol{\theta} _{i}Z_{i}^{2}}\right) &\geq
	&\sum_{i=1}^{p}\boldsymbol{\theta} _{i}\Phi \left( t\left \vert Z_{i}\right \vert \right) \\
	E\Phi \left( t\sqrt{\sum_{i=1}^{p}\boldsymbol{\theta} _{i}Z_{i}^{2}}\right) &\geq
	&\sum_{i=1}^{p}\boldsymbol{\theta} _{i}E\Phi \left( t\left \vert Z_{i}\right \vert \right)
	=E\Phi \left( t\left \vert Z_{1}\right \vert \right) =g\left( t\right) ,
\end{eqnarray*}%
since $E\Phi \left( t\left \vert Z_{1}\right \vert \right) =\dots =E\Phi
\left( t\left \vert Z_{p}\right \vert \right) .$ Then%
\begin{equation*}
	\inf_{U\in 
		\mathbb{R}
		^{p\times m}}P_{B}\left( (U^{t}X,Y-A^{t}X\right) \geq 0)=1-E\Phi \left( t%
	\sqrt{\sum_{i=1}^{p}\boldsymbol{\theta} _{i}Z_{i}^{2}}\right) \leq 1-E\Phi \left( t\sqrt{%
		Z_{1}^{2}}\right) =1-g\left( t\right).
\end{equation*}%
Then, from (\ref{bounddepthMR}) we can say, 
\begin{eqnarray*}
	D\left( \hat{A},P_{B}\right) &=&1-g\left( \hat{t}\right) \geq \frac{%
		1-3\varepsilon }{2\left( 1-\varepsilon \right) }-C_{1}\left( \sqrt{\frac{mp+p%
		}{n}}+\sqrt{\frac{\log \left( 1/\delta \right) }{2n}}\right) \\
	g\left( \hat{t}\right) &\leq &\frac{1+\varepsilon }{2\left( 1-\varepsilon
		\right) }+C_{1}\left( \sqrt{\frac{mp+p}{n}}+\sqrt{\frac{\log \left( 1/\delta
			\right) }{2n}}\right) .
\end{eqnarray*}%
Since $g$\ is increasing, we get that 
\begin{eqnarray*}
	\hat{t} &\leq &g^{-1}\left( \frac{1+\varepsilon }{2\left( 1-\varepsilon
		\right) }+C_{1}\left( \sqrt{\frac{pm+p}{n}}+\sqrt{\frac{\log \left( 1/\delta
			\right) }{2n}}\right) \right) \\
	&=&g^{-1}\left( \frac{1+\varepsilon }{2\left( 1-\varepsilon \right) }%
	+C_{1}\left( \sqrt{\frac{pm+p}{n}}+\sqrt{\frac{\log \left( 1/\delta \right) 
		}{2n}}\right) \right) -g^{-1}\left( \frac{1+\varepsilon }{2\left(
		1-\varepsilon \right) }\right) \\
	&& +g^{-1}\left( \frac{1+\varepsilon }{2\left(
		1-\varepsilon \right) }\right) \\
	&=&C_{1}\left( g^{-1}\right) ^{\prime }\left( \xi \right) \left( \sqrt{\frac{%
			pm+p}{n}}+\sqrt{\frac{\log \left( 1/\delta \right) }{2n}}\right)
	+g^{-1}\left( \frac{1+\varepsilon }{2\left( 1-\varepsilon \right) }\right) ,%
\end{eqnarray*}
where $	\xi \in \left( \frac{1+\varepsilon }{2\left( 1-\varepsilon \right) },\frac{%
	1+\varepsilon }{2\left( 1-\varepsilon \right) }+C_{1}\left( \sqrt{\frac{pm+p%
	}{n}}+\sqrt{\frac{\log \left( 1/\delta \right) }{2n}}\right) \right).$

Thus,%
\begin{equation*}
	\hat{t}\leq C\left[ \left( \frac{pm}{n}\vee g^{-1}\left( \frac{1+\varepsilon 
	}{2\left( 1-\varepsilon \right) }\right) \right) +\sqrt{\frac{\log \left(
			1/\delta \right) }{2n}}\right]
\end{equation*}%
with high probability.$\hfill \square $

\begin{remark}
	Observe that $g^{-1}\left( \frac{1+\varepsilon }{2\left( 1-\varepsilon
		\right) }\right) $ does not depend either on the number of dependent
	variables $m$ or the number of independent variables $p.$
\end{remark}

\noindent \textbf{Proof of Lemma 10. } Since we
are dealing with affine equivariant estimators we can assume that $%
\boldsymbol{\beta }=0,$ $\Sigma =I$ and $\sigma =1$. With this background we
get that%
\begin{equation*}
	\left( 
	\begin{array}{c}
		y \\ 
		\mathbf{X}%
	\end{array}%
	\right) \sim N_{p+1}\left( \left( 
	\begin{array}{c}
		0 \\ 
		\mathbf{0}_{p\times 1}%
	\end{array}%
	\right) ,\left( 
	\begin{array}{cc}
		1 & \mathbf{0}_{p\times 1}^{t} \\ 
		\mathbf{0}_{p\times 1} & I_{p\times p}%
	\end{array}%
	\right) \right)
\end{equation*}%
Then, it holds that 
\begin{equation*}
	\left( \boldsymbol{\lambda }^{t}\mathbf{X},y-\mathbf{\alpha }^{t}\mathbf{X}%
	\right) ^{t}\sim N_{2}\left( \left( 
	\begin{array}{c}
		0 \\ 
		0%
	\end{array}%
	\right) ,\left( 
	\begin{array}{cc}
		1 & -\boldsymbol{\alpha }^{t}\boldsymbol{\lambda } \\ 
		-\boldsymbol{\alpha }^{t}\boldsymbol{\lambda } & 1+\left \Vert \boldsymbol{\alpha }%
		\right \Vert ^{2}%
	\end{array}%
	\right) \right)
\end{equation*}%
for a unitary vector $\boldsymbol{\lambda }$. Therefore, 
\begin{equation*}
	P_{B}\left( \left( \boldsymbol{\lambda }^{t}\mathbf{X}\right) (y-\mathbf{X}^{t}%
	\boldsymbol{\alpha })\right) \geq 0)=1-E\Phi \left( \frac{\boldsymbol{\lambda }^{t}%
		\mathbf{XX}^{t}\boldsymbol{\alpha }}{\sigma \left \vert \mathbf{\lambda }^{t}%
		\mathbf{X}\right \vert }\right) .
\end{equation*}%
Since $h\left( -\rho \right) =1-h\left( \rho \right) $ and $\dot{h}:\left[
-1,1\right] \rightarrow \left[ 0,1\right] $ is increasing, we have that, for 
$b=\left \Vert \boldsymbol{\alpha } \right \Vert $ 
\begin{eqnarray*}
	P_{B}\left( \left( \boldsymbol{\lambda }^{t}\mathbf{X}\right) (y-\mathbf{X}^{t}%
	\boldsymbol{\alpha })\right) &\geq &0)=1-P_{B}\left( \left( \boldsymbol{\lambda }^{t}%
	\mathbf{X}\right) (y-\mathbf{X}^{t}\boldsymbol{\alpha })\right) <0) \\
	&=&1-\left[ 1-h\left( -\frac{\boldsymbol{\alpha }^{t}\boldsymbol{\lambda }}{\sqrt{%
			1+b^{2}}}\right) \right] \\
	&=&h\left( -\frac{\boldsymbol{\alpha }^{t}\boldsymbol{\lambda }}{\sqrt{1+b^{2}}}%
	\right) =1-h\left( \frac{\boldsymbol{\alpha }^{t}\boldsymbol{\lambda}}{\sqrt{1+b^{2}%
	}}\right)
\end{eqnarray*}%
This entails that, 
\begin{eqnarray*}
	\inf_{\left \Vert \boldsymbol{\lambda }\right \Vert =1}P_{B}\left( \left( \boldsymbol{%
		\lambda }^{t}\mathbf{X}\right) (y-\mathbf{X}^{t}\boldsymbol{\alpha })\right)
	&\geq &0)=1-\sup_{\left \Vert \lambda \right \Vert =1}E\Phi \left( \frac{%
		\boldsymbol{\lambda }^{t}\mathbf{XX}^{t}\boldsymbol{\alpha }}{\left \vert \boldsymbol{%
			\lambda }^{t}\mathbf{X}\right \vert }\right) \\
	&=&1-E\Phi \left( \left \Vert \boldsymbol{\alpha }\right \Vert \frac{\left \vert 
		\boldsymbol{\alpha }^{t}\mathbf{X}\right \vert }{\left \Vert \boldsymbol{\alpha }%
		\right \Vert }\right) =1-E\Phi \left( \left \Vert \boldsymbol{\alpha }\right \Vert
	\left \vert Z\right \vert \right) =1-g\left( \left \Vert \boldsymbol{\alpha }%
	\right \Vert \right) \\
	&=&1-\sup_{\left \Vert \lambda \right \Vert =1}h\left( \frac{\boldsymbol{\alpha }%
		^{t}\boldsymbol{\lambda }}{\sqrt{1+b^{2}}}\right) =1-h\left( \frac{b}{\sqrt{%
			1+b^{2}}}\right)
\end{eqnarray*}%
Then, 
\begin{equation*}
	h\left( \frac{b}{\sqrt{1+b^{2}}}\right) =g\left( b\right) =\frac{%
		1+\varepsilon }{2\left( 1-\varepsilon \right) },
\end{equation*}%
and the lemma follows. $\hfill \square $

\noindent \textbf{Proof of Remark 2.} The  implosion bias always rules the maxbias. To see this, 
observe that implosion and explosion bias
show up in the maximum bias: $f_{1}\left( \varepsilon \right) =\frac{1}{%
	\sqrt{\beta }}\Phi ^{-1}\left( \frac{3-\varepsilon }{4\left( 1-\varepsilon
	\right) }\right) $ and $f_{2}\left( \varepsilon \right) =\frac{\sqrt{\beta }%
}{\Phi ^{-1}\left( \frac{3-5\varepsilon }{4\left( 1-\varepsilon \right) }%
	\right) }$.\ We will see that  $f_{1}\left( \varepsilon \right) \leq $ $%
f_{2}\left( \varepsilon \right) $ for all $\varepsilon \in \left(
0,1/3\right)$. Observe that   
\[
\frac{3-\varepsilon }{4\left( 1-\varepsilon \right) }=\frac{3-5\varepsilon }{%
	4\left( 1-\varepsilon \right) }+\frac{\varepsilon }{\left( 1-\varepsilon
	\right) }.
\]%
$f_{1}\left( \varepsilon \right) \leq $ $f_{2}\left( \varepsilon \right) $
is equivalent to $\Phi ^{-1}\left( \frac{3-\varepsilon }{4\left(
	1-\varepsilon \right) }\right) \Phi ^{-1}\left( \frac{3-5\varepsilon }{%
	4\left( 1-\varepsilon \right) }\right) \leq \beta .$ Call $x=\frac{%
	3-5\varepsilon }{4\left( 1-\varepsilon \right) }$ which implies that $%
\varepsilon =\frac{4x-3}{4x-5}$ $\ $and $\frac{\varepsilon }{\left(
	1-\varepsilon \right) }=-2x+1.5.$ If $\varepsilon \in \left( 0,1/3\right) $
then $x\in \left( 0.5,0.75\right) .$ Then, let us consider the function $%
g\left( x\right) =\Phi ^{-1}\left( x\right) \Phi ^{-1}\left( 1.5-x\right)
=u\left( x\right) v\left( x\right) .$ We can show that $g$ is increasing in $%
\left( 0.5,0.75\right) $ and the maximum occurs at $0.75$, $g\left(
0.75\right) =\beta $ and the statement follows. Let us differentiate $g$ to
see whether the derivative is positive. Note that $u^{\prime }\left(
x\right) =\frac{1}{\varphi \left( u\left( x\right) \right) }$ and $v^{\prime
}\left( x\right) =-\frac{1}{\varphi \left( v\left( x\right) \right) }$  
\[
g^{\prime }\left( x\right) =\frac{v\left( x\right) }{\varphi \left( u\left(
	x\right) \right) }-\frac{u\left( x\right) }{\varphi \left( v\left( x\right)
	\right) }=\frac{v\left( x\right) \varphi \left( v\left( x\right) \right)
	-u\left( x\right) \varphi \left( u\left( x\right) \right) }{\varphi \left(
	u\left( x\right) \right) \varphi \left( v\left( x\right) \right) }
\]%
Then we need to show that $N\left( x\right) =v\left( x\right) \varphi \left(
v\left( x\right) \right) -u\left( x\right) \varphi \left( u\left( x\right)
\right) \geq 0.$ Note that $\varphi ^{\prime }\left( y\right) =-y\varphi
\left( y\right) .$ Then 
\begin{eqnarray*}
	&&N^{\prime }\left( x\right)  \\
	&=&v^{\prime }\left( x\right) \varphi \left( v\left( x\right) \right)
	+v\left( x\right) v^{\prime }\left( x\right) \varphi ^{\prime }\left(
	v\left( x\right) \right) -u^{\prime }\left( x\right) \varphi \left( u\left(
	x\right) \right) -u\left( x\right) u^{\prime }\left( x\right) \varphi
	^{\prime }\left( u\left( x\right) \right)  \\
	&=&-\frac{1}{\varphi \left( v\left( x\right) \right) }\varphi \left( v\left(
	x\right) \right) -\frac{1}{\varphi \left( u\left( x\right) \right) }\varphi
	\left( u\left( x\right) \right) +v^{2}\left( x\right) \frac{\varphi \left(
		v\left( x\right) \right) }{\varphi \left( v\left( x\right) \right) }%
	+u^{2}\left( x\right) \frac{\varphi \left( u\left( x\right) \right) }{%
		\varphi \left( u\left( x\right) \right) } \\
	&=&u^{2}\left( x\right) +v^{2}\left( x\right) -2.
\end{eqnarray*}%
To see the behavior of the function  $N^{\prime }\left( x\right) ,$ we
differentiate it,%
\begin{eqnarray*}
	N^{^{\prime \prime }}\left( x\right)  &=&2u\left( x\right) u^{\prime }\left(
	x\right) +2v\left( x\right) v^{\prime }\left( x\right)  \\
	&=&2\frac{u\left( x\right) }{\varphi \left( u\left( x\right) \right) }-2%
	\frac{v\left( x\right) }{\varphi \left( v\left( x\right) \right) }=2\frac{%
		u\left( x\right) \varphi \left( v\left( x\right) \right) -v\left( x\right)
		\varphi \left( u\left( x\right) \right) }{\varphi \left( u\left( x\right)
		\right) \varphi \left( v\left( x\right) \right) }.
\end{eqnarray*}%
Note that $u\left( x\right) =\Phi ^{-1}\left( x\right) <\Phi ^{-1}\left(
1.5-x\right) =v\left( x\right) $ and consequently $\varphi \left( v\left(
x\right) \right) <\varphi \left( u\left( x\right) \right) $ since $\varphi $
is decreasing. Then $N^{^{\prime \prime }}\left( x\right) <0$ and therefore $%
N^{\prime }$ is strictly decreasing. $lim_{x\rightarrow 0.5}N^{\prime
}\left( x\right) =\infty $ and $lim_{x\rightarrow 0.75}N^{\prime }\left(
x\right) =\beta -2<0.$ Therefore there exists just one point $x_{0}$ for
which  $N^{\prime }\left( x_{0}\right) =0,$ $N^{\prime }\left( x\right) >0$
if $x<x_{0}$ and negative otherwise$.$ Thus the function $N\left( x\right) $
is increasing  if $x<x_{0}$ and decreasing otherwise. Since $%
\lim_{x\rightarrow 0.5}N\left( x\right) =0$ and $lim_{x\rightarrow
	0.75}N\left( x\right) =0,$ this entails that $N\left( x\right) >0$ in $%
\left( 0.5,0.75\right) .$ Then $g^{\prime }\left( x\right) >0$ and $g\left(
x\right) $ is increasing in $\left( 0.5,0.75\right) $ which says that $%
g\left( x\right) \leq g\left( 0.75\right) =\beta $ and therefore 
\[
\frac{1}{\sqrt{\beta }}\Phi ^{-1}\left( \frac{3-\varepsilon }{4\left(
	1-\varepsilon \right) }\right) \leq \frac{\sqrt{\beta }}{\Phi ^{-1}\left( 
	\frac{3-5\varepsilon }{4\left( 1-\varepsilon \right) }\right) }.
\]%
Therefore, the implosion bias always rules the maxbias.

\subsection{Proofs in Section 6.}


\noindent \textbf{Proof of Lemma 11}. Let\textbf{\ }$\mathcal{D}%
_{MR}^{E}\left( B,P\right) {\ =}\inf_{U\mathbf{\in 
		\mathbb{R}
	}^{p\times m}}P\left( \left\Vert Y-B^{t}\mathbf{X}\right\Vert \leq
\left\Vert Y-U^{t}\mathbf{X}\right\Vert \right) .$ Thus, 
\begin{eqnarray*}
	P\left( \left\Vert Y-B^{t}\mathbf{X}\right\Vert \leq \left\Vert Y-U^{t} 
	\mathbf{X}\right\Vert \right) &=&P\left( \left\Vert Y-B^{t}\mathbf{X}
	\right\Vert ^{2}\leq \left\Vert Y-B^{t}\mathbf{X+}\left( B-U\right) ^{t} 
	\mathbf{X}\right\Vert ^{2}\right) \\
	&=&P\left( 0\leq -2\left\langle Y-B^{t}\mathbf{X,}\left( B-U\right) ^{t} 
	\mathbf{X}\right\rangle +\left\Vert \left( B-U\right) ^{t}\mathbf{X}
	\right\Vert ^{2}\right) .
\end{eqnarray*}%
Take $U=B-tV/2,$ $V$ any matrix in $\mathbf{\ 
	\mathbb{R}
}^{p\times m}$, $t\geq 0.$ Then, we get that 
\begin{eqnarray*}
	&&P\left( 0\leq -2\left\langle Y-B^{t}\mathbf{X,}\left( U-B\right) ^{t} 
	\mathbf{X}\right\rangle +\left\Vert \left( U-B\right) ^{t}\mathbf{X}
	\right\Vert ^{2}\right) \\
	&=&P\left( 0\leq t\left\langle Y-B^{t}\mathbf{X,}V^{t}\mathbf{X}
	\right\rangle +\frac{t^{2}}{4}\left\Vert V^{t}\mathbf{X}\right\Vert
	^{2}\right) \\
	&=&P\left( 0\leq t\left( \left\langle Y-B^{t}\mathbf{X,}V^{t}\mathbf{X}
	\right\rangle +\frac{t}{4}\left\Vert V^{t}\mathbf{X}\right\Vert ^{2}\right)
	\right) \\
	&=&P\left( 0\leq \left\langle Y-B^{t}\mathbf{X,}V^{t}\mathbf{X}\right\rangle
	+\frac{t}{4}\left\Vert V^{t}\mathbf{X}\right\Vert ^{2}\right) .
\end{eqnarray*}%
Since $\left[ 0\leq \left\langle Y-B^{t}\mathbf{X,}V^{t}\mathbf{X}%
\right\rangle +\frac{t}{4}\left\Vert V^{t}\mathbf{X}\right\Vert ^{2}\right]
\supseteq \left[ 0\leq \left\langle Y-B^{t}\mathbf{X,}V^{t}\mathbf{X}%
\right\rangle \right] $ for all $t\geq 0,$ we have that%
\begin{equation*}
	\inf_{U\mathbf{\in 
			\mathbb{R}
		}^{p\times m}}P\left( \left\Vert Y-B^{t}\mathbf{X}\right\Vert \leq
	\left\Vert Y-U^{t}\mathbf{X}\right\Vert \right) =\underset{V\mathbf{\in 
			\mathbb{R}
		}^{p\times m}}{\inf }\left[ 0\leq \left\langle Y-B^{t}\mathbf{X,}V^{t} 
	\mathbf{X}\right\rangle \right]
\end{equation*}%
and $\mathcal{D}_{MR}^{E}\left( B,P\right) =\mathcal{D}_{MR}\left(
B,P\right) $ as it was claimed.\ $\hfill \square $

\bigskip

\noindent \textbf{Proof of Lemma 12.} It is easily proved that $%
\min_{\lambda }P\left[ \left \vert y-\mu \right \vert \leq \left \vert
y-\lambda \right \vert \right] =\min \left \{ P(y\leq \mu ),P(y\geq \mu
)\right \} .$ Let us calculatete $\inf_{\mathit{\gamma} }P\left( \left[ \left \vert \left
\vert \frac{y-\mu }{\sigma }\right \vert -1\right \vert \leq \left \vert
\left \vert \frac{y-\mu }{\mathit{\gamma} }\right \vert -1\right
\vert \right]
\right)$. Call $z=\left \vert y-\mu \right \vert$. Then,

\begin{eqnarray*}
	\left \{ \left \vert \left \vert \frac{y-\mu }{\sigma }\right \vert -1\right
	\vert \leq \left \vert \left \vert \frac{y-\mu }{\mathit{\gamma} }\right \vert
	-1\right \vert \right \} &=&\left \{ \left \vert \left \vert \frac{y-\mu }{
		\sigma }\right \vert -1\right \vert ^{2}\leq \left \vert \left \vert \frac{
		y-\mu }{\mathit{\gamma} }\right \vert -1\right \vert ^{2}\right \} \\
	&=&\text{ }\left \{ \frac{z^{2}}{\sigma ^{2}}-2\frac{z}{\sigma }\leq \frac{
		z^{2}}{\mathit{\gamma} ^{2}}-2\frac{z}{\mathit{\gamma} }\right \} =\left \{ \frac{z}{\sigma
		^{2}}-\frac{2}{\sigma }\leq \frac{z}{\mathit{\gamma} ^{2}}-\frac{2}{\mathit{\gamma} }\right \}
	\\
	&=&\left \{ z\left( \frac{1}{\sigma ^{2}}-\frac{1}{\mathit{\Gamma} ^{2}}\right) \leq
	2\left( \frac{1}{\sigma }-\frac{1}{\mathit{\gamma} }\right) \right \} \\
	&=&\left \{ 
	\begin{array}{cc}
		z\leq \frac{2}{\frac{1}{\sigma }+\frac{1}{\mathit{\gamma} }} & \text{ si }\sigma
		<\mathit{\gamma} \\ 
		z\geq \frac{2}{\frac{1}{\sigma }+\frac{1}{\mathit{\gamma} }} & \text{si }\sigma
		>\mathit{\gamma}%
	\end{array}
	\right. .
\end{eqnarray*}%
Therefore, we get that%
\begin{eqnarray*}
	&&\min_{\mathit{\gamma} }P\left( \left \vert \frac{z}{\sigma }-1\right \vert ^{2}\leq
	\left \vert \frac{z}{\mathit{\gamma} }-1\right \vert ^{2}\right) \\
	&=&\min_{\mathit{\gamma} }\left \{ P\left( z\leq \frac{2}{\frac{1}{\sigma }+\frac{1}{
			\mathit{\gamma} }}\right) 1_{\left( \sigma ,\infty \right) }\left( \mathit{\gamma} \right)
	+P\left( z\geq \frac{2}{\frac{1}{\sigma }+\frac{1}{\mathit{\gamma} }}\right)
	1_{\left( 0,\sigma \right) }\left( \mathit{\gamma} \right) \right \} \\
	&=&\left \{ 
	\begin{array}{c}
		\min_{\text{ }\sigma <\mathit{\gamma} }P\left( z\leq \frac{2}{\frac{1}{\sigma }+\frac{
				1}{\mathit{\gamma} }}\right) =P\left( z/\sigma \leq 1\right) \\ 
		\min_{\text{ }\sigma >\mathit{\gamma} }P\left( z\geq \frac{2}{\frac{1}{\sigma }+\frac{
				1}{\mathit{\gamma} }}\right) =P\left( z/\sigma \geq 1\right)%
	\end{array}
	\right. \\
	&=&\min \left( P\left( z\leq \sigma \right) ,P\left( z\geq \sigma \right)
	\right) \\
	&=&\min \left( P\left( \mu -\sigma \leq y\leq \mu +\sigma \right) ,P\left( 
	\left[ y\geq \mu +\sigma \right] \cup \left[ y\leq \mu -\sigma \right]
	\right) \right)
\end{eqnarray*}%
and the depth turns out to be%
\begin{equation*}
	D_{LS}^{1}\left( \mu ,\sigma ,P\right) =\min \left \{ 
	\begin{array}{c}
		\min \left \{ P(y\leq \mu ),P(y\geq \mu )\right \} , \\ 
		\min \left( P\left( \mu -\sigma \leq y\leq \mu +\sigma \right) ,P\left( 
		\left[ y\geq \mu +\sigma \right] \cup \left[ y\leq \mu -\sigma \right]
		\right) \right)%
	\end{array}
	\right \} .
\end{equation*}%
If $0.5\leq D_{LS}^{1}\left( \mu ,\sigma ,P\right) $ we would say that $%
P(y\leq \mu )\geq 0.5\leq P(y\geq \mu )$ and $P\left( \left \vert y-\mu
\right \vert \leq \sigma \right) \geq 0.5\leq P\left( \left \vert y-\mu
\right \vert \geq \sigma \right) ,$which entails that $\hat{\mu}%
=med_{P}\left( Y\right) $ and $\hat{\sigma}=med_{P}\left( \left \vert y-\hat{
	\mu}\right \vert \right) .$\ $\hfill \square $

\bigskip

\section{Appendix B: Numerical study}
\subsection{The estimators}
Let us consider $\mathbf{X}=(X_{1},\ldots ,X_{p})^{T}$ a random vector and $%
\mathbf{X}_{1},\ldots ,\mathbf{X}_{n}$ a random sample from $\mathbf{X}$.
Let $\mathbb{X}=\{\mathbf{x}_{1},\ldots ,\mathbf{x}_{n}\}\subset \mathbb{R}%
^{p}$ denote a dataset. Consider $\mathbf{x},\boldsymbol{\mu }\in \mathbb{R}%
^{p}$ and $\mathit{\Sigma }\in \mathbb{R}^{p\times p}$, and let $d(\mathbf{x}%
,\boldsymbol{\mu },\mathit{\Sigma })=(\mathbf{x}-\boldsymbol{\mu })^{\prime }%
\mathit{\Sigma }^{-1}(\mathbf{x}-\boldsymbol{\mu })$ be the squared
Mahalanobis distance. Let $d_{i}=d(\mathbf{x}_{i},\boldsymbol{\mu },\mathit{%
	\Sigma })$ and $\mathbf{d}(\boldsymbol{\mu },\mathit{\Sigma })=(d_{1},\ldots
,d_{n}),$ $i=1,\dots ,n$.

We include several estimators in our simulation study. To compute them, we
use functions from R packages available on the Comprehensive R Archive
Network (CRAN), with default argument values in all cases. We selected
packages whose implementations ensure Fisher consistency of the estimators.
A brief description of them is as follows.

\begin{itemize}
	\item[\textbf{1.}] The sample covariance matrix (\textsc{SCOV)}, defined by $%
	S_{n}=\frac{1}{n}\sum_{i=1}^{n}(\boldsymbol{x}_{i}-\bar{\mathbf{x}}_{n})( 
	\mathbf{x}_{i}-\bar{\mathbf{x}}_{n})^{\prime }$, where $\bar{\mathbf{x}}_{n}$
	is the sample mean, included to provide a benchmark for comparison with the
	robust estimators that follow.
	
	\item[\textbf{2.}] The minimum volume ellipsoid estimator (\textsc{MVE}) was
	introduced by \cite{Rousseeuw1985}. Heuristically, it finds the ellipsoid
	with smallest volume that covers $h$ data points. More precisely, let $h$ be
	an integer such that $h\in \left[ \lfloor n/2\rfloor +1,n\right] $ and $c=%
	\sqrt{\chi _{p,\alpha }^{2}}$, with $\alpha =h/n$ where, given $\beta \in
	(0,1)$, $\chi _{p,\beta }^{2}$ stands for the $\beta $ quantile of the
	chi-square distribution with $p$ degrees of freedom. Let $\mathcal{E}$ as in
	(5) of Section 1 and set $C_{h,c}=\left\{ \boldsymbol{t}\in \mathbb{R}^{p},%
	\mathit{\Sigma }\in \mathcal{E}:\#\left\{ i:d\left( \mathbf{x}_{i},%
	\boldsymbol{\ \ t,}\mathit{\Sigma }\right) \leq c^{2}\right\} \geq h\right\}
	.$ Then, the \textsc{\ MVE} $(\hat{\boldsymbol{\mu }},\mathit{\Sigma })$ is
	defined to be 
	\begin{equation*}
		\left( \hat{\boldsymbol{\mu }},\mathit{\Sigma }\right) =\arg \underset{%
			\left( \boldsymbol{t},\mathit{\Sigma }\right) \in C_{h,c}}{\min }\det \left( 
		\mathit{\Sigma }\right) .
	\end{equation*}%
	Note that $h=[(n+p+1)/2]$ ensures the maximal breakdown point, as described
	in \cite{VanAelstRousseeuw2009}. 
	\textsc{MVE} is computed in \textsf{R} using the function \texttt{CovMve}
	from the package \textsc{rrcov} \citep{rrcov}.
	
	\item[\textbf{3.}] The minimum covariance determinant estimator (\textsc{MCD}%
	), was proposed by \cite{Rousseeuw1985}). Given the $n$ observations, this
	method finds $h$ observations $\mathbf{x}_{i_{1}},\ldots ,\mathbf{x}_{i_{h}}$
	with sample covariance matrix $S_{h}$ having the lowest determinant, with $%
	\lfloor (n+p+1)/2\rfloor \leq h\leq n.$ $S_{h}$ and the sample mean $\bar{ 
		\mathbf{x}}_{n}$ of the $h$ observations are the \textsc{MCD} estimators for
	multivariate scatter and location. It is computed in \textsf{R} using the 
	\texttt{CovMcd} function from the package \textsc{rrcov} developed by \cite%
	{TodorovFilzmoser2009}.
	
	\item[\textbf{4.}] S-estimators (\textsc{SE)} were introduced by \cite%
	{Davies1987}. Given the constant $\delta \in (0,1),$ an M-scale $S=S(\mathbf{%
		\ \ d}(\boldsymbol{\mu },\mathit{\Sigma }))$ is defined through the equation 
	\begin{equation*}
		\frac{1}{n}\sum_{i=1}^{n}\rho \!\left( \frac{d_{i}}{S}\right) =\delta ,
	\end{equation*}%
	with $\rho :\left[ 0,\infty \right) \rightarrow \left[ 0,1\right] $ a
	nondecreasing function, $\rho \left( 0\right) =0,$ $\sup_{x}\rho \left(
	x\right) =1,$ $\rho $ is continuous except at most for a finite set wherein
	it is right-continuous. Then the S-estimators for multivariate location and
	shape $(\widehat{\boldsymbol{\mu }},\widetilde{\mathit{\Sigma }})$ are
	defined by 
	\begin{equation*}
		(\widehat{\boldsymbol{\mu }},\widetilde{\mathit{\Sigma }})=\arg \min_{%
			\boldsymbol{\mu }\in \mathbb{R}^{p},\mathit{\Sigma }\in \mathcal{E},|\mathit{%
				\Sigma }|=1}S\bigl(\mathbf{d}(\boldsymbol{\mu },\mathit{\Sigma })\bigr),
	\end{equation*}%
	where $\delta $ controls the breakdown point, then it is taken as $\delta
	=1/2$. Hence, the S-estimator for multivariate scatter $\mathit{\Sigma }$ is
	defined as $\widehat{\mathit{\Sigma }}=S\,\widetilde{\mathit{\Sigma }}.$ 
	\textsc{SE} is computed using the function \texttt{CovSest} from the package 
	\textsc{rrcov}, choosing \texttt{method = "bisquare"}, with the bisquare
	function, $\rho (t)= 1-(1-t)^{3}  \text{if }t\leq 1$ and $1 1-(1-t)^{3}$ otherwise.
	
	\item[5.] S-estimators with non-monotonic weight functions (\textsc{Rocke} )
	were introduced by \cite{Rocke1996}. Note that S-estimators satisfy the
	system 
	\begin{eqnarray*}
		\frac{1}{n}\sum_{i=1}^{n}W\left( \frac{d_{i}}{S}\right) (\mathbf{x}_{i}-%
		\boldsymbol{\mu })(\mathbf{x}_{i}-\boldsymbol{\mu })^{\prime } &=&\mathit{%
			\Sigma }, \\
		\frac{1}{n}\sum_{i=1}^{n}W\left( \frac{d_{i}}{S}\right) (\mathbf{x}_{i}-%
		\boldsymbol{\mu }) &=&0,
	\end{eqnarray*}%
	where $W=\rho ^{\prime }$ is usually called a weight function. \textsc{Rocke}
	is computed using the function \texttt{CovSest} from the package \textsc{\
		rrcov} with option \texttt{method ="rocke"}, that by default sets $\alpha
	=0.1$. This routine considers a modification of Rocke's \textquotedblleft
	biflat\textquotedblright\ weight function given by $W(t)=\left[ 1-\left\{ (t-1)/\gamma\right\} ^{2}\right]$ if $1-\mathit{\gamma}<t<1+\mathit{\gamma}$,
	with tuning constant $\mathit{\gamma} =\min \!\left( 1,\frac{\chi _{p,1-\alpha }^{2}}{%
		p}-1\right) $ where $\alpha $ is chosen to control efficiency.
	
	\item[\textbf{6.}] MM-estimators (\textsc{MM}) were proposed by \cite{TT2000}
	. Let $\left(\widehat{\boldsymbol{\mu }}_{0},\widehat{\mathit{\Sigma }}%
	_{0}\right) $ be an initial high breakdown point estimator and the squared
	Mahalanobis distances $d_{i}^{0}=d\left( \mathbf{x}_{i},\widehat{\boldsymbol{%
			\ \mu }}_{0},\widehat{\mathit{\Sigma }}_{0}\right) $. Consider an M-scale $%
	S^{0}$ solving 
	\begin{equation*}
		\frac{1}{n}\sum_{i=1}^{n}\rho \left( \frac{d_{i}^{0}}{S}\right) =\delta .
	\end{equation*}%
	In the end, consider the location and shape estimator defined as 
	\begin{equation*}
		(\widehat{\boldsymbol{\mu}},\widetilde{\mathit{\Sigma }})=\arg \min_{u\in 
			\mathbb{R}^{p},\mathit{\Sigma }\in \mathcal{E},|\mathit{\Sigma }|=1}\sum_{i=1}^{n}\rho
		\left( \frac{d_{i}}{cS^{0}}\right) ,
	\end{equation*}%
	where $c$ is chosen to get $95\%$ asymptotic efficiency under normality. The
	MM-estimator of location and scatter is defined as $(\widehat{\boldsymbol{\
			\mu }},\widehat{\mathit{\Sigma }})$ where $\widehat{\mathit{\Sigma }}=S^{0}%
	\widetilde{\mathit{\Sigma }}$ . To compute \textsc{MM}, we use the function 
	\texttt{covRobMM} of the R package \textsc{RobStatTM}. It is initialized
	using a high-breakdown S-estimator of location and scatter (\texttt{KurtSDNew%
	}). If 	$s(d)=3.534-1.944\,d+0.864\,d^{2}-0.104\,d^{3}+0.004\,d^{4}$, then $\rho $ is the Smoothed Hard Rejection (SHR) $\rho $-function defined by 
	\begin{equation*}
		\rho _{\mathrm{SHR}}(d)=\frac{1}{6.494}%
		\begin{cases}
			d & \text{if }d\leq 4, \\ 
			s(d) & \text{if }4<d\leq 9, \\ 
			6.494 & \text{if }d>9.%
		\end{cases}%
	\end{equation*}%
	For more details, see \cite{MaronnaYohai2017}.
	
	\item[\textbf{7.}] The Stahel-Donoho location scatter estimator (\textsc{SD})
	was proposed by \cite{Stahel1981} and \cite{Donoho1982}, weighting the
	observations according to a notion of outlyingness. Given $\mathbf{u}\in 
	\mathcal{S}^{p-1}$, denote by $\mathbf{u}^{\prime }\mathbb{X}=\left\{ 
	\mathbf{u}^{\prime }\mathbf{x}_{1},\ldots ,\mathbf{u}^{\prime }\mathbf{x}%
	_{n}\right\} $. The outlyingness with respect to $\mathbb{X}$ of $\mathbf{x}%
	\in \mathbb{R}^{p}$ along $\mathbf{u}$ is defined by 
	$t(\mathbf{x},\mathbf{u})=\mathbf{x}^{\prime }\mathbf{u}-L\left( 
	\mathbf{u}^{\prime }X\right)/S\left( \mathbf{u}^{\prime }X\right)$,
	where $L$ and $S$ are robust univariate location and scale estimators, for
	example the median and the MAD respectively. Now, the outlyingness of $%
	\mathbf{x}$ is defined by 
	$t(\mathbf{x})=\sup_{\mathbf{u}\in \mathcal{S}^{p-1}}|t(\mathbf{x},\mathbf{u}%
	)|$.
	If $w_{ij}=W_{j}\left( t\left( \mathbf{x}_{i}\right) \right) ,j=1,2$ then
	the\ Stahel-Donoho location and scatter estimator is the pair $(\widehat{%
		\boldsymbol{\mu }},\widehat{\mathit{\Sigma }})$ defined by 
	\begin{eqnarray*}
		\widehat{\boldsymbol{\mu }} =\frac{1}{\sum_{i=1}^{n}w_{i1}}%
		\sum_{i=1}^{n}w_{i1}\mathbf{x}_{i}, \;\; 
		\widehat{\mathit{\Sigma }} =\frac{1}{\sum_{i=1}^{n}w_{i2}}%
		\sum_{i=1}^{n}w_{i2}\left( \mathbf{x}_{i}-\widehat{\boldsymbol{\mu }}\right)
		\left( \mathbf{x}_{i}-\widehat{\boldsymbol{\mu }}\right) ^{t}
	\end{eqnarray*}%
	It is computed using the function \texttt{CovSde} from the \textsc{rrcov}
	package. This function uses the piecewise polynomial weight $W_{\mathrm{opt}%
	} $ described in \cite{MaronnaYohai2017}.
	
	\item[\textbf{8.}] The Deepest Estimator (\textsc{MDepth}) proposed by \cite%
	{ChenGaoRen2018}. We use the function \texttt{matrix\_depth\_by\_descent}
	from the R package \textsc{DepthDescent}. For more details see \cite%
	{ChenGaoRen2018supp}.
\end{itemize}

\subsection{Complementary tables and figures}

\begin{figure}[!h]
	\centering 
	\includegraphics[width=\textwidth]{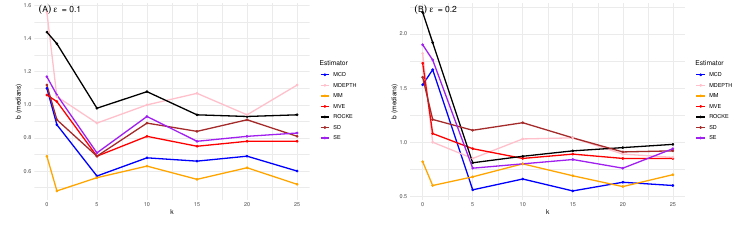}
	\caption{ $\hat{b}_{k}=\text{median}_{1\leq r\leq R}\{ \hat{b}_{k}^{(r)} \}$
		(log-scale) versus $k$ for each scatter estimator, under contamination
		levels $\protect\varepsilon$. Dimension $p=2$ and $n=20$. }
	
	\textbf{Alt text:} Two side-by-side line plots showing the median values of the bias as a function of $k$ for several scatter estimators on a log scale. Panel A corresponds to a lower contamination level and panel B to a higher contamination level.
	\label{b_K_p_2_n_20}
\end{figure}

\begin{figure}[!htbp]
	\centering
	
	\includegraphics[width=\textwidth]{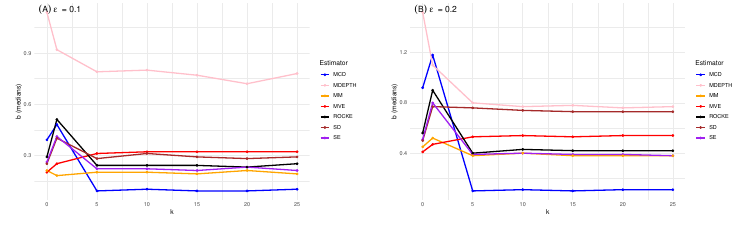}
	
	\caption{$\hat{b}_{k}=\mathrm{median}_{1\leq r\leq R}\{ \hat{b}_{k}^{(r)} \}$ (log-scale) versus $k$ for each scatter estimator, under contamination levels $\varepsilon$. Dimension $p=2$ and $n=1000$.}
	\textbf{Alt text:} Two side-by-side line plots showing the median values of the bias as a function of $k$ for several scatter estimators on a log scale. Panel A corresponds to a lower contamination level and panel B to a higher contamination level.
	\label{b_K_p_2_n_1000}
	
\end{figure}

In this subsection, we present additional results that complement those reported in Subsection 7.5 of the main paper. In Tables \ref{eff_n_50_cn_mean} and  \ref{eff_n_200_cn_mean}, the efficiencies of scatter estimators are reported for different combinations of $n$ and $p$, with respect to the CN measure.  Tables \ref{cn_p_2} to \ref{cn_p_15}  show the maximum medians of the condition number for each scatter estimator across sample sizes $n$ and contamination levels $\varepsilon $ for dimension $p=2,\ldots ,15$. 
Figures \ref{b_K_p_2_n_20} and \ref{b_K_p_2_n_1000} display the behavior of $%
\hat{b}_{k}=\text{median}_{1\leq r\leq R}\{ \hat{b}_{k}^{(r)}\}$, for each scatter estimator, as functions of $k$, for selected dimensions $p$ and sample sizes $n$, under both contamination
levels. 

\begin{table}[!h]
	\centering
	\begin{tabular}{r|rrrrrrrr}
		\hline
		& MVE & MCD & SE & ROCKE & MM & SD & MDEPTH &  \\ \hline
		$p=2$ & 0.508 & 0.688 & 0.55. & 0.359 & 0.917 & 0.589 & 0.647 &  \\ 
		$p=5$ & 0.559 & 0.736 & 0.826 & 0.565 & 0.941 & 0.601 & 0.798 &  \\ 
		$p=10$ & 0.558 & 0.700 & 0.910 & 0.524 & 0.958 & 0.585 & 0.774 &  \\ \hline
	\end{tabular}
	\caption{$\text{Eff}$ in log-scale based on $\hat{\text{B}}_{\text{CN}}$ (means) for $n=50$ over
		dimensions.}
	\label{eff_n_50_cn_mean}
\end{table}

\begin{table}[!ht]
	\centering
	\begin{tabular}{r|rrrrrrrr}
		\hline
		& MVE & MCD & SE & ROCKE & MM & SD & MDEPTH &  \\ \hline
		$p=2$ & 0.567 & 0.630 & 0.515 & 0.378 & 0.944 & 0.567 & 0.630 &  \\ 
		$p=5$ & 0.712 & 0.825 & 0.855 & 0.635 & 0.940 & 0.723 & 0.770 &  \\ 
		$p=10$ & 0.786 & 0.856 & 0.951 & 0.706 & 0.951 & 0.786 & 0.846 &  \\ 
		$p=15$ & 0.782 & 0.858 & 0.970 & 0.698 & 0.933 & 0.746 & 0.843 &  \\ \hline
	\end{tabular}
	\caption{$\text{Eff}$ in log-scale based on $\hat{\text{B}}_{\text{CN}}$ (means) for $n=200$ over
		dimensions.}
	\label{eff_n_200_cn_mean}
\end{table}

\begin{table}[!ht]
	\centering
	\begin{tabular}{ccrrrrrrrr}
		\hline
		$n$ & $\varepsilon$ & SCOV & MVE & MCD & SE & ROCKE & MM & SD & MDEPTH \\ 
		\hline
		20 & 0.10 & 4.72 & 1.45 & 1.13 & 1.58 & 1.93 & 0.69 & 1.24 & 0.97 \\ 
		20 & 0.20 & 5.36 & 2.12 & 1.95 & 2.35 & 2.38 & 0.81 & 1.69 & 1.32 \\ 
		80 & 0.10 & 4.81 & 0.81 & 0.60 & 0.79 & 0.94 & 0.40 & 0.72 & 0.62 \\ 
		80 & 0.20 & 5.54 & 1.15 & 1.26 & 1.29 & 1.40 & 0.46 & 1.14 & 1.05 \\ 
		1000 & 0.10 & 4.86 & 0.47 & 0.63 & 0.61 & 0.80 & 0.24 & 0.52 & 0.48 \\ 
		1000 & 0.20 & 5.55 & 0.79 & 1.33 & 1.17 & 1.30 & 0.45 & 1.02 & 0.96 \\ \hline
	\end{tabular}
	\caption{Empirical maximum bias $\hat{\text{B}}_{\text{CN}}$ (log scale) for
		each scatter estimator across sample sizes $n$ and contamination levels $
		\protect\varepsilon$ for dimension $p=2$.}
	\label{cn_p_2}
\end{table}

\begin{table}[!ht]
	\centering
	\begin{tabular}{llrrrrrrrr}
		\hline
		$n$ & $\varepsilon$ & SCOV & MVE & MCD & SE & ROCKE & MM & SD & MDEPTH \\ 
		\hline
		50 & 0.10 & 6.30 & 2.01 & 1.73 & 1.40 & 2.17 & 1.14 & 2.09 & 1.60 \\ 
		50 & 0.20 & 7.07 & 2.79 & 2.97 & 2.34 & 3.00 & 1.20 & 3.51 & 2.22 \\ 
		200 & 0.10 & 6.00 & 1.06 & 1.04 & 0.99 & 1.18 & 0.72 & 1.27 & 1.04 \\ 
		200 & 0.20 & 6.65 & 1.40 & 1.87 & 1.51 & 1.58 & 1.11 & 2.09 & 1.77 \\ 
		2500 & 0.10 & 5.81 & 0.70 & 0.66 & 0.76 & 0.75 & 0.52 & 0.82 & 0.92 \\ 
		2500 & 0.20 & 6.50 & 1.09 & 1.31 & 1.26 & 1.09 & 0.86 & 1.53 & 1.67 \\ \hline
	\end{tabular}
	\caption{Empirical maximum bias $\hat{\text{B}}_{\text{CN}}$ (log scale) for
		each scatter estimator under contamination levels $\protect\varepsilon$ and
		different sample sizes $n$ for dimension $p=5$.}
	\label{cn_p_5}
\end{table}

\begin{table}[!ht]
	\centering
	\begin{tabular}{ccrrrrrrrr}
		\hline
		$n$ & $\varepsilon$ & SCOV & MVE & MCD & SE & ROCKE & MM & SD & MDEPTH \\ 
		\hline
		100 & 0.10 & 7.00 & 2.30 & 2.06 & 1.64 & 1.92 & 1.40 & 2.52 & 2.52 \\ 
		100 & 0.20 & 7.80 & 3.39 & 3.35 & 2.24 & 2.48 & 1.58 & 4.41 & 3.97 \\ 
		400 & 0.10 & 6.70 & 1.20 & 1.29 & 1.30 & 1.14 & 1.09 & 1.51 & 2.03 \\ 
		400 & 0.20 & 7.41 & 1.74 & 2.33 & 1.86 & 1.39 & 1.56 & 2.64 & 3.59 \\ 
		5000 & 0.10 & 6.51 & 0.94 & 0.94 & 1.06 & 0.75 & 0.85 & 1.11 & 1.79 \\ 
		5000 & 0.20 & 7.22 & 1.46 & 1.79 & 1.64 & 0.97 & 1.34 & 2.06 & 3.34 \\ \hline
	\end{tabular}
	\caption{Empirical maximum bias $\hat{\text{B}}_{\text{CN}}$ (log scale) for
		each scatter estimator under different sample sizes $n$ and contamination
		levels $\protect\varepsilon$ for dimension $p=10$.}
	\label{cn_p_10}
\end{table}

\begin{table}[!ht]
	\centering
	\begin{tabular}{rr|rrrrrrrr}
		\hline
		$n$ & $\varepsilon$ & SCOV & MVE & MCD & SE & ROCKE & MM & SD & MDEPTH \\ 
		\hline
		150 & 0.10 & 7.48 & 2.48 & 2.41 & 1.92 & 1.95 & 1.64 & 3.03 & 3.60 \\ 
		150 & 0.20 & 8.22 & 3.36 & 3.87 & 2.67 & 2.43 & 1.93 & 5.08 & 5.66 \\ 
		600 & 0.10 & 7.15 & 1.83 & 2.41 & 1.92 & 1.95 & 1.37 & 3.03 & 3.60 \\ 
		600 & 0.20 & 7.88 & 2.04 & 2.76 & 2.16 & 1.29 & 1.87 & 3.11 & 5.19 \\ 
		7500 & 0.10 & 6.92 & 1.18 & 1.17 & 1.30 & 0.69 & 1.10 & 1.06 & 2.25 \\ 
		7500 & 0.20 & 7.62 & 1.75 & 2.15 & 1.91 & 0.84 & 1.65 & 1.73 & 4.77 \\ \hline
	\end{tabular}
	\caption{Empirical maximum bias $\hat{\text{B}}_{\text{CN}}$ (log scale) for
		each scatter estimator under different sample sizes $n$ and contamination
		levels $\protect\varepsilon$ for dimension $p=15$.}
	\label{cn_p_15}
\end{table}


\FloatBarrier

\end{document}